\documentclass[EJP,preprint]{ejpecp}
\usepackage[numbers,sort&compress]{natbib}
\usepackage[T1]{fontenc}
\usepackage[utf8]{inputenc}
\usepackage{enumitem}
\usepackage{xcolor}
\usepackage{pgfplots}
\usepackage{placeins}
\usepackage{needspace}
\usepgfplotslibrary{groupplots}
\usepackage{pgfmath-xfp}
\pgfplotsset{compat=1.18}
\definecolor{HSblue}{RGB}{0,114,178}
\definecolor{HSorange}{RGB}{213,94,0}
\definecolor{HSgreen}{RGB}{0,158,115}

\allowdisplaybreaks
\setlist[itemize]{leftmargin=1.2em}
\setlist[enumerate]{leftmargin=1.6em}

\hypersetup{
  pdfauthor={Qian Ai, Xiang Fang, Shengzhao Hou},
  pdftitle={Half-Line Hardy-Szego Point Processes: Fisher-Hartwig Asymptotics and Sparse Conditioning},
  pdfsubject={Strong count asymptotics, moving boundaries, and sparse conditioning for Hardy-Szego zeros},
  pdfkeywords={Hardy-Szego zeros, determinantal point processes, Fisher-Hartwig asymptotics, moving boundaries, large deviations, conditional Poisson limits}
}

\SHORTTITLE{Half-Line Hardy--Szeg\H{o} Point Processes}
\TITLE{Half-Line Hardy--Szeg\H{o} Point Processes: Fisher--Hartwig Asymptotics and Sparse Conditioning\support{X. F. was supported by National Science and Technology Council (No.~114-2115-M-A49-003-MY3). S. H. was supported by National Natural Science Foundation of China (No.~12371133).}}

\AUTHORS{%
  Qian~Ai\footnote{School of Mathematical Sciences, Soochow University, Suzhou 215006, P. R. China. \BEMAIL{20234007001@stu.suda.edu.cn}}
  \and
  Xiang~Fang\footnote{National Yang Ming Chiao Tung University, Hsinchu, Taiwan (R.O.C.). \BEMAIL{xfang@nycu.edu.tw}}
  \and
  Shengzhao~Hou\footnote{School of Mathematical Sciences, Soochow University, Suzhou 215006, P. R. China. \BEMAIL{shou@suda.edu.cn}}}

\AMSSUBJ{60G55; 60F10}
\AMSSUBJSECONDARY{60F05; 47B35}

\ABSTRACT{
A horizontal window extending from a fixed positive height to infinity
has infinite area but a finite expected number of Hardy--Szeg\H{o} zeros.
We study its zero count as the horizontal length $L$ grows.
Unlike for fixed bounded-height windows, the logarithmic moment
generating function has a new $\theta^2\log L/(2\pi^2)$ term;
we also identify its holomorphic order-one correction.
The logarithmic term has two equal contributions: a
Fisher--Hartwig jump in a canonical Wiener--Hopf model and
an endpoint effect in the original Fourier--Laplace operator.
This yields precise point and tail probabilities at positive
count densities.
We also determine the transition when a moving upper height
is comparable to $L$.
For sparse exact counts $n_L\to\infty$ with $n_L=o(L)$,
we obtain logarithmic probability asymptotics relative to the
exact void probability.
We determine the void exponent and sharp asymptotics relative
to the void probability for counts up to order $\log L$.
Under sparse exact-count conditioning, the empirical distribution
of rescaled horizontal positions and logarithmic heights
converges to the uniform law on the unit square.
Near distinct interior horizontal locations, the rescaled
configurations converge to independent Poisson processes
with uniform logarithmic-height marks.
}

\newcommand{\Logp}{\operatorname{Log}}
\newcommand{\Dres}{\Cc\setminus[0,1]}
\newcommand{\Qslit}{\Cc\setminus(-\infty,-1]}
\newcommand{\Hh}{\mathbb{H}}
\newcommand{\Rr}{\mathbb{R}}
\newcommand{\Cc}{\mathbb{C}}
\newcommand{\Zz}{\mathbb{Z}}
\newcommand{\Ee}{\mathbb{E}}
\newcommand{\Pp}{\mathbb{P}}
\newcommand{\ddA}{\mathrm{d}A}
\newcommand{\ii}{\mathrm{i}}
\newcommand{\Id}{\mathrm{Id}}
\DeclareMathOperator{\Tr}{Tr}
\DeclareMathOperator{\Li}{Li}
\newcommand{\dd}{\,\mathrm d}
\newcommand{\PPP}{\operatorname{PPP}}
\newcommand{\Law}{\mathcal L}
\newcommand{\Hol}{\operatorname{Hol}}
\newcommand{\Q}{\mathbb Q}

\begin{document}

\Needspace{7\baselineskip}
\section{Introduction and main results}\label{sec:intro}

A window extending to infinite height can have infinite area yet contain only finitely many Hardy--Szeg\H{o} zeros on average.
The Hardy--Szeg\H{o} Gaussian analytic function has a conformally invariant determinantal zero process, first studied in the disk by Peres and Vir\'ag~\cite{PeresVirag2005}.
In the upper half-plane this process is invariant under horizontal translations, making long horizontal windows a natural setting for studying counts.

We count the zeros with horizontal coordinate between $0$ and $L$ and height at least $a>0$.
Their intensity at height $y$ is $1/(4\pi y^2)$, so the expected count in this window is exactly $L/(4\pi a)$.
Thus the mean count grows linearly with $L$, although the window reaches arbitrarily far above the real axis.
We ask how the unbounded height range affects fluctuations and the probabilities of counts far from their mean.

Our earlier paper~\cite{AiFangHou2026} studied windows between two fixed positive heights.
It established large deviations and a strong Szeg\H{o} expansion for the logarithm of the moment generating function.
For each fixed upper height, this expansion consists of a leading term proportional to $L$ and a correction of order one.
Such a correction matters when one wants more than the leading exponential rate of a rare count.

Removing the upper boundary is not a routine limit of this formula.
The coefficient of $L$, which we call the pressure, is the limiting logarithmic moment generating function per unit horizontal length.
This coefficient has a finite limit as the upper height tends to infinity, but the order-one correction behaves differently.
Its second derivative at the origin grows like the logarithm of the upper height.
Thus the bounded-height expansion does not give a uniform order-one correction as the upper boundary recedes.
We must find the correct half-line expansion and understand what happens when the upper height grows together with $L$.

We show that the half-line expansion contains the logarithmic term $\theta^2\log L/(2\pi^2)$ together with an identified holomorphic order-one correction.
The logarithmic term contributes to the variance beyond its leading linear growth and produces a power of $L$ in precise probabilities for atypical counts of order $L$.
We also determine how the bounded-height and half-line expansions meet when the upper height moves with $L$.

A further question concerns the zeros that remain when the total count is forced far below its mean.
We condition on exactly $n_L$ zeros, where $n_L\to\infty$ and $n_L=o(L)$.
After rescaling horizontal positions and logarithmic heights, their empirical distribution converges to the uniform law on the unit square.
At the scale of their mean horizontal spacing, the local configurations near distinct interior points converge to independent Poisson processes, with heights uniform on a logarithmic scale.

\subsection{The model and the bounded-height comparison}
\label{subsec:model-scaling-intro}

We work in the upper half-plane $\Hh=\{z\in\Cc:\operatorname{Im}z>0\}$ and write $Z$ for the Hardy--Szeg\H{o} zero set.
With respect to planar Lebesgue measure $\ddA$, its correlation kernel is
\begin{equation}\label{eq:HS-kernel}
\mathcal K_{\Hh}(z,w)
=-\frac1{\pi(z-\overline w)^2},
\qquad z,w\in\Hh.
\end{equation}
For $a,L>0$, we count the zeros in the window above height $a$ by
\[
N_a^\infty(L)
:=\#\bigl(Z\cap([0,L]\times[a,\infty))\bigr),
\qquad N_\infty(L):=N_1^\infty(L).
\]
The diagonal of the kernel gives $\Ee N_a^\infty(L)=L/(4\pi a)$.

To study these counts, we use logarithmic moment generating functions normalized by the horizontal length.
For real $\theta$, set
\[
\Lambda_{a,L}^\infty(\theta)
:=\frac1L\log\Ee e^{\theta N_a^\infty(L)},
\qquad
\Lambda_{\infty,L}:=\Lambda_{1,L}^\infty.
\]
For the holomorphic strong asymptotics below, let $S_\pi:=\{\theta\in\Cc:|\operatorname{Im}\theta|<\pi\}$.
We write $\Logp$ for the principal logarithm.

An exact scaling relation lets us work at lower height one.
Indeed, the homogeneity of the kernel, $\mathcal K_{\Hh}(az,aw)=a^{-2}\mathcal K_{\Hh}(z,w)$, gives
\begin{equation}\label{eq:count-scaling-intro-halfline}
N_a^\infty(L)\stackrel d=N_\infty(L/a),
\qquad L>0,
\end{equation}
and hence
\[
\Lambda_{a,L}^\infty(\theta)
=\frac1a\Lambda_{\infty,L/a}(\theta),
\qquad\theta\in\Rr,\quad L>0.
\]
Thus it is enough to begin with lower height one, where the mean count per unit horizontal length is $\lambda_\infty:=1/(4\pi)$.

For comparison with the earlier paper, let $N_{[1,b]}(L):=\#(Z\cap([0,L]\times[1,b]))$ for $b>1$.
Theorem~1.4 of our earlier paper~\cite{AiFangHou2026} proves that, for each \emph{fixed} $b>1$,
\begin{equation}\label{eq:intro-recalled-interval-strong}
\log\Ee e^{\theta N_{[1,b]}(L)}
=L\Lambda_{[1,b]}(\theta)+\mathcal C_b(\theta)+o(1),
\qquad L\to\infty,
\end{equation}
locally uniformly for $\theta\in S_\pi$.
Here $\Lambda_{[1,b]}$ is the limiting scaled logarithmic moment generating function, and $\mathcal C_b$ is holomorphic on $S_\pi$.
For each fixed $b$, the correction $\mathcal C_b$ does not depend on $L$.
The theorem does not provide an error estimate uniform in $b$ when the upper height grows with $L$.

The earlier paper also identifies the large-height limit
\[
\Lambda_{[1,b]}(\theta)
\longrightarrow -\frac1{4\pi}\Li_2(1-e^\theta),
\qquad b\to\infty,
\]
locally uniformly on $S_\pi$, where $\Li_2$ denotes the principal dilogarithm on $\Cc\setminus[1,\infty)$.
The behavior of the strong correction is different: that paper gives
\[
\mathcal C_b''(0)
=\frac1{\pi^2}\log\frac{(b+1)^2}{4b}
\sim\frac{\log b}{\pi^2}.
\]
The second derivative at zero of a logarithmic moment generating function gives the variance, so $\mathcal C_b''(0)$ is the order-one correction to the variance in the fixed-height expansion.
This correction therefore grows without bound as the upper height increases.
Thus the bounded-height formula suggests the leading term of the half-line expansion, but it does not determine the correction that follows it.

\subsection{Half-line asymptotics and count probabilities}
\label{subsec:main-results-intro}

The leading pressure is the large-height limit just recalled.
Our first result establishes it directly for the count in the unbounded window, including its holomorphic extension.

\begin{theorem}[Thermodynamic limit]\label{thm:thermo-halfline}
For every \(\theta\in\Rr\), the limit
\[
\Lambda_\infty(\theta)
:=
\lim_{L\to\infty}\Lambda_{\infty,L}(\theta)
\]
exists, with convergence locally uniform on \(\Rr\), and is given by
\begin{equation}\label{eq:thermo-limit-halfline}
\Lambda_\infty(\theta)
=
\frac{1}{2\pi}
\int_0^\infty
\log\bigl(1+(e^\theta-1)e^{-2\xi}\bigr)\,d\xi
=
-\frac{1}{4\pi}\Li_2(1-e^\theta).
\end{equation}
Here \(\Li_2\) denotes the principal dilogarithm on \(\Cc\setminus[1,\infty)\).
Moreover, \(\Lambda_\infty\) extends holomorphically to \(S_\pi\).
At the origin,
\[
\Lambda_\infty'(0)=\frac{1}{4\pi},
\qquad
\Lambda_\infty''(0)=\frac{1}{8\pi}.
\]
\end{theorem}

The limiting pressure is strictly convex and its derivative ranges over $(0,\infty)$.
This gives the large deviation principle for the half-line count by the G\"artner--Ellis theorem.

\begin{theorem}[Large deviation principle]\label{thm:LDP-halfline}
The family $\{N_\infty(L)/L\}_{L>0}$ satisfies a large deviation principle on $[0,\infty)$ with speed $L$ and good rate function
\begin{equation}\label{eq:rate-halfline}
J_\infty(x)
:=\sup_{\theta\in\Rr}\bigl(\theta x-\Lambda_\infty(\theta)\bigr),
\qquad x\ge0.
\end{equation}
Moreover, $J_\infty$ is strictly convex on $(0,\infty)$ and vanishes only at $\lambda_\infty=1/(4\pi)$.
\end{theorem}

The large deviation principle determines probabilities on the exponential scale.
It does not determine the smaller factors in the probability of an individual atypical count.
This requires a strong expansion of the logarithmic moment generating function.
For $a,L>0$ we use the holomorphic logarithm of $\Ee e^{\theta N_a^\infty(L)}$ on $S_\pi$, normalized to vanish at $\theta=0$.
The following theorem is the main strong-asymptotic result of the paper.

\begin{theorem}[Strong half-line asymptotics]\label{thm:strong-halfline}
There exists a unique holomorphic function \(A_\infty\colon S_\pi\to\Cc\) such that, as \(L\to\infty\),
\begin{equation}\label{eq:main-strong-halfline}
\log\Ee\bigl[e^{\theta N_\infty(L)}\bigr]
=
L\Lambda_\infty(\theta)
+\frac{\theta^2}{2\pi^2}\log L
+A_\infty(\theta)
+o(1)
\end{equation}
locally uniformly for \(\theta\in S_\pi\).
The constant term decomposes as
\[
A_\infty(\theta)
=
A_{\rm FH}(\theta)+\Phi_\infty(\theta),
\]
where both summands are holomorphic on \(S_\pi\).
The canonical Fisher--Hartwig term \(A_{\rm FH}\) is given explicitly in Proposition~\ref{prop:explicit-AFH-halfline}.
The endpoint term \(\Phi_\infty\) is defined by \eqref{eq:fc-Phi-direct}, using the absolutely convergent scale integrals in \eqref{eq:fc-crossover-integrals}.
\end{theorem}

To explain the additional term, it is useful to compare the scalar profiles associated with the two height ranges:
\[
r_b(\xi)=e^{-2\xi}-e^{-2b\xi},
\qquad r_\infty(\xi)=e^{-2\xi},\qquad \xi>0.
\]
When these functions are extended by zero to negative frequencies, $r_b$ is continuous at the origin, whereas $r_\infty$ has a jump.
The corresponding canonical Wiener--Hopf comparison contributes $\theta^2\log L/(4\pi^2)$ to the half-line expansion.
But this is only half of the logarithmic correction.
The original Fourier--Laplace operator contains the additional factor
\[
\frac{2\sqrt{\xi\eta}}{\xi+\eta},
\]
whose limit at the frequency corner $(0,0)$ depends on the direction of approach.
Comparing the two operators produces another $\theta^2\log L/(4\pi^2)$.
In the earlier bounded-height analysis this comparison contributed at order one; the half-line requires both the new logarithmic term and a convergent finite correction.
The functions $A_{\rm FH}$ and $\Phi_\infty$ identify the two parts of that finite constant.

The change is visible in the variance.
A direct second-trace calculation gives
\[
\operatorname{Var}N_\infty(L)
=\frac{L}{8\pi}+\frac{\log L}{\pi^2}
 -\frac{\log2}{\pi^2}+o(1).
\]
The linear coefficient belongs to the Gaussian-scale behavior; the logarithmic and constant terms provide finer information.
Proposition~\ref{prop:exact-variance-halfline} gives an exact finite-length formula and a monotone variance-to-mean ratio.
Corollary~\ref{cor:even-cumulants-halfline} also shows that every fixed even cumulant of order at least four converges to a finite constant.
At positive atypical densities the same logarithmic correction yields a power of $L$, and $A_\infty$ determines the remaining multiplicative constant.

The strong expansion also yields the following precise asymptotics for point and tail probabilities.

\begin{theorem}[Precise positive-density deviations]\label{thm:precise-halfline}
For integers \(k\ge1\), put \(x_{L,k}:=k/L\) and let \(\theta_{L,k}\) be the unique real solution of \(\Lambda_\infty'(\theta_{L,k})=x_{L,k}\).
For every \(K\Subset(0,\infty)\), uniformly over \(k\) with \(x_{L,k}\in K\) as \(L\to\infty\),
\begin{equation}\label{eq:precise-point-halfline}
\Pp\{N_\infty(L)=k\}
=\frac{e^{A_\infty(\theta_{L,k})}L^{\theta_{L,k}^2/(2\pi^2)}}
{\sqrt{2\pi L\Lambda_\infty''(\theta_{L,k})}}
e^{-LJ_\infty(x_{L,k})}(1+o_K(1)).
\end{equation}
If \(K_+\Subset(\lambda_\infty,\infty)\) and \(K_-\Subset(0,\lambda_\infty)\), then uniformly for \(x_{L,k}\in K_+\) and \(x_{L,k}\in K_-\), respectively,
\begin{align}
\Pp\{N_\infty(L)\ge k\}
&=\frac{\Pp\{N_\infty(L)=k\}}{1-e^{-\theta_{L,k}}}(1+o_{K_+}(1)),\label{eq:precise-upper-halfline}\\
\Pp\{N_\infty(L)\le k\}
&=\frac{\Pp\{N_\infty(L)=k\}}{1-e^{\theta_{L,k}}}(1+o_{K_-}(1)).\label{eq:precise-lower-halfline}
\end{align}
\end{theorem}

The one-sided tail equivalents apply only when the density remains strictly on one side of the mean.
The point-probability equivalent is not asserted in the regime $k/L\to0$; that regime requires different arguments and is treated later in the paper.

The reference lower height one is not restrictive.
By the exact scaling relation \eqref{eq:count-scaling-intro-halfline}, the pressure, rate function and strong correction for a fixed lower height $a>0$ follow from the results above.
Write
\[
\Lambda_a^\infty(\theta)
:=\lim_{L\to\infty}\Lambda_{a,L}^\infty(\theta),
\qquad
J_a^\infty(x):=\sup_{\theta\in\Rr}
\bigl(\theta x-\Lambda_a^\infty(\theta)\bigr),
\]
and denote the corresponding strong constant by $A_{\infty,a}$.

\begin{corollary}\label{cor:general-halfline-scaling}
Let \(a>0\).
Then the corresponding limiting objects satisfy
\[
\Lambda_a^\infty(\theta)
=
a^{-1}\Lambda_\infty(\theta),
\qquad
J_a^\infty(x)
=
a^{-1}J_\infty(ax),
\qquad
A_{\infty,a}(\theta)
=
A_\infty(\theta)
-\frac{\theta^2}{2\pi^2}\log a.
\]
Moreover, the variables \(N_a^\infty(L)/L\) satisfy an LDP on \([0,\infty)\) with speed \(L\) and good rate function \(J_a^\infty\), and
\begin{equation}\label{eq:halfline-scaling-cor-strong}
\log\Ee\bigl[e^{\theta N_a^\infty(L)}\bigr]
=
L\Lambda_a^\infty(\theta)
+\frac{\theta^2}{2\pi^2}\log L
+A_{\infty,a}(\theta)
+o(1)
\end{equation}
locally uniformly for \(\theta\in S_\pi\).
\end{corollary}

\subsection{Moving upper cutoffs and the passage to the half-line}
\label{subsec:moving-intro}

The interval and half-line expansions suggest a further question.
What happens when the upper height $b=b_L$ tends to infinity together with $L$?
The expected number of zeros above height $b$ in a window of length $L$ is $L/(4\pi b)$.
Thus $b$ comparable with $L$ is a natural transition scale: at this scale the region omitted by the finite-height window has an expected count of order one.
The strong expansion detects this transition, including its effect on the constant term.

Recall the bounded-height count $N_{[1,b]}(L)$ defined above.
Its limiting pressure is
\[
\Lambda_{[1,b]}(\theta)
:=\frac1{2\pi}\int_0^\infty
\Logp\!\left(1+(e^\theta-1)(e^{-2\xi}-e^{-2b\xi})\right)\dd\xi,
\qquad \theta\in S_\pi,
\]
where the logarithm is normalized at $\theta=0$.
The next theorem distinguishes three regimes according to the ratio $b_L/L$.

\begin{theorem}[Moving-cutoff crossover]\label{thm:intro-moving-cutoff}
There exist $A_{\mathrm{cross}}(\,\cdot\,,c)\in\Hol(S_\pi)$,
$c>0$, and $H_*\in\Hol(S_\pi)$ such that, for every sequence
$b_L\to\infty$, the following asymptotics hold locally uniformly
on $S_\pi$ as $L\to\infty$.
\begin{enumerate}[label=\textup{(\roman*)}]
\item \emph{Subcritical regime.}
If $b_L/L\to0$, then
\[
\log\Ee e^{\theta N_{[1,b_L]}(L)}
=L\Lambda_{[1,b_L]}(\theta)
+\frac{\theta^2}{2\pi^2}\log(b_L/2)
+A_\infty(\theta)+H_*(\theta)+o(1).
\]

\item \emph{Critical regime.}
If $b_L/L\to c\in(0,\infty)$, then
\[
\log\Ee e^{\theta N_{[1,b_L]}(L)}
=L\Lambda_{[1,b_L]}(\theta)
+\frac{\theta^2}{2\pi^2}\log L
+A_{\mathrm{cross}}(\theta,c)+o(1).
\]

\item \emph{Supercritical regime.}
If $b_L/L\to\infty$, then
\[
\log\Ee e^{\theta N_{[1,b_L]}(L)}
=\log\Ee e^{\theta N_\infty(L)}+o(1)
=L\Lambda_\infty(\theta)
+\frac{\theta^2}{2\pi^2}\log L
+A_\infty(\theta)+o(1).
\]
\end{enumerate}
Moreover, for every $K\Subset S_\pi$, uniformly for $\theta\in K$,
$A_{\mathrm{cross}}(\theta,c)=A_\infty(\theta)+O_K(c^{-1})$
as $c\to\infty$, and
\[
A_{\mathrm{cross}}(\theta,c)
-\frac{\theta^2}{2\pi^2}\log(c/2)
\longrightarrow A_\infty(\theta)+H_*(\theta),
\qquad c\downarrow0.
\]
\end{theorem}

The function $A_{\mathrm{cross}}$ is given by \eqref{MC:eq:Across}.
Theorem~\ref{BC:thm:local} gives the local series \eqref{BC:eq:main} for $H_*$, with $H_*(0)=H_*'(0)=H_*''(0)=0$, and a power bound for the small-$c$ limit on compact subsets of $S_\pi$.
Lemma~\ref{SB:lem:boundary-tail} gives the uniform remainder $O_K(b/L)$ in \eqref{eq:intro-recalled-interval-strong} for $K\Subset S_\pi$ and $L\ge b\ge2$. In the subcritical and critical regimes, the leading pressure must be evaluated at the actual height $b_L$.
For example, if $b_L=L^\alpha$ with fixed $\alpha>0$, the theorem and \eqref{eq:pressure-height} give
\[
\frac{\log\Ee e^{\theta N_{[1,L^\alpha]}(L)}
-L\Lambda_{[1,L^\alpha]}(\theta)}{\log L}
\longrightarrow\frac{\theta^2}{2\pi^2}\min\{\alpha,1\}
\]
locally uniformly on $S_\pi$. Thus the coefficient of $\theta^2\log L$ increases linearly with $\alpha$ up to $\alpha=1$ and then remains constant.

There is also a direct relation with the correction in the earlier bounded-height theorem~\eqref{eq:intro-recalled-interval-strong}.
The large-height limit proved in Section~\ref{sec:cutoffs} (Theorem~\ref{SB:thm:static-bridge}) gives
\begin{equation}\label{eq:intro-static-largeheight}
\mathcal C_b(\theta)-\frac{\theta^2}{2\pi^2}\log(b/2)
\longrightarrow A_\infty(\theta)+H_*(\theta),
\qquad b\to\infty,
\end{equation}
locally uniformly on $S_\pi$.
In particular, the bounded-height constants have a well-defined limit after subtraction of their logarithmic divergence, but this limit includes $H_*$.
Formula~\eqref{eq:intro-static-largeheight} makes the connection between the two papers precise.
The corresponding statements for other fixed lower heights follow by dilation.

\subsection{Sparse counts and conditioned geometry}
\label{subsec:sparse-intro}

The precise point-probability formula above applies when the target density stays bounded away from zero.
A different question arises when the total number of zeros is forced to be much smaller than its usual order $L$.
Fix $a>0$, let $D_{a,L}=[0,L]\times[a,\infty)$, and write $N_{a,L}=N_a^\infty(L)$.
We consider an exact count $n_L\to\infty$ with $n_L=o(L)$.
Besides the cost of this rare event, we ask where the remaining zeros lie and what their local arrangement looks like.
Because there is no upper height bound, the conditioning can alter their vertical distribution on a scale not present in a fixed strip.

The probability of such a count can be compared with the exact void probability at logarithmic precision.

\begin{theorem}[Dilute relative logarithmic asymptotics]\label{thm:intro-dilute-relative}
Fix $a>0$.
Let $n_L$ be positive integers such that
\[
n_L\longrightarrow\infty,
\qquad
\frac{n_L}{L}\longrightarrow0.
\]
Then
\begin{equation}\label{eq:intro-dilute-relative}
\log\frac{\Pp\{N_{a,L}=n_L\}}{\Pp\{N_{a,L}=0\}}
=
L\left[J_a^\infty(0)-J_a^\infty(n_L/L)\right]+o(n_L).
\end{equation}
Equivalently,
\begin{equation}\label{eq:intro-dilute-relative-explicit}
\log\frac{\Pp\{N_{a,L}=n_L\}}{\Pp\{N_{a,L}=0\}}
=
n_L\left[
\log\frac{L}{n_L}
+\log\log\frac{L}{n_L}
-\log(4\pi a)+1+o(1)
\right].
\end{equation}
\end{theorem}

The error in this theorem is $o(n_L)$ in the logarithm of the \emph{ratio} to the exact void probability.
It is not a multiplicative asymptotic formula for the absolute point probability.
The result covers every diverging sublinear sequence $n_L$, even if it tends to infinity arbitrarily slowly; fixed positive counts are treated separately below.

There is a corresponding change in the geometry under exact conditioning.
The relevant horizontal scale is the conditioned mean spacing, while the vertical scale is logarithmic in the inverse count density.
The following theorem describes the macroscopic distribution of the remaining zeros and their local configuration near finitely many distinct interior points of the horizontal window.
We write $\PPP(\nu)$ for the law of a Poisson point process with intensity measure $\nu$.

\begin{theorem}[Sparse count-conditioned geometry]\label{thm:intro-conditioned-poisson}
Fix $a>0$.
Let $n_L$ satisfy the assumptions of Theorem~\ref{thm:intro-dilute-relative}, and set $d_L=n_L/L$ and $h_L=\log(L/n_L)$.
On the event $\{N_{a,L}=n_L\}$, define the empirical probability measure
\[
\mathfrak M_L
:=
\frac1{n_L}
\sum_{x+\ii y\in Z\cap D_{a,L}}
\delta_{\left(x/L,\,\log(y/a)/h_L\right)}.
\]
Under $\Pp(\,\cdot\mid N_{a,L}=n_L)$, $\mathfrak M_L$ converges in probability, for the weak topology on probability measures on $[0,1]\times[0,\infty)$, to $\mathbf1_{[0,1]^2}(u,v)\,\dd u\,\dd v$.
For the local limit, fix pairwise distinct $\tau_1,\ldots,\tau_r\in(0,1)$, with $r$ fixed, and define point measures on $E=\Rr\times(0,1)$ by
\[
\eta_{L,j}
=\left.
\sum_{x+\ii y\in Z\cap D_{a,L}}
\delta_{\left(d_L(x-\tau_jL),\,\log(y/a)/h_L\right)}
\right|_E.
\]
Then, under the same conditional law,
\begin{equation}\label{eq:intro-conditioned-poisson}
(\eta_{L,1},\ldots,\eta_{L,r})
\Longrightarrow
\bigotimes_{j=1}^r\PPP(\dd u\,\dd v)
\end{equation}
in the product vague topology.
Moreover, for every fixed finite family of relatively compact continuity sets, all mixed factorial moments of the corresponding counts converge to those of the limiting independent Poisson processes.
\end{theorem}

Thus, under the conditioning $N_{a,L}=n_L$, the empirical distribution of the rescaled horizontal positions and logarithmic heights converges in probability to the uniform law on the unit square.
At the horizontal mean-spacing scale $L/n_L$, local configurations near distinct bulk centers become independent Poisson processes with uniform height marks.
These conclusions concern the \emph{conditional} process, not merely a stationary process under an exponential tilt.
The proof first uses weighted Bergman kernels to identify a stationary comparison process at low density.
We then choose the exponential tilt so that the mean count in the finite window is exactly $n_L$.
Reduced-Palm estimates allow us to transfer the local limits to the law conditioned on that count.

The zero-count endpoint exhibits different subleading behavior.
Theorem~\ref{VC:thm:void} gives the void exponent
\[
\log\Pp\{N_a^\infty(L)=0\}
=-\frac{\pi}{24a}L+o(L).
\]
The remainder after subtracting $-\pi L/(24a)$ is nonnegative and $o(L)$, but its ratio to $\log(L/a)$ tends to infinity.
We do not identify its first subleading equivalent.

For very small counts a spectral quantity $W_a(L)$, defined in \eqref{VC:eq:W-def}, satisfies
\[
W_a(L)\sim \frac{L}{4\pi a}\log(L/a).
\]
For every fixed $M>0$,
\[
\Pp\{N_{a,L}=k\}
=\Pp\{N_{a,L}=0\}\frac{W_a(L)^k}{k!}(1+o(1))
\]
uniformly for integers $0\le k\le M\log(L/a)$.
For the growing-count range treated in the paper, the exact quantity $W_a(L)$ must be kept inside the power; its leading equivalent alone is not sufficient.
The small-count asymptotics and the near-zero generating-function limit are proved in Theorem~\ref{VC:thm:near-zero}.

\subsection{Earlier work and mathematical context}

The closest predecessor is our bounded-height paper~\cite{AiFangHou2026}, whose strong expansion is recalled in \eqref{eq:intro-recalled-interval-strong}.
Its large-height pressure expansion is recalled in Lemma~\ref{lem:pressure-height}.
The Gaussian-scale results of Ai, Guo and Zhou~\cite{AiGuoZhou2026} already include half-line heights.
For fixed bounded-height intervals, Zhou~\cite{Zhou2026} obtains precise count probabilities by the exponential-tilting and local-limit method used in Section~\ref{sec:entropy-precise-halfline}.
Chen~\cite{Chen2026} obtains absolute point-probability asymptotics, including fixed counts and vanishing densities, in the same height regime.
Our dilute formula concerns the half-line and gives logarithmic precision relative to the exact void probability.

For fixed bounded-height intervals, Ai~\cite{AiPalm2026} proves joint local marked Poisson limits under sparse exact-count conditioning.
Our half-line result also describes logarithmic height escape and the uniform empirical profile of Theorem~\ref{thm:intro-conditioned-poisson}.
For a fixed bounded-height interval the Fourier symbol is separated from one, whereas the half-line symbol approaches one at zero frequency.
The loss of this gap requires additional estimates in the sparse regime.
The height scale and limiting profile are determined by the weighted Bergman comparison and localization estimates in Section~\ref{sec:conditioning}.

Fenzl and Lambert~\cite{FenzlLambert2022} obtain precise deviations for disk counts of the same zero process.
Comparing their formulas with \eqref{eq:thermo-limit-halfline} shows that the uncentered logarithmic moment generating functions have the same limit after division by their respective mean counts.
The logarithmic correction and its finite part for the half-line windows require the boundary analysis developed here.

The analytic comparison belongs to the theory of Szeg\H{o}, Wiener--Hopf and Fisher--Hartwig determinants; see, for example,
\cite{Szego1952,Widom1982,BasorWidom1983,
BottcherSilbermannWidom1994a,BottcherSilbermannWidom1994b,
BasorWidom2004,BottcherWidom2005,Kozlowski2022}.
Discontinuous symbols can produce logarithmic boundary corrections \cite{Widom1990}.
Sobolev~\cite{Sobolev2017} studies two-parameter transitions from smooth to discontinuous Wiener--Hopf symbols.
A related transition from strong Szeg\H{o} to Fisher--Hartwig asymptotics for Toeplitz determinants is studied in \cite{ClaeysItsKrasovsky2011}.
In our problem, the classical scalar theory describes the jump contribution to the canonical Wiener--Hopf comparison.
For the original Fourier--Laplace operator, the additional logarithmic term and its finite part require the endpoint comparison developed below.
Section~\ref{sec:cutoffs} determines the order-one terms in the moving-height transition.

For general background on Gaussian analytic zeros and determinantal
point processes, see~\cite{Hough2009,Sodin2000,Soshnikov2000,
ShiraiTakahashi2003}.  We use the determinantal Bernoulli
representation for count probabilities, the G\"artner--Ellis theorem for the leading large deviations, and local count and Palm methods for conditioning; relevant references include
\cite{HKPV2006,Gartner1977,Ellis1984,ForresterLebowitz2014,
BarthelmeAmblardTremblay2019}.

\subsection{Proof strategy and organization}\label{subsec:organization}

The planar determinantal formula is first reduced to a one-dimensional Fredholm determinant, as in the bounded-height analysis.
We then compare its operator with a canonical Wiener--Hopf model.
The canonical determinant has a Fisher--Hartwig expansion, while the difference between the logarithms of the two determinants contributes an additional $\log L$ term.
The main task is not merely to bound this difference, but to identify its finite part after subtracting the logarithmic term.
To do this we analyze the Fourier--Laplace endpoint on dyadic frequency scales.
Local trace estimates, resolvent comparisons and a continuous-scale identity yield a convergent expression for $\Phi_\infty$.
Together with the canonical result, this proves Theorem~\ref{thm:strong-halfline} and identifies its holomorphic constant term.
The endpoint part has a power remainder, while the canonical expansion supplies an $o(1)$ remainder.

The further results use additional arguments.
For moving upper cutoffs, relative Fredholm determinants describe the critical regime.
A uniform boundary comparison controls the fixed-height approximation, while trace-coefficient estimates and holomorphic bounds identify the small-ratio correction throughout $S_\pi$.
For sparse exact-count conditioning, weighted Bergman comparison identifies the height distribution and the low-density Poisson limit; local kernel estimates and a reduced-Palm comparison then handle the conditioning.
The void exponent and small-count behavior follow from separate spectral and counting estimates.
In particular, the statements about exact sparse conditioning do not follow from applying the fixed-$\theta$ strong expansion with a parameter tending to minus infinity.

Section~\ref{sec:fredholm-pressure} establishes the pressure and large deviations.
Sections~\ref{sec:WH-halfline}--\ref{sec:finite-constant} prove the two parts of the strong expansion; Sections~\ref{sec:all-order}--\ref{sec:entropy-precise-halfline} develop renormalized trace coefficients and precise count probabilities.
Section~\ref{sec:cutoffs} treats moving upper heights, Section~\ref{sec:conditioning} sparse conditional geometry, and Section~\ref{sec:endpoint} voids, dilute relative count asymptotics, and very small counts.

\section{Fredholm reduction and large deviations}\label{sec:fredholm-pressure}

By scaling, it is enough to work with heights in $[1,\infty)$.
We establish trace class directly for the unbounded height interval and use the Fourier--Laplace reduction of \cite[Section~2]{AiFangHou2026}.
A trace-norm comparison with bounded-height windows then gives the pressure.
The value of the limiting symbol at the frequency endpoint is $r_\infty(0+)=1$; we will treat the corresponding zero-density endpoint of the rate function explicitly.

\subsection{Fredholm determinant and Fourier--Laplace reduction}\label{subsec:reduction-halfline}

For $L>0$ write $B_L^\infty=[0,L]\times[1,\infty)\subset\Hh$.

\begin{lemma}\label{lem:fredholm-rectangles-halfline}
Let \(L>0\), and let \(K_L^\infty\) be the integral operator on \(L^2(B_L^\infty,\ddA)\) with kernel \(\mathcal K_{\Hh}|_{B_L^\infty\times B_L^\infty}\).
Then \(K_L^\infty\) is a positive trace-class contraction, and
\begin{equation}\label{eq:fredholm-halfline-rectangle}
\Ee\bigl[e^{\theta N_\infty(L)}\bigr]
=
\det\!\bigl(\Id+(e^\theta-1)K_L^\infty\bigr),
\qquad \theta\in\Cc.
\end{equation}
Moreover, \(\Tr K_L^\infty=\Ee N_\infty(L)=L/(4\pi)\).
\end{lemma}

\begin{proof}
For compactly supported $g\in L^2(0,\infty)$, define
\[
(\mathsf Jg)(z)=\frac1{\sqrt\pi}
\int_0^\infty\sqrt\xi\,e^{\ii\xi z}g(\xi)\,d\xi.
\]
Horizontal Plancherel and $2\int_0^\infty\xi e^{-2\xi y}\,dy=1$ show that $\mathsf J$ extends to an isometry from $L^2(0,\infty)$ to $L^2(\Hh,\ddA)$.
Moreover,
\[
(\mathsf J\mathsf J^*)(z,w)
=\frac1\pi\int_0^\infty\xi e^{\ii\xi(z-\overline w)}\,d\xi
=-\frac1{\pi(z-\overline w)^2}
=\mathcal K_{\Hh}(z,w).
\]
Put $\mathsf J_L=\mathbf1_{B_L^\infty}\mathsf J$.
Then $K_L^\infty=\mathsf J_L\mathsf J_L^*$ and
\[
\|\mathsf J_L\|_{\mathrm{HS}}^2
=\int_{B_L^\infty}\mathcal K_{\Hh}(z,z)\,\ddA(z)
=\frac{L}{4\pi}.
\]
Thus $K_L^\infty$ is a positive trace-class contraction, with the claimed trace and mean.

Its eigenvalues $(\lambda_j)$ lie in $[0,1]$ and sum to $L/(4\pi)$.
The Bernoulli decomposition for a trace-class determinantal process \cite[Theorem~7]{HKPV2006} gives $N_\infty(L)\stackrel d=\sum_jX_j$, with independent $X_j\sim\operatorname{Bernoulli}(\lambda_j)$.
For $t>0$,
\[
\Ee e^{tN_\infty(L)}
=\prod_j\bigl(1+(e^t-1)\lambda_j\bigr)
\le\exp\!\left(\frac{(e^t-1)L}{4\pi}\right).
\]
All exponential moments are therefore finite.
The Bernoulli product extends to every complex $\theta$, and the trace-class Fredholm product gives \eqref{eq:fredholm-halfline-rectangle}.
\end{proof}

The same factorization gives a one-dimensional determinant.

\begin{lemma}\label{lem:reduction-halfline}
For \(L>0\), set \(X_L(s):=\int_0^L e^{\ii sx}\,dx\), and let \(T_L^\infty\) be the integral operator on \(L^2(0,\infty)\) with kernel
\[
\mathcal T_L^\infty(\xi,\eta)
=
\frac{\sqrt{\xi\eta}\,e^{-(\xi+\eta)}}
     {\pi(\xi+\eta)}
X_L(\xi-\eta).
\]
Then \(T_L^\infty\) is a positive trace-class contraction and has the same nonzero eigenvalues as \(K_L^\infty\), counting multiplicities.
Consequently,
\[
\Ee\bigl[e^{\theta N_\infty(L)}\bigr]
=
\det\!\bigl(\Id+(e^\theta-1)T_L^\infty\bigr),
\qquad \theta\in\Cc.
\]
\end{lemma}

\begin{proof}
Use $\mathsf J_L$ from the preceding proof.
Its adjoint product has kernel
\[
(\mathsf J_L^*\mathsf J_L)(\xi,\eta)
=\frac{\sqrt{\xi\eta}}\pi
\int_0^L\!\int_1^\infty
 e^{\ii(\eta-\xi)x}e^{-(\xi+\eta)y}\,dy\,dx
=\frac{\sqrt{\xi\eta}\,e^{-(\xi+\eta)}}
{\pi(\xi+\eta)}X_L(\eta-\xi).
\]
Complex conjugation $\mathcal Cf=\overline f$ on $L^2(0,\infty)$ changes the phase to $X_L(\xi-\eta)$, so $T_L^\infty=\mathcal C\mathsf J_L^*\mathsf J_L\mathcal C$.
The nonzero eigenvalues of $\mathsf J_L^*\mathsf J_L$ and $\mathsf J_L\mathsf J_L^*$ agree, with multiplicities; antiunitary conjugation preserves their nonnegative eigenvalues.
The Hilbert--Schmidt bound in the preceding proof gives positivity, trace class, and contractivity.
The determinant identity follows from Lemma~\ref{lem:fredholm-rectangles-halfline}.
\end{proof}

\subsection{Height truncation and the thermodynamic limit}\label{subsec:thermodynamic-halfline}

Recall $r_\infty(\xi)=e^{-2\xi}$ and put $r_b(\xi)=e^{-2\xi}-e^{-2b\xi}$ for $b>1$ and $\xi>0$.
Let $T_L^{[1,b]}$ be the Fourier--Laplace compression associated with heights in $[1,b]$.
The factorization in Lemma~\ref{lem:reduction-halfline}, applied to $[1,b]$ and $[b,\infty)$, gives
\[
0\le T_L^{[1,b]}\le T_L^\infty\le\Id,
\qquad
\|T_L^\infty-T_L^{[1,b]}\|_1=\frac{L}{4\pi b}.
\]

\begin{lemma}
For every $K\Subset S_\pi$, there is $C_K<\infty$ such that, for $b>1$ and $L>0$,
\begin{equation}\label{eq:height-truncation-log}
\sup_{\theta\in K}
\left|\log\det(\Id+(e^\theta-1)T_L^\infty)
-\log\det(\Id+(e^\theta-1)T_L^{[1,b]})\right|
\le C_K L/b,
\end{equation}
with both logarithms normalized at zero.
\end{lemma}

\begin{proof}
Set $A=T_L^\infty$, $B=T_L^{[1,b]}$, $A_s=B+s(A-B)$, and $q_\theta=e^\theta-1$.
Each $A_s$ is a positive contraction, so $(\Id+q_\theta A_s)^{-1}$ has uniformly bounded norm for $\theta\in K$ and $s\in[0,1]$.
Integration of
\[
\partial_s\log\det(\Id+q_\theta A_s)
=q_\theta\Tr[(\Id+q_\theta A_s)^{-1}(A-B)]
\]
proves \eqref{eq:height-truncation-log}.
\end{proof}

\begin{proof}[Proof of Theorem~\ref{thm:thermo-halfline}]
For fixed $b$, \eqref{eq:intro-recalled-interval-strong} gives the limiting pressure $\Lambda_{[1,b]}$ locally uniformly on $S_\pi$.
The logarithm is uniformly Lipschitz in the scalar spectral variable on compact tilt sets, so $r_\infty-r_b=e^{-2b\xi}$ implies
\[
\sup_{\theta\in K}\left|
\Lambda_{[1,b]}(\theta)
-\frac1{2\pi}\int_0^\infty\Logp(1+(e^\theta-1)e^{-2\xi})\,d\xi
\right|\le C_K/b.
\]
Together with \eqref{eq:height-truncation-log}, taking first $L\to\infty$ and then $b\to\infty$ proves the asserted limit.
The substitution $u=e^{-2\xi}$ yields
\[
\Lambda_\infty(\theta)
=\frac1{4\pi}\int_0^1\frac{\Logp(1+(e^\theta-1)u)}u\,du
=-\frac1{4\pi}\Li_2(1-e^\theta).
\]
For $\theta\in S_\pi$, the segments $(1-u)+ue^\theta$, $0\le u\le1$, avoid $(-\infty,0]$.
After removal of the singularity at $u=0$, the integrand is holomorphic and locally uniformly bounded, together with its derivatives.
Thus $\Lambda_\infty$ is holomorphic on $S_\pi$, with
$\Lambda_\infty'(0)=1/(4\pi)$ and $\Lambda_\infty''(0)=1/(8\pi)$.
\end{proof}

\subsection{Convex duality and large deviations}\label{subsec:convexity-halfline}

We next identify the slope range and convex dual of $\Lambda_\infty$.
The zero-density endpoint deserves attention because the symbol reaches one at the origin.

\begin{lemma}\label{lem:Lambda-halfline-smooth}
The function \(\Lambda_\infty\) defined by \eqref{eq:thermo-limit-halfline} belongs to \(C^\infty(\Rr)\).
For \(\theta\in\Rr\setminus\{0\}\),
\begin{equation}\label{eq:Lambda-derivatives-halfline}
\Lambda_\infty'(\theta)
=\frac{\theta e^\theta}{4\pi(e^\theta-1)},
\qquad
\Lambda_\infty''(\theta)
=\frac{e^\theta(e^\theta-1-\theta)}
{4\pi(e^\theta-1)^2}>0.
\end{equation}
Moreover, \(\Lambda_\infty'(0)=1/(4\pi)\) and \(\Lambda_\infty''(0)=1/(8\pi)\), so \(\Lambda_\infty\) is strictly convex on \(\Rr\), and
\[
\lim_{\theta\to-\infty}\Lambda_\infty'(\theta)=0,
\qquad
\lim_{\theta\to+\infty}\Lambda_\infty'(\theta)=+\infty.
\]
\end{lemma}

\begin{proof}
Smoothness and the values at zero follow from Theorem~\ref{thm:thermo-halfline}.
Differentiating the integral in \eqref{eq:thermo-limit-halfline} gives \eqref{eq:Lambda-derivatives-halfline}, justified on each compact \(K\Subset\Rr\) by the bound \(C_K e^{-2\xi}\) for the absolute values of the first two \(\theta\)-derivatives of the integrand.
The inequality \(e^\theta>1+\theta\) for \(\theta\ne0\), together with \(\Lambda_\infty''(0)>0\), proves strict convexity.
The slope limits follow from the formula for \(\Lambda_\infty'\).
\end{proof}

\begin{proposition}\label{prop:Legendre-halfline}
Extend \(J_\infty\) to \(\Rr\) using the Legendre--Fenchel formula in \eqref{eq:rate-halfline}.
Then \(J_\infty\) is lower semicontinuous on \(\Rr\), equals \(+\infty\) on \((-\infty,0)\), and is finite on \([0,\infty)\).
For every \(x>0\), there is a unique \(\theta_x\in\Rr\) satisfying \(\Lambda_\infty'(\theta_x)=x\), and
\[
J_\infty(x)=\theta_xx-\Lambda_\infty(\theta_x),
\qquad
J_\infty'(x)=\theta_x.
\]
Moreover, \(J_\infty\) is strictly convex on \((0,\infty)\), with \(J_\infty(0)=\lim_{x\downarrow0}J_\infty(x)=\pi/24\).
Its restriction to \([0,\infty)\) is a good rate function with unique zero \(\lambda_\infty=1/(4\pi)\).
\end{proposition}

\begin{proof}
By Lemma~\ref{lem:Lambda-halfline-smooth}, $\Lambda_\infty'$ is strictly increasing from $0$ to $+\infty$.
Thus for $x>0$ there is a unique saddle $\theta_x$, and it maximizes $\theta x-\Lambda_\infty(\theta)$.
Standard differentiation of this maximized expression gives $J_\infty'(x)=\theta_x$ and $J_\infty''(x)=1/\Lambda_\infty''(\theta_x)>0$.

The behavior at zero follows from the integrable logarithmic singularity of $\log(1-e^{-2\xi})$ at $\xi=0$:
\[
J_\infty(0)
=-\lim_{\theta\to-\infty}\Lambda_\infty(\theta)
=-\frac1{2\pi}\int_0^\infty\log(1-e^{-2\xi})\,d\xi
=\frac1{4\pi}\sum_{n\ge1}\frac1{n^2}
=\frac\pi{24}.
\]
For $x<0$, the expression $\theta x-\Lambda_\infty(\theta)$ tends to $+\infty$ as $\theta\to-\infty$, so $J_\infty(x)=+\infty$.
As $x\downarrow0$, one has $\theta_x\to-\infty$ and $\theta_xx=\theta_x\Lambda_\infty'(\theta_x)\to0$ by \eqref{eq:Lambda-derivatives-halfline}.
Hence the maximizing formula gives $J_\infty(x)\to J_\infty(0)$.

The supremum representation makes $J_\infty$ lower semicontinuous.
For any fixed $\theta_0>0$, $J_\infty(x)\ge\theta_0x-\Lambda_\infty(\theta_0)$, so its sublevel sets in $[0,\infty)$ are compact.
Finally $J_\infty\ge0$, and the saddle $\theta_x=0$ occurs only at $x=\Lambda_\infty'(0)=1/(4\pi)$, where $J_\infty$ vanishes.
Strict convexity and $J_\infty(0)>0$ give uniqueness.
\end{proof}

\begin{proof}[Proof of Theorem~\ref{thm:LDP-halfline}]
The normalized logarithmic moment generating functions of $Y_L=N_\infty(L)/L$ converge to $\Lambda_\infty$ by Theorem~\ref{thm:thermo-halfline}.
The limit is finite and differentiable on all of $\Rr$; its effective domain has no boundary, so the steepness condition in the G\"artner--Ellis theorem is vacuous.
The full theorem \cite[Theorem~2.3.6]{DemboZeitouni2010} gives an LDP on $\Rr$ with the good rate function of Proposition~\ref{prop:Legendre-halfline}.
Since $Y_L\ge0$ and the rate is infinite on $(-\infty,0)$, this restricts to the relative topology of $[0,\infty)$, including its endpoint.
\end{proof}

The rate function admits the following global parametrization.

\begin{corollary}\label{cor:rate-parametrization-halfline}
For \(s>0\), set \(x(s):=s\log s/[4\pi(s-1)]\), with the continuous value \(x(1):=1/(4\pi)\).
Then \(x:(0,\infty)\to(0,\infty)\) is a strictly increasing bijection, and
\[
J_\infty(x(s))
=
\frac{1}{4\pi}
\left[
\frac{s}{s-1}(\log s)^2+\Li_2(1-s)
\right],
\qquad s>0,
\]
where the right-hand side is understood by continuity at \(s=1\), with value \(0\).
Moreover, \(J_\infty(x)=2\pi x^2-\pi/24+o(1)\) as \(x\to\infty\).
\end{corollary}

\begin{proof}
With $s=e^\theta$, Lemma~\ref{lem:Lambda-halfline-smooth} gives $x(s)=\Lambda_\infty'(\log s)$.
The same lemma implies strict monotonicity and the stated range.
Substitution in the maximizing formula of Proposition~\ref{prop:Legendre-halfline} gives the displayed expression for $J_\infty(x(s))$, and its value at $s=1$ follows by continuity.

For $\theta\to\infty$, the dilogarithm inversion formula gives
\[
\Lambda_\infty(\theta)
=\frac{\theta^2}{8\pi}+\frac\pi{24}
+O(\theta e^{-\theta}).
\]
The saddle relation $x=\theta_x/[4\pi(1-e^{-\theta_x})]$ yields $\theta_x=4\pi x+o(1)$.
Completing the square in $x\theta_x-\theta_x^2/(8\pi)$ gives the stated high-density asymptotic.
\end{proof}

We next describe the behavior of the rate function near zero.
The resulting expansion will also be used to make the dilute-count asymptotics in Section~\ref{sec:endpoint} explicit.

\begin{corollary}\label{cor:low-density-rate}
Fix $a>0$, and set $\ell_{a,x}:=\log(1/(4\pi ax))$.
As $x\downarrow0$,
\[
J_a^\infty(0)-J_a^\infty(x)
=x\left[
\ell_{a,x}+\log\ell_{a,x}+1
+O\!\left(\frac{\log\ell_{a,x}}{\ell_{a,x}}\right)
\right].
\]
\end{corollary}

\begin{proof}
For sufficiently small $x>0$, let $t_x>0$ satisfy
$(\Lambda_a^\infty)'(-t_x)=x$.
The saddle equation gives $4\pi ax=t_x/(e^{t_x}-1)$, so
$\ell_{a,x}=t_x-\log t_x+O(e^{-t_x})$.
Since $t_x\to\infty$, this implies
$t_x=\ell_{a,x}+O(\log\ell_{a,x})$, and substitution yields
\[
t_x=\ell_{a,x}+\log\ell_{a,x}
+O\!\left(\frac{\log\ell_{a,x}}{\ell_{a,x}}\right).
\]
Legendre duality and continuity at zero give
$J_a^\infty(0)-J_a^\infty(x)=\int_0^x t_u\,du$.
Using $\ell_{a,u}=\ell_{a,x}+\log(x/u)$, we obtain
$\int_0^x\ell_{a,u}\,du=x(\ell_{a,x}+1)$ and
$\int_0^x\log\ell_{a,u}\,du
=x\log\ell_{a,x}+O(x/\ell_{a,x})$.
The error term integrates to
$O(x\log\ell_{a,x}/\ell_{a,x})$, proving the claim.
\end{proof}

\section{Canonical Fisher--Hartwig asymptotics}\label{sec:WH-halfline}

The preceding section gives the leading term of the logarithmic moment generating function.
We now determine the logarithmic and constant terms for its scalar Wiener--Hopf comparison.
The additional contribution from the two-frequency kernel is treated in the next sections.

\subsection{Wiener--Hopf representation and pure-jump reference}

Extend \(r_\infty(\xi)=e^{-2\xi}\) by zero to \(\xi\le0\).
We use the Fourier convention
\[
\widehat b(x)
:=
\frac{1}{2\pi}\int_\Rr b(\xi)e^{-\ii x\xi}\,d\xi.
\]
For a bounded symbol \(a\), let \(a(D)\) denote the Fourier multiplier with convolution kernel \(\widehat a(x-y)\), interpreted distributionally, and set \(\widetilde a(\xi):=a(-\xi)\), so that \(a(-D)=\widetilde a(D)\).
Let \(W_L(a)\) and \(W_+(a)\) denote the compressions of \(a(D)\) to \(L^2(0,L)\) and \(L^2(0,\infty)\), respectively.
We write $\mathfrak S_1,\mathfrak S_2$ for the trace-class and Hilbert--Schmidt ideals, with norms $\|\cdot\|_1,\|\cdot\|_2$, and use $\|\cdot\|$ for the operator norm.
If \(a-1\in L^2(\Rr)\), then \(W_L(a)-\Id\in\mathfrak S_2\).
For \(K\in\mathfrak S_2\), define
\[
\operatorname{det}_{2}(\Id+K)
:=
\det\bigl((\Id+K)e^{-K}\bigr);
\]
see \cite[Lemma~9.1 and Theorem~9.2]{Simon2005} for its construction and continuity properties.

\begin{lemma}\label{lem:unitary-WH-halfline}
Let \(\widetilde W_L^\infty\) be the integral operator on \(L^2(0,\infty)\) with kernel
\[
\widetilde{\mathcal W}_L^\infty(\xi,\eta)
=
\frac{e^{-(\xi+\eta)}}{2\pi}X_L(\xi-\eta).
\]
Then \(\widetilde W_L^\infty\) is a positive trace-class contraction with \(\Tr\widetilde W_L^\infty=L/(4\pi)\).
Moreover, for every \(q\in\Cc\),
\begin{equation}\label{eq:det-unitary-halfline}
\det(\Id+q\widetilde W_L^\infty)
=
\det(\Id+qW_L(r_\infty)).
\end{equation}
\end{lemma}

\begin{proof}
Define \(V:L^2(0,L)\to L^2(0,\infty)\) by
\[
(Vf)(\xi)
:=
\frac{e^{-\xi}}{\sqrt{2\pi}}
\int_0^L e^{\ii\xi x}f(x)\,dx.
\]
By direct computation, \(V^*V=W_L(r_\infty)\), \(VV^*=\widetilde W_L^\infty\), and \(\|V\|_2^2=L/(4\pi)\).
Since \(0\le r_\infty\le1\), the operator \(W_L(r_\infty)=V^*V\) is a positive contraction.
Hence \(\widetilde W_L^\infty=VV^*\) is a positive trace-class contraction with \(\Tr\widetilde W_L^\infty=\|V\|_2^2=L/(4\pi)\).
Finally, \(V^*V\) and \(VV^*\) have the same nonzero eigenvalues, counting multiplicities, which gives \eqref{eq:det-unitary-halfline}.
\end{proof}

Set
\[
\sigma_\theta(\xi)
:=
1+(e^\theta-1)e^{-2\xi}\mathbf1_{(0,\infty)}(\xi),
\qquad
\ell_1(\xi):=\Logp\sigma_\theta(\xi),
\qquad
\theta\in S_\pi.
\]
For \(\xi>0\), the value \(\sigma_\theta(\xi)\) lies on the segment joining \(1\) and \(e^\theta\), which avoids \((-\infty,0]\) for \(\theta\in S_\pi\).
Consequently, \(\ell_1\) is well defined and holomorphic in \(\theta\); on every compact \(K\Subset S_\pi\), the distance of \(\sigma_\theta(\xi)\) from \((-\infty,0]\) is uniformly positive for \(\theta\in K\) and \(\xi\in\Rr\).

For \(\xi\ne0\), set
\[
p(\xi):=\frac{1}{\pi}\arctan(1/\xi),
\qquad
\ell_0(\xi):=\theta p(\xi),
\qquad
a_0:=e^{\ell_0},
\qquad
\beta:=\frac{\theta}{2\pi\ii}.
\]
The function \(p\) is odd, \(p(0+)-p(0-)=1\), and \(p(\xi)\to0\) as \(\lvert\xi\rvert\to\infty\).
We use the following specialization of the pure-jump Wiener--Hopf asymptotic of Basor and Widom~\cite{BasorWidom2004}, where \(G\) denotes the Barnes \(G\)-function, normalized by \(G(1)=1\) and \(G(z+1)=\Gamma(z)G(z)\).

\begin{lemma}[Pure-jump Wiener--Hopf asymptotic]\label{lem:new-pure-FH}
As \(L\to\infty\), locally uniformly for \(\theta\in S_\pi\),
\begin{equation}\label{eq:new-pure-FH}
\begin{aligned}
\log\operatorname{det}_{2}W_L(a_0)
={}&
\frac{L}{2\pi}\int_\Rr(\ell_0-a_0+1)\,d\xi
-\beta^2\log L+\beta^2\log2\\
&+\log G(1+\beta)+\log G(1-\beta)+o(1).
\end{aligned}
\end{equation}
The logarithms of the determinant and the Barnes \(G\)-factors are the holomorphic branches normalized to vanish at \(\theta=0\).
\end{lemma}

\begin{proof}
Since \(p(\xi)=O(|\xi|^{-1})\), we have
\(a_0-1\in L^2(\Rr)\) and \(\ell_0-a_0+1\in L^1(\Rr)\).
For \(\xi\ne0\), use
\(\Logp(\xi\pm\ii0^+)=\log|\xi|\pm\ii\pi\mathbf1_{\{\xi<0\}}\).
Then
\[
\beta\Logp\frac{\xi+\ii}{\xi+\ii0^+}
-\beta\Logp\frac{\xi-\ii}{\xi-\ii0^+}
=2\ii\beta\arctan(1/\xi)=\ell_0(\xi).
\]
Thus \(a_0\) is the Basor--Widom symbol with parameters
\((\beta,-\beta)\), and their admissibility condition becomes
\(\lvert\operatorname{Re}\beta\rvert<1/2\), equivalently
\(\theta\in S_\pi\).

Let \(\varphi_{\beta,-\beta}(e^{\ii t})
=e^{\ii\beta(t-\pi)}\) for \(0<t<2\pi\), and put
\(n=\lfloor L/2\rfloor\).
Write \(T_n(\varphi)\) for the \(n\times n\) Toeplitz matrix with symbol \(\varphi\).
Write
$G_2(a_0)=\exp\{(2\pi)^{-1}\int_\Rr(\ell_0-a_0+1)\,d\xi\}$
and $G_2(\varphi)=\exp((\log\varphi)_0-\varphi_0+1)$.
The normalized regularized determinants in
\cite[Theorem, p.~133]{BasorWidom2004} satisfy
\[
\frac{\operatorname{det}_2 W_L(a_0)}{G_2(a_0)^L}
=\frac{\operatorname{det}_2 T_n(\varphi_{\beta,-\beta})}
{G_2(\varphi_{\beta,-\beta})^n}(1+o(1))
=\det T_n(\varphi_{\beta,-\beta})(1+o(1)).
\]
Indeed, $\operatorname{det}_2 T_n(\varphi)
=\det T_n(\varphi)e^{-n(\varphi_0-1)}$
and $(\log\varphi_{\beta,-\beta})_0=0$.
The convergence is locally uniform on \(S_\pi\), by the
compact-parameter estimates and analytic continuation in the cited proof.
The B\"ottcher--Widom identity
\cite[Eq.~(2), p.~593]{BottcherWidom2005}
and the large-argument expansion of the Barnes \(G\)-function give
\[
\det T_n(\varphi_{\beta,-\beta})
=\frac{G(n+1)^2G(1+\beta)G(1-\beta)}
        {G(n+1+\beta)G(n+1-\beta)}
=n^{-\beta^2}G(1+\beta)G(1-\beta)(1+o(1)),
\]
locally uniformly for \(\lvert\operatorname{Re}\beta\rvert<1/2\). Finally,
\[
\operatorname{Re}a_0(\xi)
\ge e^{-|\operatorname{Re}\theta|/2}
\cos(|\operatorname{Im}\theta|/2)>0,
\qquad \theta\in S_\pi,
\]
so \(W_L(a_0)\) is invertible.
Also, \(G(1\pm\beta)\ne0\) on \(S_\pi\), so all the normalized
holomorphic logarithms exist.
Taking these logarithms and using
\(\log n=\log L-\log2+O(L^{-1})\) proves
\eqref{eq:new-pure-FH}, locally uniformly on \(S_\pi\).
\end{proof}

\begin{remark}
The pure-jump formula can also be read from \cite[Theorem~1.1]{Kozlowski2022}.
In that notation the choice \((\nu_1,\bar\nu_1)=(\beta,-\beta)\) gives \((\gamma_1,\delta_1)=(0,-\beta)\) and the factor \((L/2)^{-\beta^2}G(1+\beta)G(1-\beta)\), in agreement with \eqref{eq:new-pure-FH}.
\end{remark}

\subsection{Regular deformation and boundary operators}

The symbols \(\sigma_\theta\) and \(a_0\) have the same multiplicative jump at the origin. To compare the normalized determinants of \(a_0\) and \(\sigma_\theta\), set
\[
h:=\ell_1-\ell_0,
\qquad
\ell_u:=\ell_0+uh,
\qquad
a_u:=e^{\ell_u},
\qquad
0\le u\le1.
\]
Then \(a_1=\sigma_\theta\), and \(h\) extends continuously across the origin with \(h(0)=\theta/2\), although \(h'\) may jump there.
For a kernel \(k\), write \(\check k(x):=k(-x)\), and let \(\mathcal H(k)\) denote the Hankel operator on \(L^2(0,\infty)\) with kernel \(k(x+y)\).
Define \(s_u:=\widehat{a_u-1}\) and \(k_h:=\widehat h\) in the \(L^2\)-sense.
The boundary factorization rests on the following uniform coercivity and kernel-decay estimates.

\begin{lemma}\label{lem:new-canonical-bounds}
For every \(K\Subset S_\pi\), there exist \(c_K,C_K>0\) such that,
uniformly for \(\theta\in K\) and \(u\in[0,1]\),
\begin{equation}\label{eq:new-coercivity}
\operatorname{Re}\!\left(
e^{-\ii u\operatorname{Im}\theta/2}a_u(\xi)
\right)
\ge c_K
\qquad\text{for a.e. }\xi\in\Rr.
\end{equation}
All half-line and finite-interval compressions of \(a_u(D)\) and
\(\widetilde a_u(D)\) are invertible, with inverse norm at most \(c_K^{-1}\).
Moreover, \(s_u\) and \(k_h\) admit representatives satisfying
\begin{equation}\label{eq:new-canonical-kernels}
|s_u(x)|\le\frac{C_K}{1+|x|},
\qquad
|k_h(x)|\le\frac{C_K}{(1+|x|)^2},
\qquad x\ne0.
\end{equation}
For each \(x\ne0\), these representatives depend continuously on
\((u,\theta)\).
\end{lemma}

\begin{proof}
Put \(v=\operatorname{Im}\theta\).
The definitions give
\(\lvert\operatorname{Im}\ell_u-uv/2\rvert\le|v|/2\),
while \(\operatorname{Re}\ell_u\) is uniformly bounded below.
This proves \eqref{eq:new-coercivity}, also for the reflected symbol.
The rotated compressions and their adjoints satisfy the same
coercive bound, giving invertibility and the stated inverse norms.

Uniformly in \((u,\theta)\in[0,1]\times K\),
\[
a_u(\xi)-1=\frac{(1-u)\theta}{\pi\xi}+O_K(\xi^{-2}),
\qquad
h(\xi)=-\frac{\theta}{\pi\xi}+O_K(\xi^{-2}),
\qquad |\xi|\to\infty.
\]
Subtracting the corresponding multiples of \(\xi/(1+\xi^2)\)
leaves remainders uniformly bounded and continuous in \(L^1(\Rr)\).
Since the rational tail has a bounded Fourier transform, this gives
uniform boundedness near zero and the asserted parameter continuity
for \(s_u\) and \(k_h\).
Direct differentiation, using continuity of \(h\) at zero, shows
that \(Da_u\) and \(D^2h\) are finite measures with uniformly bounded
total variation.
Fourier differentiation then gives the decay in
\eqref{eq:new-canonical-kernels} for \(|x|\ge1\), completing the proof.
\end{proof}

Define the boundary products
\[
K_-(u,\theta)
:=
\mathcal H(s_u)\mathcal H(\check k_h),
\qquad
K_+(u,\theta)
:=
\mathcal H(\check s_u)\mathcal H(k_h).
\]

\begin{lemma}\label{lem:new-boundary-nuclear}
For every \(K\Subset S_\pi\), the maps \(K_\pm:[0,1]\times K\to\mathfrak S_1\) are continuous and have uniformly bounded trace norms.
\end{lemma}

\begin{proof}
Fix \(K\Subset S_\pi\) and put \(w(x)=(1+x)^{1/2}\).
Factor the first product as
\[
K_-=[\mathcal H(s_u)M_{w^{-1}}]
[M_w\mathcal H(\check k_h)].
\]
By \eqref{eq:new-canonical-kernels}, both factors are Hilbert--Schmidt, uniformly on \([0,1]\times K\):
\[
\begin{aligned}
\|\mathcal H(s_u)M_{w^{-1}}\|_2^2
&\le C_K\int_0^\infty\!\int_0^\infty
\frac{dx\,dy}{(1+x+y)^2(1+y)}<\infty,\\
\|M_w\mathcal H(\check k_h)\|_2^2
&\le C_K\int_0^\infty\!\int_0^\infty
\frac{1+x}{(1+x+y)^4}\,dx\,dy<\infty.
\end{aligned}
\]
The same factorization applies to \(K_+\), with both Fourier kernels reflected.
The common integrable bounds and the parameter continuity in Lemma~\ref{lem:new-canonical-bounds} give Hilbert--Schmidt continuity of each factor.
Their products therefore form bounded, continuous families in \(\mathfrak S_1\).
\end{proof}

Let \(P_L\) denote multiplication by \(\mathbf1_{(0,L)}\) on the ambient \(L^2\)-space.
For \(a=a_u\) or \(\widetilde a_u\), let \(\mathcal R_L(a)\) denote \(W_L(a)^{-1}\), extended by zero to \(L^2(0,\infty)\).

\begin{lemma}\label{lem:new-canonical-derivative}
For \((u,\theta)\in[0,1]\times S_\pi\), define
\[
E_L(u,\theta)
:=
\log\operatorname{det}_{2}W_L(a_u)
-\frac{L}{2\pi}\int_\Rr(\ell_u-a_u+1)\,d\xi.
\]
Then
\begin{equation}\label{eq:new-canonical-derivative}
\partial_uE_L(u,\theta)
=
\Tr\!\bigl[\mathcal R_L(a_u)K_-\bigr]
+
\Tr\!\bigl[\mathcal R_L(\widetilde a_u)K_+\bigr].
\end{equation}
The traces on the right are ordinary traces of trace-class operators.
\end{lemma}

\begin{proof}
Write
\[
A=W_L(a_u),\quad C=A-\Id,\quad B=W_L(h),\quad
D=W_L((a_u-1)h).
\]
Here \(B,C\in\mathfrak S_2\) and \(D\in\mathfrak S_1\): indeed \((a_u-1)h\in L^1(\Rr)\), and factoring this symbol into two \(L^2\) functions proves the trace-class assertion and
\[
\Tr D=\frac L{2\pi}\int_\Rr(a_u-1)h\,d\xi.
\]
Also \(\partial_uA=B+D=W_L(a_uh)\).
Since \(\ell_u-a_u+1=O_K(\xi^{-2})\), the scalar integral in \(E_L\) is absolutely convergent.
Differentiation gives
\[
\begin{aligned}
\partial_uE_L
&=\Tr[(A^{-1}-\Id)(B+D)]+\Tr D=\Tr\bigl[A^{-1}(D-CB)\bigr]\\
&=\Tr\!\left[A^{-1}
\bigl(W_L(a_uh)-W_L(a_u)W_L(h)\bigr)\right].
\end{aligned}
\]
Every trace above is of a trace-class product; in particular no separate trace is assigned to \(B\).
Since
\[
W_L(a_uh)-W_L(a_u)W_L(h)
=P_La_u(D)(\Id-P_L)h(D)P_L,
\]
split the complementary projection into \(( -\infty,0)\) and \((L,\infty)\).
The first piece is the compression of \(K_-\), with kernel \(\int_0^\infty s_u(x+w)k_h(-w-y)\,dw\).
Reflection \(x\mapsto L-x\) changes the second piece into the compression of \(K_+\) and changes \(a_u\) to \(\widetilde a_u\).
Inserting the extended inverses gives \eqref{eq:new-canonical-derivative}.
\end{proof}

It remains to pass from finite intervals to the half-line in the boundary traces.

\begin{lemma}\label{lem:new-canonical-limit}
For every \(K\Subset S_\pi\), as \(L\to\infty\), uniformly for \((u,\theta)\in[0,1]\times K\),
\begin{equation}\label{eq:new-boundary-limit}
\partial_uE_L(u,\theta)
\longrightarrow
b(u,\theta)
:=
\Tr\!\bigl[W_+(a_u)^{-1}K_-\bigr]
+
\Tr\!\bigl[W_+(\widetilde a_u)^{-1}K_+\bigr].
\end{equation}
Consequently,
\(
E_L(1,\theta)-E_L(0,\theta)
\longrightarrow
\int_0^1 b(u,\theta)\,du
\)
locally uniformly for \(\theta\in S_\pi\).
\end{lemma}

\begin{proof}
Fix \(K\Subset S_\pi\) and \(g\in L^2(0,\infty)\).
Set \(v=W_+(a_u)^{-1}g\).
Then
\[
\mathcal R_L(a_u)g-P_Lv
=\mathcal R_L(a_u)P_LW_+(a_u)(\Id-P_L)v.
\]
The inverse norms are bounded uniformly on \([0,1]\times K\).
The resolvent identity makes \(v\) continuous in \((u,\theta)\); its image of this compact set is compact in \(L^2\).
Hence \((\Id-P_L)v\to0\) uniformly over the parameters, and \(\mathcal R_L(a_u)g\to W_+(a_u)^{-1}g\) uniformly for each fixed \(g\).
The same holds for the reflected symbols.

By Lemma~\ref{lem:new-boundary-nuclear}, the families \(K_\pm\) are compact in trace norm.
Approximate each uniformly by finitely many finite-rank operators.
The strong convergence just proved, together with the uniform inverse bounds, gives convergence of the traces in \eqref{eq:new-canonical-derivative}, uniformly on \([0,1]\times K\).
This is \eqref{eq:new-boundary-limit}; integration in \(u\) proves the final assertion.
\end{proof}

\subsection{Evaluation of the canonical constant}

Set \(\tau_j:=\widehat{\ell_j}\), \(j=0,1\), and define the relative energy
\begin{equation}\label{eq:new-relative-energy}
\mathcal E(\ell_1,\ell_0)
:=
\int_0^\infty
x\bigl[\tau_1(x)\tau_1(-x)-\tau_0(x)\tau_0(-x)\bigr]\,dx.
\end{equation}
The integral is absolutely convergent since \(\tau_1-\tau_0=k_h=O((1+|x|)^{-2})\) and \(\tau_j=O((1+|x|)^{-1})\).

\begin{lemma}\label{lem:new-relative-energy}
For every \(\theta\in S_\pi\), the function \(b(\cdot,\theta)\) defined by \eqref{eq:new-boundary-limit} satisfies
\begin{equation}\label{eq:new-b-energy}
\int_0^1 b(u,\theta)\,du
=
\mathcal E(\ell_1,\ell_0).
\end{equation}
\end{lemma}

\begin{proof}
Fix \(\theta\in S_\pi\).
Let \(\varrho_\epsilon\) be an even, nonnegative, compactly supported smooth approximate identity and set \(\ell_{u,\epsilon}=\varrho_\epsilon*\ell_u\), \(h_\epsilon=\varrho_\epsilon*h\), and \(a_{u,\epsilon}=e^{\ell_{u,\epsilon}}\).
Define \(K_{\pm,\epsilon}\) as above from \(\widehat{a_{u,\epsilon}-1}\) and \(\widehat h_\epsilon\), and put
\[
 b_\epsilon(u,\theta)
 =\Tr[W_+(a_{u,\epsilon})^{-1}K_{-,\epsilon}]
 +\Tr[W_+(\widetilde a_{u,\epsilon})^{-1}K_{+,\epsilon}].
\]
Smoothing the logarithms preserves the common coercive sector and the \(1/\xi\) tails.
The total variations of \(D\ell_{u,\epsilon}\) and \(D^2h_\epsilon\) are bounded uniformly in \(\epsilon\) and \(u\).
The difference \(h_\epsilon\) has Fourier transform \(\widehat h_\epsilon=m_\epsilon k_h\), where \(m_\epsilon(x)\to1\) and \(|m_\epsilon(x)|\le1\).
Consequently the bounds in \eqref{eq:new-canonical-kernels} hold for these smoothed kernels with constants independent of \(\epsilon\).
The bounded-variation translation estimate, after subtracting the common rational tails, gives uniformly in \(u\),
\[
\|\ell_{u,\epsilon}-\ell_u\|_1
+\|a_{u,\epsilon}-a_u\|_1=O_\theta(\epsilon).
\]
The weighted Hilbert--Schmidt factorization in Lemma~\ref{lem:new-boundary-nuclear} now gives \(K_{\pm,\epsilon}\to K_\pm\) in trace norm.
The Wiener--Hopf compressions converge strongly; uniform coercivity gives strong convergence of their inverses.
As in Lemma~\ref{lem:new-canonical-limit}, finite-rank approximation then yields \(b_\epsilon(u,\theta)\to b(u,\theta)\), uniformly in \(u\in[0,1]\).

For each fixed \(\epsilon>0\), the regularized continuous Borodin--Okounkov identity \cite[Appendix]{BasorWidom2004} applies to \(a_{0,\epsilon}\) and \(a_{1,\epsilon}\).
Indeed, \(a_{j,\epsilon}-1\in L^2\), its Fourier transform belongs to \(L^1\) and has a finite \(|x|\)-weighted squared integral: the smoothed symbol has integrable second derivative, while subtraction of its rational tail bounds the transform near zero.
Moreover, \(\ell_{j,\epsilon}\) is a continuous logarithm vanishing at both infinities, so the symbol is nonzero and has index zero.
These are the hypotheses of the cited identity.

Define \(E_{L,\epsilon}\) by the formula for \(E_L\), with
\((\ell_u,a_u)\) replaced by \((\ell_{u,\epsilon},a_{u,\epsilon})\).
The identity gives, for \(j=0,1\),
\[
E_{L,\epsilon}(j,\theta)\longrightarrow
\int_0^\infty x\widehat{\ell_{j,\epsilon}}(x)
\widehat{\ell_{j,\epsilon}}(-x)\,dx.
\]
Its remaining determinant tends to one: the corresponding shifted Hankel product has trace norm controlled by tails of integrable \(|x|\)-weighted Fourier squares.
Integrating the smoothed version of \eqref{eq:new-canonical-derivative} therefore gives
\[
\int_0^1 b_\epsilon(u,\theta)\,du
=\int_0^\infty x\bigl[
\widehat{\ell_{1,\epsilon}}(x)\widehat{\ell_{1,\epsilon}}(-x)
-\widehat{\ell_{0,\epsilon}}(x)\widehat{\ell_{0,\epsilon}}(-x)
\bigr]dx.
\]
Here \(\widehat{\ell_{j,\epsilon}}=m_\epsilon\tau_j\), with \(m_\epsilon\) real and even.
The difference of the unsmoothed quadratic products is absolutely integrable after multiplication by \(x\), by \eqref{eq:new-relative-energy}.
Thus \(|m_\epsilon|\le1\) permits passage to the limit on the right, while the uniform convergence of \(b_\epsilon\) permits it on the left.
We obtain \eqref{eq:new-b-energy}.
\end{proof}

\begin{theorem}[Canonical Fisher--Hartwig asymptotic]
\label{thm:canonical-FH-halfline}
There exists a holomorphic function \(A_{\rm FH}\) on \(S_\pi\) such that
\begin{equation}\label{eq:canonical-main-global-halfline}
\log\det(\Id+(e^\theta-1)\widetilde W_L^\infty)
=
L\Lambda_\infty(\theta)
+\frac{\theta^2}{4\pi^2}\log L
+A_{\rm FH}(\theta)+o(1)
\end{equation}
locally uniformly for \(\theta\in S_\pi\).
\end{theorem}

\begin{proof}
Lemma~\ref{lem:new-pure-FH} gives the locally uniform limit of \(E_L(0,\theta)+\beta^2\log L\), while Lemmas~\ref{lem:new-canonical-limit} and \ref{lem:new-relative-energy} give that of \(E_L(1,\theta)-E_L(0,\theta)\).
At \(u=1\), the symbol \(a_1-1\in L^1(\Rr)\), so the ordinary Fredholm determinant is defined and
\[
\begin{aligned}
\log\det W_L(a_1)
&=\log\operatorname{det}_{2}W_L(a_1)
+\frac L{2\pi}\int_\Rr(a_1-1)\,d\xi\\
&=E_L(1,\theta)+\frac L{2\pi}\int_\Rr\ell_1(\xi)\,d\xi
=E_L(1,\theta)+L\Lambda_\infty(\theta).
\end{aligned}
\]
For the reference symbol \(a_0\), we continue to use the regularized determinant.
Since \(-\beta^2=\theta^2/(4\pi^2)\), the determinant identity \eqref{eq:det-unitary-halfline} proves \eqref{eq:canonical-main-global-halfline}.
The logarithms are holomorphic on \(S_\pi\), normalized to vanish at \(0\).
Their locally uniform limit defines the asserted holomorphic constant.
\end{proof}

Set \(\tau_\theta:=\widehat{\ell_1}\) and define
\begin{equation}\label{eq:Eren-halfline}
\mathcal E_{\rm ren}(\theta)
:=
\int_0^1 x\tau_\theta(x)\tau_\theta(-x)\,dx
+
\int_1^\infty
\left[
x\tau_\theta(x)\tau_\theta(-x)+\frac{\beta^2}{x}
\right]\,dx.
\end{equation}
The jump of \(\ell_1\) gives, locally uniformly in \(\theta\in S_\pi\), \(\tau_\theta(x)=\beta/x+O(x^{-2})\) and \(\tau_\theta(-x)=-\beta/x+O(x^{-2})\) as \(x\to\infty\).
Both transforms are bounded near zero.
Thus the subtraction in \eqref{eq:Eren-halfline} removes the only logarithmic divergence, and \(\mathcal E_{\rm ren}\) is holomorphic on \(S_\pi\).
The constant has the following explicit representation.

\begin{proposition}\label{prop:explicit-AFH-halfline}
For \(\theta\in S_\pi\),
\begin{equation}\label{eq:explicit-AFH-halfline}
A_{\rm FH}(\theta)
=
\mathcal E_{\rm ren}(\theta)
+\log G(1+\beta)+\log G(1-\beta)
+\gamma_{\rm E}\beta^2,
\qquad
\beta=\frac{\theta}{2\pi\ii},
\end{equation}
where the logarithms of the Barnes \(G\)-factors are the holomorphic
branches normalized to vanish at \(\theta=0\), and \(\gamma_{\rm E}\)
is the Euler--Mascheroni constant.
Moreover, \(A_{\rm FH}''(0)=(1-\log2)/(2\pi^2)\).
\end{proposition}

\begin{proof}
The distributional identity
\(Dp=\delta_0-(\pi(1+\xi^2))^{-1}d\xi\) gives
\(\tau_0(x)=\beta(1-e^{-|x|})/x\) for \(x\ne0\).
Using
\(\int_0^R(1-e^{-x})^2\,dx/x
=\log R+\gamma_{\rm E}-\log2+o(1)\),
we obtain
\[
\mathcal E(\ell_1,\ell_0)
=\mathcal E_{\rm ren}(\theta)
-\beta^2(\log2-\gamma_{\rm E}).
\]
Substitution into the constant obtained in the proof of
Theorem~\ref{thm:canonical-FH-halfline} proves
\eqref{eq:explicit-AFH-halfline}.

For the second derivative, direct integration of the kernel
\((2\pi(2+\ii(x-y)))^{-1}\) of \(W_L(r_\infty)\) gives
\[
\begin{aligned}
\Tr\bigl(W_L(r_\infty)(\Id-W_L(r_\infty))\bigr)
&=\frac{L}{4\pi}
-\frac{L}{4\pi^2}\arctan(L/2)
+\frac{1}{4\pi^2}\log(1+L^2/4)\\
&=\frac{L}{8\pi}+\frac{\log L}{2\pi^2}
+\frac{1-\log2}{2\pi^2}+o(1).
\end{aligned}
\]
The left side is the second derivative at \(\theta=0\) of the
canonical determinant logarithm.
Cauchy's formula permits differentiation of
\eqref{eq:canonical-main-global-halfline}; comparison of constant
terms, using \(\Lambda_\infty''(0)=1/(8\pi)\), gives the stated
value of \(A_{\rm FH}''(0)\).
\end{proof}

\section{Endpoint correction to the canonical Wiener--Hopf operator}
\label{sec:hardedge-halfline}

We compare the Fourier--Laplace operator with its Wiener--Hopf counterpart by localizing their difference to dyadic frequency blocks.
This permits trace-class estimates without assuming that the unweighted global difference is trace class.

\subsection{Homogeneous model and trace-class dyadic blocks}

We use homogeneous finite-band and one-edge operators so that the dyadic estimates remain uniform under exponential weighting.
Write \(T_L:=T_L^\infty\) and \(W_L:=\widetilde W_L^\infty\).
For \(t>0\), define
\begin{align*}
\mathsf W_t(u,v)&=\frac{1}{2\pi}X_t(u-v),&
\mathsf T_t(u,v)&=\frac{\sqrt{uv}}{\pi(u+v)}X_t(u-v).
\end{align*}
Set \(\mathsf A_{t,s}:=(1-s)\mathsf W_t+s\mathsf T_t\), \(0\le s\le1\), and \(\mathsf D_t:=\mathsf T_t-\mathsf W_t\).
For \(r\ge0\), let \(M_r\) denote multiplication by \(e^{-ru}\), and set
\[
A_{r,t,s}:=M_r\mathsf A_{t,s}M_r,\qquad
D_{r,t}:=M_r\mathsf D_tM_r.
\]
If \(r>0\) and \(t=Lr\), the dilation \((\mathcal U_rf)(u)=r^{1/2}f(ru)\) makes \(M_r\mathsf T_tM_r\) and \(M_r\mathsf W_tM_r\) unitarily equivalent to \(T_L\) and \(W_L\), respectively.

For \(\mathcal M\in\{T,W\}\), let \(\mathsf A_\pm^{\mathcal M}\) be the corresponding one-edge operators obtained by replacing \(X_t\) with
\[
X_\pm(v)=\pi\delta_0(v)\pm\ii\operatorname{pv}\frac1v.
\]
These models describe a single endpoint of the horizontal interval: the \(+\) model corresponds to \([0,\infty)\), and the \(-\) model to \((-\infty,0]\).
Set \(\mathsf A_{\pm,s}:=(1-s)\mathsf A_\pm^W+s\mathsf A_\pm^T\), and
\[
A_{r,\pm}^{\mathcal M}:=M_r\mathsf A_\pm^{\mathcal M}M_r,\qquad
A_{r,\pm,s}:=M_r\mathsf A_{\pm,s}M_r,\qquad
D_{r,\pm}:=A_{r,\pm}^T-A_{r,\pm}^W.
\]

\begin{lemma}\label{lem:new-dephasing}
For \(t>0\),
\begin{equation}\label{eq:new-dephasing}
\mathsf T_t
=
\int_{\Rr}\operatorname{sech}(\pi y)\,
V_y\mathsf W_tV_y^*\,dy,
\qquad
(V_yf)(u)=u^{\ii y}f(u),
\end{equation}
in the weak operator sense.
The same identity holds with \((\mathsf T_t,\mathsf W_t)\) replaced by \((\mathsf A_\pm^T,\mathsf A_\pm^W)\).
Consequently, for \(s\in[0,1]\) and \(r\ge0\), \(\mathsf A_{t,s}\), \(\mathsf A_{\pm,s}\), \(A_{r,t,s}\), and \(A_{r,\pm,s}\) are positive contractions.
\end{lemma}

\begin{proof}
The operators \(\mathsf W_t\) and \(\mathsf A_\pm^W\) are compressions of Fourier projections.
Identity \eqref{eq:new-dephasing} follows from
\[
\int_{\Rr}\operatorname{sech}(\pi y)e^{\ii yv}\,dy
=
\operatorname{sech}(v/2),
\qquad
\operatorname{sech}\!\left(\tfrac12\log(u/v)\right)
=
\frac{2\sqrt{uv}}{u+v}.
\]
The same computation gives the one-edge identity.
Since \(\operatorname{sech}(\pi y)\) is nonnegative and has integral one, \(\mathsf T_t\) and \(\mathsf A_\pm^T\) are convex averages of unitary conjugates of the corresponding \(W\)-model operators.
The assertions for the affine interpolants and their \(M_r\)-weighted versions follow immediately.
\end{proof}

Let \(P_j\) denote multiplication by \(\mathbf1_{[2^j,2^{j+1})}\), \(j\in\Zz\).
The next estimates control both resolvents and localized differences.

\begin{lemma}\label{lem:new-blocks}
There are \(a,C>0\) such that, uniformly in \(r\ge0\), \(t>0\), \(s\in[0,1]\), and \(j,l\in\Zz\),
\begin{equation}\label{eq:new-block-nuclear}
\|P_jA_{r,t,s}P_l\|\le Ce^{-a|j-l|},\qquad
\|P_jD_{r,t}P_l\|_1\le Ce^{-a|j-l|}e^{-r(2^j+2^l)}.
\end{equation}
The estimate \eqref{eq:new-block-nuclear} remains valid for the weighted one-edge operators $A_{r,\pm,s}$ and $D_{r,\pm}$.
\end{lemma}

\begin{proof}
On separated dyadic intervals, $|j-l|\ge3$ implies $|u-v|\asymp2^{\max(j,l)}$.
Since $|X_t(u-v)|\le2/|u-v|$ and $2\sqrt{uv}/(u+v)\le1$, the finite-band blocks of both models have Hilbert--Schmidt norm at most $C2^{-|j-l|/2}$.
For neighboring blocks, the contraction property gives a uniform bound.
Multiplication by $M_r$ only decreases these bounds.
The same argument applies to the one-edge operators; their delta terms vanish between separated intervals.

The trace-norm estimates use the cancellation
\[
\frac{\sqrt{uv}}{\pi(u+v)}-\frac1{2\pi}
=-\frac{(u-v)^2}{2\pi(u+v)(\sqrt u+\sqrt v)^2}.
\]
With $b(u,v):=\ii(u-v)/[2\pi(u+v)(\sqrt u+\sqrt v)^2]$, the kernel of $\mathsf D_t$ is $(e^{\ii t(u-v)}-1)b(u,v)$; the one-edge discrepancy kernels are $-b$ and $b$.
There is no remaining diagonal distribution in either difference.

Choose $\zeta\in C_c^\infty((1/2,4))$ equal to one on $[1,2]$.
On comparable dyadic intervals, dilation transforms $\zeta(2^{-j}u)b(u,v)\zeta(2^{-l}v)$ into a smooth kernel on a fixed rectangle.
If $j\ge l+3$, put $\varepsilon=2^{l-j}\le1/8$ and dilate the variables separately.
The transformed kernel is
\[
\sqrt\varepsilon\,
\frac{\ii(x-\varepsilon y)}
{2\pi(x+\varepsilon y)(\sqrt x+\sqrt\varepsilon\sqrt y)^2}
\zeta(x)\zeta(y).
\]
After removing $\sqrt\varepsilon$, its mixed derivatives of order at most two in each variable are bounded uniformly in $\varepsilon$.
For a smooth kernel supported on a fixed rectangle, its trace norm is bounded by these derivative norms: extend it smoothly to a torus and expand in a double Fourier series, whose rank-one coefficients are bounded by $C(1+m^2)^{-1}(1+n^2)^{-1}$.
Since $P_j\zeta(2^{-j}\cdot)=P_j$, the ideal property transfers this estimate to the sharp dyadic blocks.
It gives $C2^{-(j-l)/2}$ for the trace norm; the opposite ordering is symmetric.
Multiplication by $e^{\ii t(u-v)}$ is unitary conjugation, so it does not introduce any $t$-dependent derivatives.
Since $\|M_rP_j\|\le e^{-r2^j}$, the nuclear estimate in \eqref{eq:new-block-nuclear} follows for all three discrepancies.

\end{proof}

The following block resolvent estimate is obtained by a standard Combes--Thomas argument; see \cite{CombesThomas1973}.
We include the short adaptation needed here.

\begin{lemma}\label{lem:new-resolvent-blocks}
Let \(\Sigma\Subset\Cc\setminus[0,1]\), and let \(A\) be a positive contraction satisfying \(\|P_jAP_l\|\le Ce^{-a|j-l|}\) for \(j,l\in\Zz\).
There are \(b,C_\Sigma>0\) such that \(Q(\lambda):=(\lambda\Id-A)^{-1}\) satisfies
\begin{equation}\label{eq:new-Qblocks}
\|P_jQ(\lambda)P_l\|
\le C_\Sigma e^{-b|j-l|},
\qquad
\lambda\in\Sigma,\quad j,l\in\Zz.
\end{equation}
The constants may be chosen uniformly for \(A_{r,t,s}\), \(A_{r,\pm,s}\), and all their multiplication compressions \(cAc\) with \(0\le c\le1\).
\end{lemma}

\begin{proof}
Let $d=\operatorname{dist}(\Sigma,[0,1])>0$.
Write $P^{(n)}=\sum_{|j|\le n}P_j$, $A_n=P^{(n)}AP^{(n)}$, and
\[
Z_{\gamma,n}=\Id-P^{(n)}+\sum_{|j|\le n}e^{\gamma j}P_j.
\]
This truncated weight is bounded and invertible.
The block Schur test gives, for $|\gamma|\le a/2$,
\[
\|Z_{\gamma,n}A_nZ_{\gamma,n}^{-1}-A_n\|
\le C\sum_{k\in\Zz}e^{-a|k|}|e^{\gamma k}-1|
\le C'|\gamma|,
\]
uniformly in $n$.
Choose $b>0$ with $C'b<d/2$.
As $0\le A_n\le\Id$, the resolvent identity bounds $\|Z_{\pm b,n}(\lambda\Id-A_n)^{-1}Z_{\pm b,n}^{-1}\|$ by $2/d$ on $\Sigma$.
Taking the sign according to the order of $j,l$ gives
\[
\|P_j(\lambda\Id-A_n)^{-1}P_l\|
\le (2/d)e^{-b|j-l|}.
\]
Now $A_n\to A$ strongly, and all resolvents have norm at most $d^{-1}$.
Passing to the strong resolvent limit proves \eqref{eq:new-Qblocks}.
If $c$ is multiplication by a function with $0\le c\le1$, then $cAc$ remains a positive contraction and its dyadic block norms do not increase.
This proves the stated uniformity for all the model families.
\end{proof}

\subsection{Localized traces and multiplication weights}

Fix a nonzero real-valued $\chi\in C_c^\infty((0,\infty))$.
Let $(\mathsf B^T,\mathsf B^W)$ be either
$(\mathsf T_t,\mathsf W_t)$, $t>0$, or
$(\mathsf A_\sigma^T,\mathsf A_\sigma^W)$, $\sigma=\pm$.
For a measurable weight $m:(0,\infty)\to[0,1]$ and
$f$ holomorphic near $[0,1]$, define
\[
\mathfrak D_m^\chi[f]
:=\chi\bigl[f(m\mathsf B^Tm)-f(m\mathsf B^Wm)\bigr]\chi.
\]

\begin{lemma}\label{lem:new-local-trace}
The operator $\mathfrak D_m^\chi[f]$ is trace class.
There exist $c,C_f>0$, independent of the model pair and the
measurable weights $m,n:(0,\infty)\to[0,1]$, such that
\begin{equation}\label{eq:unified-weight-comparison}
\|\mathfrak D_m^\chi[f]-\mathfrak D_n^\chi[f]\|_1
\le C_f\sum_{j\in\Zz}e^{-c|j|}\|(m-n)P_j\|.
\end{equation}
The constants are uniform for bounded holomorphic families
on a common neighborhood of $[0,1]$.
\end{lemma}

\begin{proof}
Choose a positively oriented contour $\Gamma$ enclosing $[0,1]$
whose closed interior lies in the domain of $f$.
For $\mathcal M\in\{T,W\}$, put
$A_m^{\mathcal M}=m\mathsf B^{\mathcal M}m$,
$Q_m^{\mathcal M}(\lambda)=(\lambda\Id-A_m^{\mathcal M})^{-1}$,
and $D_m=A_m^T-A_m^W$.
Since multiplication weights commute with $P_j$,
Lemmas~\ref{lem:new-blocks} and~\ref{lem:new-resolvent-blocks} give
\[
\|P_jQ_m^{\mathcal M}(\lambda)P_l\|
\le C_\Gamma e^{-b|j-l|},
\qquad
\|P_jD_mP_l\|_1\le Ce^{-a|j-l|},
\]
uniformly in the model pair, $m$, and $\lambda\in\Gamma$.
The set $K_\chi=\{j:P_j\chi\ne0\}$ is finite, so the dyadic
expansion of $\chi Q_m^T(\lambda)D_mQ_m^W(\lambda)\chi$
converges absolutely in trace norm.
The resolvent identity and holomorphic functional calculus yield
\[
\mathfrak D_m^\chi[f]
=\frac1{2\pi\ii}\int_\Gamma
f(\lambda)\chi Q_m^T(\lambda)D_mQ_m^W(\lambda)\chi\,d\lambda,
\]
proving the trace-class assertion.

Set $d_j=\|(m-n)P_j\|$,
$E_{\mathcal M}=A_m^{\mathcal M}-A_n^{\mathcal M}$,
and $G=D_m-D_n$.
The blocks of $E_{\mathcal M}$ and $G$ are bounded in operator
and trace norm, respectively, by $C(d_j+d_l)e^{-a|j-l|}$.
Suppressing $\lambda$, the resolvent identity gives
\[
Q_m^TD_mQ_m^W-Q_n^TD_nQ_n^W
=Q_m^TE_TQ_n^TD_mQ_m^W+Q_n^TGQ_m^W
+Q_n^TD_nQ_m^WE_WQ_n^W.
\]
Insert dyadic projections, estimating one block of $D_m$, $G$,
or $D_n$ in trace norm and all other factors in operator norm.
Both end indices lie in $K_\chi$.
For each occurrence of $d_j$, the sum of successive index
distances is therefore at least $2\operatorname{dist}(j,K_\chi)$.
Retaining half the exponential decay for this distance and
summing the remaining geometric series gives
\[
\|\chi(Q_m^TD_mQ_m^W-Q_n^TD_nQ_n^W)\chi\|_1\le C_\Gamma\sum_{j\in\Zz}
e^{-c\operatorname{dist}(j,K_\chi)}d_j.
\]
Since $K_\chi$ is fixed and finite,
$e^{-c\operatorname{dist}(j,K_\chi)}\le C_\chi e^{-c|j|}$.
Integration over $\Gamma$ proves
\eqref{eq:unified-weight-comparison}.
A common contour gives the asserted uniformity in $f$.
\end{proof}

For $\Omega\Subset\Dres$, fix a positively oriented contour $\Sigma$ enclosing $[0,1]$ whose closed interior is disjoint from $\Omega$.
For $z\in\Omega$, put $\eta_z(x)=x(1-x)/(z-x)$ and define
\begin{equation}\label{eq:new-Jrt}
\mathcal J_{r,t}^\chi(z)
:=\Tr\!\left(\chi[\eta_z(M_r\mathsf T_tM_r)
-\eta_z(M_r\mathsf W_tM_r)]\chi\right).
\end{equation}
Define $\mathcal J_{r,\sigma}^\chi(z)$ by the same formula for the weighted one-edge pair, and write $\mathcal J_\sigma^\chi=\mathcal J_{0,\sigma}^\chi$.
Set
\begin{equation}\label{eq:weighted-upper-profile}
\mathcal P_r^\chi(z):=\mathcal J_{r,+}^\chi(z)+\mathcal J_{r,-}^\chi(z).
\end{equation}
Lemma~\ref{lem:new-local-trace} justifies these paired traces and their contour representations on $\Sigma$.

\begin{proposition}\label{prop:new-radial}
For every $\Omega\Subset\Dres$, there are $\delta,C_\Omega>0$ such that, uniformly in $t>0$, $z\in\Omega$, and $\sigma=\pm$,
\[
\begin{aligned}
|\mathcal J_{r,t}^\chi(z)-\mathcal J_{0,t}^\chi(z)|
+|\mathcal J_{r,\sigma}^\chi(z)-\mathcal J_\sigma^\chi(z)|
&\le C_\Omega r^\delta,&0<r\le1,\\
|\mathcal J_{r,t}^\chi(z)|+|\mathcal P_r^\chi(z)|
&\le C_\Omega r^{-\delta},&r\ge1.
\end{aligned}
\]
\end{proposition}

\begin{proof}
Apply \eqref{eq:unified-weight-comparison} with $m=M_r$ and $n=\Id$ in the first case, and $n=0$ in the second.
Here $\mathfrak D_0^\chi[f]=0$.
For $0<\delta<\min\{1,c/\log2\}$, the bounds
\[
\|(M_r-\Id)P_j\|\le\min\{1,2^{j+1}r\},
\qquad
\|M_rP_j\|\le e^{-r2^j}
\]
give the required estimates by geometric summation.
The functions $\eta_z$, $z\in\Omega$, are bounded and holomorphic on a common neighborhood of $[0,1]$, so the constants are uniform in $z$.
\end{proof}

\subsection{Two-edge decoupling on growing dyadic truncations}

We first localize the homogeneous models to a growing dyadic window.
This changes the localized traces by an exponentially small amount and prepares the separation of the two edge contributions.

Choose \(\vartheta\in C^\infty(\Rr;[0,1])\) such that \(\vartheta=1\) on \((-\infty,0]\) and \(\vartheta=0\) on \([1,\infty)\).
For integers \(J\ge3\), set
\[
c_J(u)
=
\vartheta(\log_2u-J)\vartheta(-\log_2u-J).
\]
Then \(c_J=1\) on \([2^{-J},2^J]\) and \(\operatorname{supp}c_J\subset[2^{-J-1},2^{J+1}]\).
Fix \(J_0\ge3\) so that \(c_J\chi=\chi\) for all \(J\ge J_0\), and henceforth take \(J\ge J_0\).
For an operator \(B\), write \(B^J:=c_JBc_J\).
Let \(\mathcal J_{r,t}^{\chi,J}\) and \(\mathcal J_{r,\pm}^{\chi,J}\) denote the corresponding paired traces for the compressed finite-band and one-edge pairs.

\begin{lemma}\label{lem:new-cluster-locality}
For every $\Omega\Subset\Dres$, there are $b_0,C_\Omega>0$ such that, uniformly in $r\ge0$, $J\ge J_0$, $t>0$, $z\in\Omega$, and $\sigma=\pm$,
\begin{equation}\label{eq:new-cluster-locality}
\bigl|\mathcal J_{r,t}^{\chi}(z)-\mathcal J_{r,t}^{\chi,J}(z)\bigr|
\le C_\Omega e^{-b_0J},
\qquad
\bigl|\mathcal J_{r,\sigma}^{\chi}(z)-\mathcal J_{r,\sigma}^{\chi,J}(z)\bigr|
\le C_\Omega e^{-b_0J}.
\end{equation}
\end{lemma}

\begin{proof}
Apply \eqref{eq:unified-weight-comparison} with $m=M_r$, $n=c_JM_r$, and $f=\eta_z$ to each model pair.
The weights agree on the dyadic intervals with $|j|\le J-1$, and their difference has norm at most one elsewhere.
Summing $e^{-c|j|}$ over the remaining indices proves \eqref{eq:new-cluster-locality}, uniformly in $r$ and the model parameters.
\end{proof}

For \(r\ge0\), set \(h_J(u):=c_J(u)e^{-ru}\).
For \(\mathcal M\in\{T,W\}\), define
\[
H_J:=M_{h_J^2},\qquad
L_{\mathcal M}:=M_{h_J}\mathsf A_-^{\mathcal M}M_{h_J},\qquad
R_{\mathcal M}:=M_{h_J}\mathsf A_+^{\mathcal M}M_{h_J}.
\]
For \(t>0\), let \(U_t\) be the unitary multiplication operator on \(L^2(0,\infty)\) defined by
\[
(U_tg)(u):=e^{\ii tu}g(u).
\]
Set \(R_{\mathcal M,t}:=U_tR_{\mathcal M}U_t^*\) and \(G_J(\lambda):=(\lambda\Id-H_J)^{-1}\).
The dependence on \(r\) is suppressed in this notation.
The following estimate controls the interaction between the two truncated one-edge operators.

\begin{lemma}\label{lem:new-cross}
There exists \(C_\Sigma>0\) such that, uniformly in \(r\ge0\), \(J\ge J_0\), \(\lambda\in\Sigma\), \(t\ge1\), and \(\mathcal M,\mathcal M'\in\{T,W\}\),
\[
\|R_{\mathcal M,t}G_JL_{\mathcal M'}\|_1
+\|L_{\mathcal M}G_JR_{\mathcal M',t}\|_1
\le C_\Sigma2^{64J}t^{-2}.
\]
\end{lemma}

\begin{proof}
We treat $R_{\mathcal M,t}G_JL_{\mathcal M'}$; the reversed product has the same proof.
Put $b_{\mathcal M}(u,v)=h_J(u)c_{\mathcal M}(u,v)h_J(v)$, where $c_W=1/(2\pi)$ and $c_T=\sqrt{uv}/[\pi(u+v)]$.
The support of $h_J$ lies in $I_J=[2^{-J-1},2^{J+1}]$.
As $r^ke^{-ru}\le C_ku^{-k}$ and $\|c_J^{(k)}\|_\infty\le C_k2^{k(J+1)}$, we have, uniformly in $r\ge0$,
\[
\|\partial_u^k\partial_v^\ell b_{\mathcal M}\|_\infty
\le C_{k,\ell}2^{(k+\ell)(J+1)}
\quad(k+\ell\le12).
\]
The same estimate holds for derivatives through order eight of $g_J(u)=(\lambda-h_J(u)^2)^{-1}$, with constants depending on $\Sigma$.

Regularize both half-line Fourier integrals by the Abel factor $e^{-\varepsilon(p+q)}$.
The resulting mixed-product kernel is
\[
K_\varepsilon(u,v)
=\int_t^\infty\!\int_0^\infty
 e^{-\varepsilon(p+q)}e^{\ii pu+\ii qv}\left[\int_0^\infty
 b_{\mathcal M}(u,w)g_J(w)b_{\mathcal M'}(w,v)
 e^{-\ii(p+q)w}\,dw\right]dq\,dp.
\]
The frequency sum $p+q$, rather than a difference, is the reason that the two edges separate.
Differentiate at most twice in each external variable and integrate eight times by parts in $w$.
The $w$-amplitude is smooth and compactly supported, so no boundary terms occur.
The preceding estimates give
\[
\|\partial_u^k\partial_v^\ell K_\varepsilon\|_\infty
\le C_\Sigma2^{13(J+1)}
\int_t^\infty\!\int_0^\infty(p+q)^{k+\ell-8}\,dq\,dp
\le C_\Sigma2^{13(J+1)}t^{-2},\qquad 0\le k,\ell\le2.
\]
Rescaling $I_J$ to a fixed rectangle costs at most $C2^{5(J+1)}$ in these derivative bounds.
The double Fourier-series nuclear criterion used in Lemma~\ref{lem:new-blocks} therefore gives
\[
\|K_\varepsilon\|_1
\le C_\Sigma2^{18(J+1)}t^{-2}
\le C_\Sigma2^{64J}t^{-2}.
\]
This estimate is uniform in $\varepsilon>0$ and $r\ge0$.

After the integrations by parts, the same integrable $(p+q)^{k+\ell-8}$ majorant applies to differences of Abel factors.
Dominated convergence in the $C^{2,2}$ kernel norm, followed by the nuclear criterion, gives a limit in $\mathfrak S_1$ as $\varepsilon\downarrow0$.
The regularized one-edge operators converge strongly with uniform operator bounds: for the $W$-model by their Fourier multipliers and for the $T$-model by \eqref{eq:new-dephasing}.
The shift from $p\ge0$ to $p\ge t$ introduces only a scalar factor $e^{-\varepsilon t}$, which tends to one for each fixed $t$.
Hence the trace-norm limit is $R_{\mathcal M,t}G_JL_{\mathcal M'}$.
Reversing the frequency orientations proves the other estimate.
\end{proof}

The cross-edge estimate gives the corresponding separation of the compressed paired traces.

\begin{lemma}\label{lem:new-decouple}
For every \(\Omega\Subset\Dres\), there is \(C_\Omega>0\) such that, uniformly in \(r\ge0\), integers \(J\ge J_0\), \(t\ge1\), and \(z\in\Omega\),
\begin{equation}\label{eq:new-cluster-paired}
\bigl|
\mathcal J_{r,t}^{\chi,J}(z)
-\mathcal J_{r,+}^{\chi,J}(z)
-\mathcal J_{r,-}^{\chi,J}(z)
\bigr|
\le C_\Omega2^{64J}t^{-2}.
\end{equation}
\end{lemma}

\begin{proof}
Write $A_{r,t}^{T,J}:=(M_r\mathsf T_tM_r)^J$ and
$A_{r,t}^{W,J}:=(M_r\mathsf W_tM_r)^J$.
For either model $\mathcal M\in\{T,W\}$, the distributional identity $X_t(v)=2\pi\delta_0(v)-X_-(v)-e^{\ii tv}X_+(v)$, together with $2\pi c_{\mathcal M}(u,u)=1$, gives
\[
A_{r,t}^{\mathcal M,J}=H_J-L_{\mathcal M}-R_{\mathcal M,t},
\quad A_{r,+}^{\mathcal M,J}=H_J-L_{\mathcal M},
\quad A_{r,-}^{\mathcal M,J}=H_J-R_{\mathcal M}.
\]
Let $Q_t,Q_+$ be the resolvents of the first two operators, and let $Q_{-,t}$ be the conjugated resolvent of the third.
Write $G=G_J$ and $Q_{\rm app}=Q_++Q_{-,t}-G$.
Since
\[
Q_+-G=-GL_{\mathcal M}Q_+,
\qquad Q_{-,t}-G=-GR_{\mathcal M,t}Q_{-,t},
\]
direct multiplication by $G^{-1}+L_{\mathcal M}+R_{\mathcal M,t}$ yields
\[
Q_t-Q_{\rm app}
=Q_t\bigl(R_{\mathcal M,t}GL_{\mathcal M}Q_+
+L_{\mathcal M}GR_{\mathcal M,t}Q_{-,t}\bigr).
\]
All resolvents have norm bounded uniformly on $\Sigma$.
Lemma~\ref{lem:new-cross} therefore shows that the difference on the left is trace class, with norm at most $C_\Sigma2^{64J}t^{-2}$.

Subtract the identities for $\mathcal M=T$ and $W$ before taking traces: the common term $G$ cancels.
As $U_t$ commutes with $\chi$, the contour formula for the paired localized differences contains only the trace-class resolvent errors.
Uniform boundedness of $\eta_z$ on $\Omega\times\Sigma$ gives \eqref{eq:new-cluster-paired}.
No separate trace of $G$ or of an individual one-edge resolvent is used.
\end{proof}

\begin{theorem}
\label{thm:new-homogeneous-limit}
For every \(\Omega\Subset\Dres\), there are \(\epsilon,C_\Omega>0\) such that
\[
\sup_{r\ge0}\sup_{z\in\Omega}
\bigl|\mathcal J_{r,t}^{\chi}(z)-\mathcal P_r^{\chi}(z)\bigr|
\le C_\Omega t^{-\epsilon},
\qquad t\ge1.
\]
\end{theorem}

\begin{proof}
Lemmas~\ref{lem:new-cluster-locality} and~\ref{lem:new-decouple} give, uniformly in $r\ge0$ and $z\in\Omega$,
\[
\bigl|\mathcal J_{r,t}^\chi(z)-\mathcal P_r^\chi(z)\bigr|
\le C_\Omega\bigl(e^{-b_0J}+2^{64J}t^{-2}\bigr).
\]
For $t\ge2^{128J_0}$ choose $J=\lfloor(\log_2t)/128\rfloor$.
The right side is $O_\Omega(t^{-b_0/(128\log2)}+t^{-3/2})$.
For the remaining bounded range of $t\ge1$, the uniform localized trace bound from Lemma~\ref{lem:new-local-trace} applies.
Take $\epsilon=\min\{1,b_0/(128\log2)\}>0$.
\end{proof}

\subsection{One-edge coefficient and continuous localization}

The two-edge estimate reduces the large-$t$ analysis to the one-edge traces.
We compute their common coefficient and then integrate over scales.

\begin{lemma}\label{lem:new-normal-density}
Set \(w_\chi:=\int_0^\infty\chi(u)^2\,du/u\).
For \(\sigma\in\{+,-\}\) and \(z\in\Dres\),
\[
\mathcal J_\sigma^\chi(z)
=
\frac{w_\chi}{4\pi^2}
\int_0^1\frac{dx}{z-x}.
\]
\end{lemma}

\begin{proof}
Under the unitary logarithmic-coordinate map $(\mathcal Vg)(x)=e^{x/2}g(e^x)$, followed by Fourier transform, the positive-edge operators become multipliers with symbols
\[
m_T(\alpha)=\frac{1+\tanh(\pi\alpha/2)}2,
\qquad
m_W(\alpha)=\frac{1+\tanh(\pi\alpha)}2.
\]
Indeed, their kernels in logarithmic coordinates are respectively $\frac12\delta_0+\frac{\ii}{2\pi}\operatorname{pv}\operatorname{csch}x$ and
$\frac12\delta_0+\frac{\ii}{4\pi}
\operatorname{pv}\operatorname{csch}(x/2)$;
the Fourier identity
$\operatorname{pv}\int_{\Rr}e^{-\ii\alpha x}(\sinh x)^{-1}\,dx
=-\ii\pi\tanh(\pi\alpha/2)$ gives the displayed multipliers.

Multiplication by $\chi$ becomes multiplication by $\chi(e^x)$, whose squared $L^2$ norm is $w_\chi$.
Since $\eta_z(0)=\eta_z(1)=0$ and both multipliers approach their endpoint values exponentially, $\eta_z(m_T)$ and $\eta_z(m_W)$ are integrable.
Their $L^2$ factorizations justify the localized trace formula
\[
\mathcal J_+^\chi(z)=\frac{w_\chi}{2\pi}
\int_{\Rr}[\eta_z(m_T(\alpha))-
\eta_z(m_W(\alpha))]d\alpha.
\]
Now $m_T'=\pi m_T(1-m_T)$ and $m_W'=2\pi m_W(1-m_W)$.
Changing variables gives
\[
\mathcal J_+^\chi(z)=\frac{w_\chi}{4\pi^2}
\int_0^1\frac{\eta_z(x)}{x(1-x)}\,dx
=\frac{w_\chi}{4\pi^2}\int_0^1\frac{dx}{z-x}.
\]
The negative-edge symbols are obtained by replacing $\alpha$ with $-\alpha$, so they give the same answer.
\end{proof}

The small-$t$ estimate below complements the weight bounds in Proposition~\ref{prop:new-radial}.

\begin{lemma}\label{lem:new-exterior}
For every $\Omega\Subset\Dres$, uniformly in $r\ge0$, $0<t\le1$, and $z\in\Omega$,
\[
|\mathcal J_{r,t}^\chi(z)|\le C_\Omega t.
\]
\end{lemma}

\begin{proof}
Write $\eta_z(A)=A^{1/2}g_z(A)A^{1/2}$ with $g_z(x)=(1-x)/(z-x)$.
For either positive contraction $A=M_r\mathsf T_tM_r$ or $M_r\mathsf W_tM_r$,
\[
\|\chi\eta_z(A)\chi\|_1
\le C_\Omega\Tr(\chi A\chi)
=\frac{C_\Omega t}{2\pi}\int_0^\infty\chi(u)^2e^{-2ru}\,du.
\]
The equality follows from the common kernel diagonal, proving the claim.
\end{proof}

The following continuous localization identity will convert the paired traces into a global scale integral.

\begin{lemma}\label{lem:continuous-scale-resolution}
Set $\chi_r(\xi)=\chi(\xi/r)$.
For every trace-class operator $K$ on $L^2(0,\infty)$,
\[
\int_0^\infty\|\chi_rK\chi_r\|_1\,\frac{dr}{r}
\le w_\chi\|K\|_1,
\qquad
\int_0^\infty\Tr(\chi_rK\chi_r)\,\frac{dr}{r}
=w_\chi\Tr K.
\]
\end{lemma}

\begin{proof}
The identity $\int_0^\infty\chi(\xi/r)^2\,dr/r=w_\chi$ gives both assertions for a rank-one operator $K=a\otimes b^*$, using Cauchy--Schwarz for the trace-norm bound.
An absolutely convergent singular-value expansion extends them to $\mathfrak S_1$.
\end{proof}

\section{Endpoint finite part and strong asymptotics}
\label{sec:finite-constant}

We integrate the localized estimates of Section~\ref{sec:hardedge-halfline} to identify the logarithmic endpoint correction and its finite part.
Combined with the canonical asymptotic \eqref{eq:canonical-main-global-halfline}, this completes the proof of Theorem~\ref{thm:strong-halfline}.

\subsection{The scale integral}

We retain the notation \(T_L:=T_L^\infty\) and \(W_L:=\widetilde W_L^\infty\) from Section~\ref{sec:hardedge-halfline}.
Fix a nonzero real-valued \(\chi\in C_c^\infty((0,\infty))\) and write
\[
w_\chi:=\int_0^\infty\chi(u)^2\,\frac{du}{u}.
\]
For \(f\) holomorphic in a neighborhood of \([0,1]\), with \(f(0)=f(1)=0\), define, for \(r\ge0\) and \(t>0\),
\begin{equation}\label{eq:fc-scale-profiles}
\begin{aligned}
\mathcal J_{r,t}^\chi[f]
&:=\Tr\!\bigl(
\chi\,[f(M_r\mathsf T_tM_r)-f(M_r\mathsf W_tM_r)]\,\chi
\bigr),\\
\mathcal P_r^\chi[f]
&:=\sum_{\sigma=\pm}
\Tr\!\bigl(
\chi\,[f(A_{r,\sigma}^T)-f(A_{r,\sigma}^W)]\,\chi
\bigr).
\end{aligned}
\end{equation}
For \(f=\eta_z\), these definitions recover \(\mathcal J_{r,t}^\chi(z)\) and \(\mathcal P_r^\chi(z)\) from \eqref{eq:new-Jrt} and \eqref{eq:weighted-upper-profile}, respectively, by Lemma~\ref{lem:new-local-trace}.
Set
\[
\mathfrak c(f)
:=
\frac{1}{2\pi^2}
\int_0^1\frac{f(x)}{x(1-x)}\,dx,
\qquad
c_\chi(f):=w_\chi\mathfrak c(f).
\]
The quotient \(f(z)/[z(1-z)]\) extends holomorphically across \(z=0\) and \(z=1\).
The functional $\mathfrak c(f)$ is $1/(2\pi^2)$ times the usual jump functional in Wiener--Hopf trace asymptotics; see \cite[Eq.~(2.12)]{Sobolev2017}.
Theorem~\ref{thm:fc-direct-scale} identifies this coefficient and the finite part of the trace difference considered here.

\begin{lemma}\label{lem:fc-weighted-transitions}
The localized differences in \eqref{eq:fc-scale-profiles} are trace class.
There exist \(\delta,\delta_+,\epsilon>0\) such that
\[
\begin{alignedat}{2}
\sup_{r\ge0}|\mathcal J_{r,t}^\chi[f]|
&\le Ct,
&\qquad& 0<t\le1,\\
\sup_{t>0}
|\mathcal J_{r,t}^\chi[f]-\mathcal J_{0,t}^\chi[f]|
+|\mathcal P_r^\chi[f]-c_\chi(f)|
&\le Cr^\delta,
&& 0<r\le1,\\
\sup_{t>0}|\mathcal J_{r,t}^\chi[f]|
+|\mathcal P_r^\chi[f]|
&\le Cr^{-\delta_+},
&& r\ge1,\\
\sup_{r\ge0}
|\mathcal J_{r,t}^\chi[f]-\mathcal P_r^\chi[f]|
&\le Ct^{-\epsilon},
&& t\ge1.
\end{alignedat}
\]
The bounds are uniform for bounded holomorphic families on a common neighborhood of \([0,1]\).
\end{lemma}

\begin{proof}
Choose a positively oriented contour \(\Gamma\) around \([0,1]\) on which \(g_f(z)=f(z)/[z(1-z)]\) is holomorphic.
Since
\[
f(x)=\frac{1}{2\pi\ii}\int_\Gamma g_f(z)\eta_z(x)\,dz,
\]
Lemma~\ref{lem:new-local-trace} and holomorphic functional calculus give, with trace-class differences,
\[
\mathcal J_{r,t}^\chi[f]
=\frac{1}{2\pi\ii}\int_\Gamma g_f(z)\mathcal J_{r,t}^\chi(z)\,dz,
\qquad
\mathcal P_r^\chi[f]
=\frac{1}{2\pi\ii}\int_\Gamma g_f(z)\mathcal P_r^\chi(z)\,dz.
\]
Lemma~\ref{lem:new-exterior} gives the first estimate, Proposition~\ref{prop:new-radial} the second and third, and Theorem~\ref{thm:new-homogeneous-limit} the fourth.
Here $\mathcal P_0^\chi[f]=c_\chi(f)$ follows from Lemma~\ref{lem:new-normal-density}.
For a bounded holomorphic family on a common neighborhood of \([0,1]\), one may fix \(\Gamma\) and bound \(g_f\) there uniformly.
The same contour therefore gives uniform constants in all four estimates.
\end{proof}

Define
\begin{equation}\label{eq:fc-crossover-integrals}
\begin{aligned}
C_{\rm low}^\chi(f)
&:=
\int_0^1\mathcal J_{0,t}^\chi[f]\,\frac{dt}{t}
+\int_1^\infty
\bigl(\mathcal J_{0,t}^\chi[f]-c_\chi(f)\bigr)\,\frac{dt}{t},\\
C_{\rm up}^\chi(f)
&:=
\int_0^1
\bigl(\mathcal P_r^\chi[f]-c_\chi(f)\bigr)\,\frac{dr}{r}
+\int_1^\infty\mathcal P_r^\chi[f]\,\frac{dr}{r}.
\end{aligned}
\end{equation}
By Lemma~\ref{lem:fc-weighted-transitions}, all four integrals converge absolutely, uniformly for bounded holomorphic families on a common neighborhood of \([0,1]\).
Set
\[
\mathcal B_\chi(f)
:=
\frac{C_{\rm low}^\chi(f)+C_{\rm up}^\chi(f)}{w_\chi}.
\]

\begin{theorem}[Continuous-scale finite part]
\label{thm:fc-direct-scale}
The functional \(\mathcal B_\chi\) is independent of \(\chi\); write \(\mathcal B\) for the common functional.
For every \(f\) holomorphic near \([0,1]\) with \(f(0)=f(1)=0\), there exists \(\kappa>0\) such that
\[
\Tr[f(T_L)-f(W_L)]
=
\mathfrak c(f)\log L+\mathcal B(f)+O(L^{-\kappa}).
\]
For bounded holomorphic families on a common neighborhood of \([0,1]\), \(\kappa\) and the implied constant may be chosen uniformly.
\end{theorem}

\begin{proof}
Both \(f(T_L)\) and \(f(W_L)\) are trace class because \(f(0)=0\).
For \(\chi_r(u)=\chi(u/r)\), Lemma~\ref{lem:continuous-scale-resolution}, together with dilation \(u\mapsto ru\), gives the exact, absolutely convergent identity
\begin{equation}\label{eq:fc-continuous-resolution}
\Tr[f(T_L)-f(W_L)]
=\frac1{w_\chi}\int_0^\infty
\mathcal J_{r,Lr}^\chi[f] \,\frac{dr}{r}.
\end{equation}
For \(L^{-1}\le\rho\le1\), use the small-\(r\) estimate below \(\rho\), and the weighted large-\(t\) estimate above \(\rho\).
After the change of variables \(t=Lr\) in the first part, this yields
\[
\int_0^\infty\mathcal J_{r,Lr}^\chi[f]\,\frac{dr}{r}
=\int_0^{L\rho}\mathcal J_{0,t}^\chi[f]\,\frac{dt}{t}
+\int_\rho^\infty\mathcal P_r^\chi[f]\,\frac{dr}{r}
+O\bigl(\rho^\delta+(L\rho)^{-\epsilon}\bigr).
\]
By \eqref{eq:fc-crossover-integrals}, the sum of the two integrals is
\[
c_\chi(f)\log L+C_{\rm low}^\chi(f)+C_{\rm up}^\chi(f)
-\int_{L\rho}^\infty
  (\mathcal J_{0,t}^\chi[f]-c_\chi(f))\,\frac{dt}{t}
-\int_0^\rho(\mathcal P_r^\chi[f]-c_\chi(f))\,\frac{dr}{r}.
\]
By Lemma~\ref{lem:fc-weighted-transitions}, the two tails are \(O((L\rho)^{-\epsilon})\) and \(O(\rho^\delta)\), respectively.
Choose \(\rho=L^{-\epsilon/(\delta+\epsilon)}\).
Both errors are then \(O(L^{-\kappa})\), where \(\kappa=\delta\epsilon/(\delta+\epsilon)>0\).
The result follows from \eqref{eq:fc-continuous-resolution} and \(c_\chi(f)=w_\chi\mathfrak c(f)\).

For two choices \(\chi\) and \(\psi\), the trace and its logarithmic coefficient are unchanged.
Subtract their expansions and let \(L\to\infty\) to obtain \(\mathcal B_\chi(f)=\mathcal B_\psi(f)\).
Uniformity in \(f\) follows from the common-contour estimates in Lemma~\ref{lem:fc-weighted-transitions}.
\end{proof}

\subsection{The determinant and resolvent finite parts}

For $\theta\in S_\pi$, define
\begin{equation}\label{eq:FL-hardedge-def-halfline}
F_L(\theta)
:=\log\det(\Id+(e^\theta-1)T_L)
-\log\det(\Id+(e^\theta-1)W_L),
\end{equation}
using the holomorphic logarithms normalized at zero, and set
\begin{equation}\label{eq:fc-HL}
\Phi_L(\theta)
:=
F_L(\theta)-\frac{\theta^2}{4\pi^2}\log L,
\qquad \theta\in S_\pi.
\end{equation}

\begin{corollary}\label{cor:fc-Phi}
For \(\theta\in S_\pi\), define
\begin{equation}\label{eq:fc-Phi-direct}
\Phi_\infty(\theta)
:=
\mathcal B(f_\theta),
\qquad
f_\theta(x)
:=
\Logp\!\bigl(1+(e^\theta-1)x\bigr)-\theta x.
\end{equation}
Then \(\Phi_\infty\) is holomorphic on \(S_\pi\), \(\Phi_\infty(0)=\Phi_\infty'(0)=0\), and for every \(K\Subset S_\pi\) there exists \(\kappa_K>0\) such that
\begin{equation}\label{eq:fc-Phi-rate}
\sup_{\theta\in K}
|\Phi_L(\theta)-\Phi_\infty(\theta)|
=
O_K(L^{-\kappa_K}).
\end{equation}
\end{corollary}

\begin{proof}
For \(K\Subset S_\pi\), the segments \(1+(e^\theta-1)[0,1]\), \(\theta\in K\), form a compact set avoiding \((-\infty,0]\).
Thus the \(f_\theta\) are bounded and holomorphic on a common neighborhood of \([0,1]\), with \(f_\theta(0)=f_\theta(1)=0\).
Since \(\Tr T_L=\Tr W_L\), we have the exact identity \(F_L(\theta)=\Tr[f_\theta(T_L)-f_\theta(W_L)]\).
Moreover,
\[
\frac{d^2}{d\theta^2}\mathfrak c(f_\theta)
=\frac1{2\pi^2}\int_0^1
\frac{e^\theta\,dx}{[1+(e^\theta-1)x]^2}
=\frac1{2\pi^2}.
\]
Since \(f_0=\partial_\theta f_\theta|_{\theta=0}=0\), we have \(\mathfrak c(f_\theta)=\theta^2/(4\pi^2)\).
Theorem~\ref{thm:fc-direct-scale} now gives \eqref{eq:fc-Phi-rate}, uniformly for \(\theta\in K\).
Each \(\Phi_L\) is holomorphic on \(S_\pi\), so the locally uniform limit \(\Phi_\infty\) is holomorphic.
Finally, \(\Phi_L(0)=\Phi_L'(0)=0\) for all \(L\); Cauchy's formula passes these identities to the limit.
\end{proof}

\begin{proof}[Proof of Theorem~\ref{thm:strong-halfline}]
By Lemma~\ref{lem:reduction-halfline} and \eqref{eq:FL-hardedge-def-halfline}, the normalized holomorphic logarithms satisfy
\[
\log\Ee\bigl[e^{\theta N_\infty(L)}\bigr]
=\log\det\bigl(\Id+(e^\theta-1)\widetilde W_L^\infty\bigr)
+F_L(\theta).
\]
The canonical expansion \eqref{eq:canonical-main-global-halfline} contributes
\(L\Lambda_\infty(\theta)+\theta^2\log L/(4\pi^2)
+A_{\rm FH}(\theta)+o(1)\).
Corollary~\ref{cor:fc-Phi} gives
\[
F_L(\theta)=\frac{\theta^2}{4\pi^2}\log L
+\Phi_\infty(\theta)+O_K(L^{-\kappa_K}),
\qquad \theta\in K\Subset S_\pi.
\]
Adding the two expressions proves \eqref{eq:main-strong-halfline} with \(A_\infty=A_{\rm FH}+\Phi_\infty\), holomorphic on \(S_\pi\).
The canonical remainder is \(o(1)\), while the endpoint remainder has a power bound; together they are locally uniformly \(o(1)\).
Uniqueness follows by subtracting two such expansions and letting \(L\to\infty\).
\end{proof}

\begin{proof}[Proof of Corollary~\ref{cor:general-halfline-scaling}]
Exact scaling \eqref{eq:count-scaling-intro-halfline} gives
\(\Lambda_a^\infty=a^{-1}\Lambda_\infty\) on \(S_\pi\), and Legendre duality yields
\(J_a^\infty(x)=a^{-1}J_\infty(ax)\).
Theorem~\ref{thm:LDP-halfline}, applied at length \(L/a\), gives an LDP for \(N_a^\infty(L)/L\) with speed \(L/a\) and rate \(x\mapsto J_\infty(ax)\), equivalently with speed \(L\) and rate \(J_a^\infty\).
Finally, substituting \(L/a\) for \(L\) in \eqref{eq:main-strong-halfline} and using \(\log(L/a)=\log L-\log a\) proves \eqref{eq:halfline-scaling-cor-strong} with
\(A_{\infty,a}(\theta)=A_\infty(\theta)-\theta^2\log a/(2\pi^2)\).
The remainder is locally uniform on \(S_\pi\) since \(a>0\) is fixed.
\end{proof}

Define
\[
R_L(z)
:=
\Tr[\eta_z(T_L)-\eta_z(W_L)]
-\frac{\log L}{2\pi^2}\int_0^1\frac{dx}{z-x},
\qquad z\in\Dres.
\]

\begin{corollary}\label{thm:fc-resolvent}
The function \(R_\infty(z):=\mathcal B(\eta_z)\) is holomorphic on \(\Dres\).
For every \(\Omega\Subset\Dres\), there exists \(\kappa_\Omega>0\) such that
\begin{equation}\label{eq:fc-Rconv}
\sup_{z\in\Omega}|R_L(z)-R_\infty(z)|
=
O_\Omega(L^{-\kappa_\Omega}).
\end{equation}
With \(\theta(z):=\Logp(1-z^{-1})\),
\begin{equation}\label{eq:fc-Rphi}
R_\infty(z)=-\Phi_\infty'(\theta(z)).
\end{equation}
If \(f\) is holomorphic in a neighborhood of \([0,1]\) and \(f(0)=f(1)=0\), then, for every positively oriented piecewise \(C^1\) simple closed contour \(\Gamma\) whose interior contains \([0,1]\) and whose closed interior is contained in that neighborhood,
\begin{equation}\label{eq:fc-Rfun}
\mathcal B(f)
=
\frac{1}{2\pi\ii}
\int_\Gamma
\frac{f(z)}{z(1-z)}R_\infty(z)\,dz.
\end{equation}
\end{corollary}

\begin{proof}
For \(\Omega\Subset\Dres\), the functions \(\eta_z\), \(z\in\Omega\), form a bounded holomorphic family on a common neighborhood of \([0,1]\).
Theorem~\ref{thm:fc-direct-scale} gives \eqref{eq:fc-Rconv}; since each \(R_L\) is holomorphic, its locally uniform limit \(R_\infty\) is holomorphic.

For \(z\in\Dres\), put \(\theta(z)=\Logp(1-z^{-1})\).
Then \(\theta(z)\in S_\pi\setminus\{0\}\), \(\int_0^1(z-x)^{-1}\,dx=-\theta(z)\), and
\[
\eta_z(x)
=x-\frac{e^{\theta(z)}x}{1+(e^{\theta(z)}-1)x}.
\]
Differentiate the normalized determinant difference \eqref{eq:FL-hardedge-def-halfline} and use \(\Tr T_L=\Tr W_L\).
Together with \eqref{eq:fc-HL}, these identities give \(R_L(z)=-\Phi_L'(\theta(z))\) exactly.
By Cauchy's formula, the locally uniform convergence \eqref{eq:fc-Phi-rate} also holds for derivatives, proving \eqref{eq:fc-Rphi}.
Finally, holomorphic functional calculus gives
\[
\Tr[f(T_L)-f(W_L)]-\mathfrak c(f)\log L
=\frac1{2\pi\ii}\int_\Gamma
\frac{f(z)}{z(1-z)}R_L(z)\,dz.
\]
Passing to the limit by \eqref{eq:fc-Rconv} and Theorem~\ref{thm:fc-direct-scale} proves \eqref{eq:fc-Rfun}.
\end{proof}

\section{Trace coefficients and normalization}\label{sec:all-order}

We identify the renormalized trace coefficients, evaluate the finite-length variance and quadratic constant, and describe the higher-order terms and cumulants.

\subsection{Fixed-order coefficients}

Set \(h_0:=0\) and \(h_m:=\sum_{j=1}^m j^{-1}\) for \(m\ge1\).
For \(n\ge1\), define
\[
\mu_n(L)
:=
\Tr(T_L^n-W_L^n)
+\frac{h_{n-1}}{2\pi^2}\log L.
\]
Since \(\Tr T_L=\Tr W_L\), one has \(\mu_1(L)=0\).

\begin{lemma}\label{lem:fc-coeff}
For every \(n\ge1\),
\[
\mu_{n,\infty}
:=
\mathcal B(x^n-x)
=
\lim_{L\to\infty}\mu_n(L).
\]
Moreover, for every \(0<r<1\), there exists \(C_r<\infty\) such that, for all sufficiently large \(L\), every \(0<s<r\), and every integer \(N\ge0\),
\[
\sup_{|q|\le s}
\left|
\sum_{n>N}\frac{(-1)^{n-1}}{n}\mu_n(L)q^n
\right|
\le
C_r\frac{(s/r)^{N+1}}{1-s/r},
\]
and the same bound holds with \(\mu_{n,\infty}\) in place of \(\mu_n(L)\).
\end{lemma}

\begin{proof}
For $n\ge2$,
\(
(x^n-x)/(x(1-x))=-\sum_{j=0}^{n-2}x^j
\), so
\(\mathfrak c(x^n-x)=-h_{n-1}/(2\pi^2)\).
Theorem~\ref{thm:fc-direct-scale}, together with \(\Tr T_L=\Tr W_L\), gives the asserted limit; the case $n=1$ is identically zero.

For $|q|<1$, let $\Psi_L(q)=\Phi_L(\Logp(1+q))$.
The two trace-logarithm expansions and
\(
\Logp^2(1+q)
=2\sum_{n\ge2}(-1)^nh_{n-1}q^n/n
\)
give
\begin{equation}\label{eq:fc-Jseries}
\Psi_L(q)
=
\sum_{n\ge1}
\frac{(-1)^{n-1}}{n}\mu_n(L)q^n.
\end{equation}
Fix $0<r<1$.
The image of $|q|\le r$ under $q\mapsto\Logp(1+q)$ is compact in $S_\pi$.
Corollary~\ref{cor:fc-Phi} therefore bounds $\Psi_L$ uniformly on this disk for all sufficiently large $L$, and gives locally uniform convergence to $\Phi_\infty(\Logp(1+q))$.
Cauchy's coefficient estimate yields
\[
|\mu_n(L)|+|\mu_{n,\infty}|\le nC_rr^{-n},\qquad n\ge1,
\]
with the limiting coefficients identified by the first part.
Summing over $n>N$ on $|q|\le s<r$ proves both tail bounds.
\end{proof}

The trace powers also admit exact finite-dimensional integral representations.
For \(n\ge2\), set
\[
\Delta_{n-1}^{\circ}
:=
\left\{
\boldsymbol\alpha:\alpha_j>0,\ 
\sum_{j<n}\alpha_j<1
\right\},
\qquad
\alpha_n:=1-\sum_{j<n}\alpha_j,
\qquad
\alpha_{n+1}:=\alpha_1,
\]
and define
\[
p_n(\boldsymbol\alpha)
:=
\prod_{j=1}^n
\frac{\alpha_j}{\alpha_j+\alpha_{j+1}},
\qquad
b_n(\boldsymbol\alpha,\boldsymbol\tau)
:=
\sum_{j=1}^n
(\alpha_j-\alpha_{j+1})\tau_j.
\]

\begin{lemma}\label{lem:fc-parameter}
For every fixed \(n\ge2\),
\begin{align*}
\Tr(T_L^n)
&=
\frac{(n-1)!L^n}{\pi^n}
\int_{\Delta_{n-1}^{\circ}\times[0,1]^n}
\frac{p_n(\boldsymbol\alpha)}
     {(2-\ii Lb_n(\boldsymbol\alpha,\boldsymbol\tau))^n}\,
d\boldsymbol\alpha\,d\boldsymbol\tau,
\\
\Tr(W_L^n)
&=
\frac{(n-1)!L^n}{\pi^n}
\int_{\Delta_{n-1}^{\circ}\times[0,1]^n}
\frac{2^{-n}}
     {(2-\ii Lb_n(\boldsymbol\alpha,\boldsymbol\tau))^n}\,
d\boldsymbol\alpha\,d\boldsymbol\tau.
\end{align*}
The integrals are absolutely convergent and real.
\end{lemma}

\begin{proof}
The cyclic kernel products are absolutely integrable, since $|X_L|\le L$ and $\prod_{j=1}^n\xi_j/(\xi_j+\xi_{j+1})\le2^{-n}$.
Insert $X_L(v)=L\int_0^1e^{\ii Lv\tau}\,d\tau$ and set $s=\sum_j\xi_j$, $\xi_j=s\alpha_j$.
The Jacobian is $s^{n-1}$, and the exponential factor becomes $e^{-(2-\ii Lb_n)s}$.
Thus
\[
\int_0^\infty s^{n-1}e^{-(2-\ii Lb_n)s}\,ds
=\frac{(n-1)!}{(2-\ii Lb_n)^n},
\]
which gives both formulas.
Their integrals are absolutely convergent because $p_n\le2^{-n}$ and $|2-\ii Lb_n|\ge2$.
Finally, $\tau_j\mapsto1-\tau_j$ sends $b_n$ to $-b_n$ and conjugates the integrand.
Hence both integrals are real.
\end{proof}

\subsection{Quadratic normalization}

The exact variance determines the quadratic constants.

\begin{proposition}\label{prop:exact-variance-halfline}
For $a,L>0$, put $T=L/a$.
Then
\begin{equation}\label{eq:exact-variance-halfline}
\operatorname{Var}N_a^\infty(L)
=\frac{T}{4\pi}
-\frac1{2\pi^2}
\left[
\left(\frac T2-\frac2T\right)\arctan\frac T2
-\log\left(1+\frac{T^2}{4}\right)+1
\right].
\end{equation}
Moreover, $T\mapsto\operatorname{Var}N_\infty(T)/\Ee N_\infty(T)$ is strictly decreasing, with limits $1$ as $T\downarrow0$ and $1/2$ as $T\to\infty$.
\end{proposition}

\begin{proof}
By dilation, it suffices to take $a=1$ and $L=T$.
The determinantal variance formula gives $\operatorname{Var}N_\infty(T)=\Tr K_T^\infty-\Tr((K_T^\infty)^2)$.
Writing $y$ for the sum of the two heights and integrating over the horizontal coordinates gives
\[
\Tr\bigl((K_T^\infty)^2\bigr)
=\frac{T}{\pi^2}\int_2^\infty
\frac{(y-2)\arctan(T/y)}{y^3}\dd y
=\frac{T}{2\pi^2}\int_0^1(1-v)\arctan(Tv/2)\dd v.
\]
Evaluation of the last integral and $\Tr K_T^\infty=T/(4\pi)$ prove \eqref{eq:exact-variance-halfline}.
The variance-to-mean ratio is therefore $1-(2/\pi)\int_0^1(1-v)\arctan(Tv/2)\dd v$.
The integral is strictly increasing in $T$, and dominated convergence gives the stated limits.
\end{proof}

\begin{corollary}\label{cor:fc-series}
For \(\theta\in S_\pi\) with \(|e^\theta-1|<1\),
\begin{equation}\label{eq:fc-PhiSeries}
\Phi_\infty(\theta)
=
\sum_{n=1}^\infty
\frac{(-1)^{n-1}}{n}
\mu_{n,\infty}(e^\theta-1)^n.
\end{equation}
Moreover,
\[
\mu_{2,\infty}
=
\frac{1+\log 2}{2\pi^2},
\qquad
\Phi_\infty''(0)
=
-\frac{1+\log 2}{2\pi^2}.
\]
Consequently,
\[
A_\infty''(0)
=
-\frac{\log 2}{\pi^2},
\qquad
\operatorname{Var}(N_\infty(L))
=
\frac{L}{8\pi}
+\frac{\log L}{\pi^2}
-\frac{\log 2}{\pi^2}
+o(1).
\]
\end{corollary}

\begin{proof}
Lemma~\ref{lem:fc-coeff} allows passage to the limit in \eqref{eq:fc-Jseries} on compact subsets of $|q|<1$, giving \eqref{eq:fc-PhiSeries}.
For $a=1$, expanding \eqref{eq:exact-variance-halfline} as $L\to\infty$ gives the stated variance asymptotic.
Cauchy's formula applied to \eqref{eq:main-strong-halfline} therefore yields $A_\infty''(0)=-\log2/\pi^2$.
Since $A_\infty=A_{\rm FH}+\Phi_\infty$, Proposition~\ref{prop:explicit-AFH-halfline} gives $\Phi_\infty''(0)=-(1+\log2)/(2\pi^2)$.
Finally, \eqref{eq:fc-PhiSeries} and $\mu_{1,\infty}=0$ give $\mu_{2,\infty}=-\Phi_\infty''(0)$.
\end{proof}

For the quadratic test function, both crossover contributions can be evaluated explicitly.

\begin{proposition}
Let \(f_2(x):=x(1-x)\), and let \(C_{\rm low}^\chi(f_2)\) and \(C_{\rm up}^\chi(f_2)\) be defined by \eqref{eq:fc-crossover-integrals}.
Set \(m_\chi:=\int_0^\infty\chi(u)^2\log u\,du/u\).
Then
\[
C_{\rm low}^\chi(f_2)
=
\frac{\gamma_{\rm E}w_\chi+m_\chi}{2\pi^2},
\qquad
C_{\rm up}^\chi(f_2)
=
-\frac{\gamma_{\rm E}w_\chi+m_\chi+(1+\log 2)w_\chi}{2\pi^2}.
\]
Consequently, \(\mathcal B(f_2)=\Phi_\infty''(0)=-(1+\log 2)/(2\pi^2)\).
\end{proposition}

\begin{proof}
The linear traces cancel in $f_2(A)=A-A^2$.
Directly from the kernels,
\[
\begin{aligned}
\mathcal J_{r,t}^{\chi}[f_2]
&=\frac{1}{2\pi^2}\int_0^\infty\chi(u)^2
\int_0^\infty
\frac{e^{-2r(u+v)}[1-\cos(t(u-v))]}{(u+v)^2}\,dv\,du,\\
\mathcal P_r^{\chi}[f_2]
&=\frac{1}{2\pi^2}\int_0^\infty\chi(u)^2
\int_0^\infty\frac{e^{-2r(u+v)}}{(u+v)^2}\,dv\,du.
\end{aligned}
\]
The second expression is the large-$t$ limit of the first by Theorem~\ref{thm:new-homogeneous-limit} and the Riemann--Lebesgue lemma.
In particular, $c_\chi(f_2)=w_\chi/(2\pi^2)$.

For the lower crossover, put $r=0$ and $v=uy$.
Combine the two $t$-integrals in \eqref{eq:fc-crossover-integrals} at a finite upper limit $R$ before taking $R\to\infty$.
For $a\ne0$,
\[
\int_0^R\frac{1-\cos(at)}t\,dt-\log R
\longrightarrow\gamma_{\rm E}+\log|a|.
\]
The bound $C+|\log|a||$ permits integration against $\chi(u)^2(1+y)^{-2}$ with $a=u(1-y)$.
Using $\int_0^\infty\log|1-y|(1+y)^{-2}\,dy=0$ gives $C_{\rm low}^\chi(f_2)=(\gamma_{\rm E}w_\chi+m_\chi)/(2\pi^2)$.

For the upper crossover, combine the two $r$-integrals before interchanging integrations and use, for $a>0$,
\[
\int_0^1\frac{e^{-ar}-1}{r}\,dr
+\int_1^\infty\frac{e^{-ar}}r\,dr
=-\gamma_{\rm E}-\log a.
\]
The logarithmic moment is integrable because $\chi$ is supported away from zero.
With $a=2(u+v)$ and
\(
\int_0^\infty\log(2(u+v))(u+v)^{-2}\,dv
=[\log(2u)+1]/u
\), this gives the stated upper term.
Adding both terms cancels $m_\chi$ and $\gamma_{\rm E}w_\chi$, leaving $\mathcal B(f_2)=-(1+\log2)/(2\pi^2)$.
\end{proof}

\subsection{Higher-order representation}

To isolate the contribution beyond quadratic order, define, for
\(q\in\Qslit\setminus\{0\}\),
\[
\widehat{\mathcal K}_\infty(q)
:=
\frac{-R_\infty(-1/q)+\mu_{2,\infty}\Logp(1+q)}
     {q(1+q)}.
\]

\begin{proposition}
The function \(\widehat{\mathcal K}_\infty\) extends holomorphically
to \(\Qslit\), with \(\widehat{\mathcal K}_\infty(0)=0\).
For \(|q|<1\),
\[
\widehat{\mathcal K}_\infty(q)
=
\sum_{m\ge1}(-1)^{m+1}
\bigl[\mu_{m+2,\infty}-h_{m+1}\mu_{2,\infty}\bigr]q^m.
\]
For every \(\theta\in S_\pi\),
\[
\Phi_\infty(\theta)
=
-\frac{1+\log2}{4\pi^2}\theta^2
+\int_0^{e^\theta-1}
q\widehat{\mathcal K}_\infty(q)\,dq,
\]
with the integral independent of the path in \(\Qslit\).
Moreover, \(\widehat{\mathcal K}_\infty(q)=-\mathcal B(\psi_q)\)
for \(q\in\Qslit\), where \(\psi_0:=0\) and
\[
\psi_q(x)
:=
\frac{x(1-x)}{1+q}
\left[
\frac{\Logp(1+q)}q-\frac1{1+qx}
\right],
\qquad q\ne0.
\]
\end{proposition}

\begin{proof}
Define the holomorphic function
\[
H(q):=\Phi_\infty(\Logp(1+q))
+\frac{\mu_{2,\infty}}2\Logp^2(1+q),
\qquad q\in\Qslit.
\]
By \eqref{eq:fc-Rphi},
\(H'(q)=q\widehat{\mathcal K}_\infty(q)\) for \(q\ne0\).
Equation~\eqref{eq:fc-PhiSeries} and \(\mu_{1,\infty}=0\)
give \(H(q)=O(q^3)\) at zero, proving the asserted extension
and its value there.
Expanding \(H'(q)/q\) using \eqref{eq:fc-PhiSeries} and
\(\Logp(1+q)/(q(1+q))
=\sum_{m\ge0}(-1)^mh_{m+1}q^m\)
gives the stated series.
Since \(H(0)=0\), integration gives the formula for \(\Phi_\infty\)
and path independence, using \(\Logp(e^\theta)=\theta\) and
\(\mu_{2,\infty}=(1+\log2)/(2\pi^2)\).

For \(q\in\Qslit\setminus\{0\}\), the function \(\psi_q\)
is holomorphic near \([0,1]\), vanishes at both endpoints, and satisfies
\[
q(1+q)\psi_q
=\Logp(1+q)f_2+\eta_{-1/q}.
\]
Applying \(\mathcal B\), with
\(\mathcal B(f_2)=-\mu_{2,\infty}\) and
\(\mathcal B(\eta_{-1/q})=R_\infty(-1/q)\), proves the final identity.
Both sides vanish at \(q=0\).
\end{proof}

The symmetry of the limiting pressure determines which higher cumulants can grow linearly.
We write $\kappa_r(X)$ for the ordinary cumulant of order $r$ of a random variable $X$.

\begin{corollary}\label{cor:even-cumulants-halfline}
The limiting pressure satisfies
\[
\Lambda_\infty(\theta)+\Lambda_\infty(-\theta)
=\frac{\theta^2}{8\pi},\qquad \theta\in S_\pi.
\]
For every fixed $a>0$ and integer $m\ge2$,
\[
\kappa_{2m}\bigl(N_a^\infty(L)\bigr)
=A_\infty^{(2m)}(0)+o(1),
\qquad L\to\infty.
\]
Moreover, $\kappa_3(N_a^\infty(L))=L/(24\pi a)+A_\infty'''(0)+o(1)$.
\end{corollary}

\begin{proof}
Cauchy's formula applied to \eqref{eq:halfline-scaling-cor-strong} reduces the claims to computing $\Lambda_\infty^{(2m)}(0)$ and $\Lambda_\infty'''(0)$, since the logarithmic term and $A_{\infty,a}-A_\infty$ are quadratic in $\theta$.
Equation~\eqref{eq:Lambda-derivatives-halfline} gives $\Lambda_\infty'''(0)=1/(24\pi)$ and $\Lambda_\infty'(\theta)-\Lambda_\infty'(-\theta)=\theta/(4\pi)$ for real $\theta$.
Integration from zero, using $\Lambda_\infty(0)=0$, and the identity theorem yield $\Lambda_\infty(\theta)+\Lambda_\infty(-\theta)=\theta^2/(8\pi)$ on $S_\pi$.
Differentiation at zero gives $\Lambda_\infty^{(2m)}(0)=0$ for $m\ge2$.
\end{proof}

\section{Entropy and precise deviations}\label{sec:entropy-precise-halfline}

In this section, we prove Theorem~\ref{thm:precise-halfline} and give a Bernoulli relative-entropy representation of the rate function.
For $p\in[0,1]$, let $\mathrm{Ber}(p)$ denote the Bernoulli distribution assigning probabilities $p$ and $1-p$ to $1$ and $0$, respectively.
For $p\in(0,1)$ and $q\in[0,1]$, define the Bernoulli relative entropy (Kullback--Leibler divergence) by
\[
\operatorname{KL}(\mathrm{Ber}(q)\Vert\mathrm{Ber}(p))
:=q\log\frac{q}{p}
+(1-q)\log\frac{1-q}{1-p},
\]
where terms with zero coefficients are interpreted as zero.
This relative entropy is nonnegative and vanishes if and only if $q=p$.

\begin{proposition}\label{prop:entropy-halfline}
For $x\ge0$, let $\mathcal Q_x$ be the measurable functions $q:(0,\infty)\to[0,1]$ satisfying $(2\pi)^{-1}\int_0^\infty q_\xi\,d\xi=x$.
Then
\[
J_\infty(x)
=\inf_{q\in\mathcal Q_x}\frac1{2\pi}\int_0^\infty
\operatorname{KL}\bigl(\mathrm{Ber}(q_\xi)\Vert\mathrm{Ber}(e^{-2\xi})\bigr)\,d\xi.
\]
For $x>0$, the unique minimizer up to null sets is
\[
q_\xi^{\theta_x}
=\frac{e^{\theta_x}e^{-2\xi}}{1+(e^{\theta_x}-1)e^{-2\xi}},
\qquad \Lambda_\infty'(\theta_x)=x.
\]
For $x=0$, the unique minimizer is $q=0$ almost everywhere, and $J_\infty(0)=\pi/24$.
\end{proposition}

\begin{proof}
Put $p_\xi=e^{-2\xi}$, $M_\xi(\theta)=1+(e^\theta-1)p_\xi$, and $q_\xi^\theta=e^\theta p_\xi/M_\xi(\theta)$.
The Bernoulli identity
\[
\operatorname{KL}(\mathrm{Ber}(q_\xi)\Vert\mathrm{Ber}(p_\xi))
=\operatorname{KL}(\mathrm{Ber}(q_\xi)\Vert\mathrm{Ber}(q_\xi^\theta))
+\theta q_\xi-\log M_\xi(\theta)
\]
gives a cost at least $\theta x-\Lambda_\infty(\theta)$ for $q\in\mathcal Q_x$.
Taking the supremum proves the lower bound.
For $x>0$, $q^{\theta_x}\in\mathcal Q_x$ since $(2\pi)^{-1}\int_0^\infty q_\xi^{\theta_x}\,d\xi=\Lambda_\infty'(\theta_x)=x$.
It has finite entropy and attains equality, uniquely by strictness of relative entropy.
At $x=0$, $q=0$ almost everywhere, with cost $-(2\pi)^{-1}\int_0^\infty\log(1-e^{-2\xi})\,d\xi=\pi/24$.
\end{proof}

The proof of Theorem~\ref{thm:precise-halfline} uses the exponential-tilting and local-limit argument of the bounded-height analysis in \cite{Zhou2026}, with the strong half-line expansion supplied by Theorem~\ref{thm:strong-halfline}.

\begin{proof}[Proof of Theorem~\ref{thm:precise-halfline}]
Let $K_L(z)=\log\Ee e^{zN_\infty(L)}$, with its normalized holomorphic branch, and set $d\Pp_{\theta,L}=e^{\theta N_\infty(L)-K_L(\theta)}d\Pp$ for real $\theta$.
At the saddle $\theta=\theta_{L,k}$, the exact change of measure gives
\[
\Pp\{N_\infty(L)=k\}
=e^{K_L(\theta)-\theta k}\Pp_{\theta,L}\{N_\infty(L)=k\}.
\]
To estimate the tilted probability, fix $\Theta\Subset\mathbb R$.
By Lemma~\ref{lem:fredholm-rectangles-halfline}, writing $(\lambda_{j,L})_j$ for the eigenvalues of the restricted determinantal operator, the tilted count is a sum of independent Bernoulli variables with parameters $p_{j,L}(\theta)=e^\theta\lambda_{j,L}/[1+(e^\theta-1)\lambda_{j,L}]$ and variance $v_L(\theta)=K_L''(\theta)$.
These parameters are the eigenvalues of the positive trace-class contraction $e^\theta K_L^\infty(\Id+(e^\theta-1)K_L^\infty)^{-1}$.
The strong expansion and Cauchy's estimates give, uniformly for $\theta\in\Theta$,
\[
m_L(\theta):=K_L'(\theta)=L\Lambda_\infty'(\theta)+O_\Theta(\log L),
\qquad
v_L(\theta)=L\Lambda_\infty''(\theta)+O_\Theta(\log L).
\]
In particular, $v_L(\theta)\ge c_\Theta L$.
After omitting any deterministic Bernoulli summands, the local central limit theorem \cite[Theorem~1]{ForresterLebowitz2014} applies to every sequence $\theta_L\in\Theta$.
Its conclusion is therefore uniform on $\Theta$ and gives
$\sup_n\Pp_{\theta,L}\{N_\infty(L)=n\}=O_\Theta(L^{-1/2})$.
For fixed $M\ge0$, if $|n-L\Lambda_\infty'(\theta)|\le M$, then
$(n-m_L(\theta))/\sqrt{v_L(\theta)}=O_{\Theta,M}(\log L/\sqrt L)\to0$.
The same theorem consequently gives
\begin{equation}\label{eq:tilted-local-halfline}
\Pp_{\theta,L}\{N_\infty(L)=n\}
=\frac{1+o_{\Theta,M}(1)}{\sqrt{2\pi L\Lambda_\infty''(\theta)}},
\qquad |n-L\Lambda_\infty'(\theta)|\le M.
\end{equation}

For $\theta=\theta_{L,k}$, so that $k=L\Lambda_\infty'(\theta)$, the strong expansion and Legendre duality give
\[
K_L(\theta)-\theta k
=-LJ_\infty(k/L)+\frac{\theta^2}{2\pi^2}\log L
+A_\infty(\theta)+o(1).
\]
The change-of-measure identity and \eqref{eq:tilted-local-halfline} prove \eqref{eq:precise-point-halfline}, uniformly for $k/L\in K$.
For $k/L\in K_+$, the same change of measure gives
\[
\frac{\Pp\{N_\infty(L)\ge k\}}{\Pp\{N_\infty(L)=k\}}
=\sum_{m\ge0}e^{-m\theta_{L,k}}
\frac{\Pp_{\theta_{L,k},L}\{N_\infty(L)=k+m\}}
{\Pp_{\theta_{L,k},L}\{N_\infty(L)=k\}}.
\]
For each fixed $m$, the fraction tends uniformly to one by \eqref{eq:tilted-local-halfline}.
The uniform point-mass bound and the local lower bound control these fractions independently of $m$.
Since $\theta_{L,k}$ is uniformly bounded away from zero on $K_+$, geometric domination proves \eqref{eq:precise-upper-halfline}.
On $K_-$, the analogous sum has weights $e^{m\theta_{L,k}}$ and is truncated at $m=k$.
The negative saddles are uniformly bounded away from zero, and $k\to\infty$ uniformly, so the same argument proves \eqref{eq:precise-lower-halfline}.
\end{proof}

\section{Moving upper cutoffs and boundary crossover}\label{sec:cutoffs}

We prove Theorem~\ref{thm:intro-moving-cutoff}, which establishes strong asymptotics for the logarithmic moment generating function under moving upper cutoffs $b=b_L\to\infty$.
The critical regime $b_L\asymp L$ is treated by a relative Fredholm determinant, while $1\ll b_L\ll L$ requires a separate boundary analysis.
Recall that, for $b>1$, the limiting pressure of $N_{[1,b]}(L)$ is
\[
\Lambda_{[1,b]}(\theta)
=\frac1{2\pi}\int_0^\infty
\Logp\!\left(1+(e^\theta-1)(e^{-2\xi}-e^{-2b\xi})\right)\dd\xi,
\qquad \theta\in S_\pi.
\]
The logarithm is normalized at $\theta=0$.
On every compact subset of $S_\pi$, the integrand has a common integrable majorant in $\xi$, so $\Lambda_{[1,b]}$ is holomorphic on $S_\pi$.
Write
\[
\mathcal F_{b,L}(\theta)
:=\log\Ee e^{\theta N_{[1,b]}(L)}-L\Lambda_{[1,b]}(\theta).
\]
The logarithm of $\Ee e^{\theta N_{[1,b]}(L)}$ is taken on its holomorphic branch on $S_\pi$ normalized at zero.
It exists because the eigenvalues of the Fourier--Laplace compression lie in $[0,1]$, so its Fredholm determinant has no zeros on $S_\pi$.
Thus $\mathcal F_{b,L}$ is holomorphic on $S_\pi$.

We recall the large-height expansion of \cite{AiFangHou2026}, written in terms of $\Lambda_\infty$, and include its short proof for convenience.

\begin{lemma}\label{lem:pressure-height}
For every $K\Subset S_\pi$, as $b\to\infty$,
\begin{equation}\label{eq:pressure-height}
\Lambda_{[1,b]}(\theta)
=\Lambda_\infty(\theta)+b^{-1}\Lambda_\infty(-\theta)
+O_K(b^{-2}),
\end{equation}
uniformly for $\theta\in K$.
\end{lemma}

\begin{proof}
Put $q=e^\theta-1$ and change variables $t=b\xi$ in $\Lambda_{[1,b]}(\theta)-\Lambda_\infty(\theta)$.
With $z=e^{-2t/b}$ and $w=e^{-2t}$, the resulting difference of logarithms is $\Logp(1+q(z-w))-\Logp(1+qz)$, with branches normalized at $\theta=0$.
Since $0\le w\le z\le1$ and $1+qr$ is uniformly bounded away from zero for $\theta\in K$ and $r\in[0,1]$, the derivative of this difference with respect to $z$ is $q^2w/[(1+q(z-w))(1+qz)]=O_K(w)$.
Thus replacing $z$ by $1$ changes the difference of pressures by at most
\[
\frac{C_K}{b}\int_0^\infty e^{-2t}(1-e^{-2t/b})\dd t
=O_K(b^{-2}).
\]
At $z=1$, the normalized logarithms give
$\Logp(1+q(1-w))-\Logp(1+q)
=\Logp(1+(e^{-\theta}-1)w)$.
The resulting integral equals $b^{-1}\Lambda_\infty(-\theta)$, proving \eqref{eq:pressure-height}.
\end{proof}

\subsection{Critical-ratio crossover}

For $h>0$, let $\mathsf P_h$ be the integral operator on $L^2(0,\infty)$ with kernel
\begin{equation}\label{MC:eq:Ph-kernel}
\mathsf P_h(u,v)
=\frac{\sqrt{uv}}{\pi(u+v)}e^{-h(u+v)}X_1(u-v),
\qquad
X_1(t)=\int_0^1e^{\ii tx}\dd x.
\end{equation}
Lemma~\ref{MC:lem:compression} shows that $\mathsf P_h$ is a positive trace-class contraction and converges strongly, as $h\downarrow0$, to a bounded positive contraction $\mathsf P_0$.
For $c>0$ and $\theta\in S_\pi$, put $q_\theta=e^\theta-1$ and define
\begin{equation}\label{MC:eq:D-definition}
D(\theta,c)
:=\det\!\left[
\Id-q_\theta(\Id+q_\theta\mathsf P_0)^{-1}\mathsf P_c
\right].
\end{equation}
The determinant is well defined since $\mathsf P_c$ is trace class and $\Id+q_\theta\mathsf P_0$ is invertible on $S_\pi$.

\begin{theorem}[Critical-ratio crossover]\label{MC:thm:critical}
For every $c>0$, $D(\cdot,c)$ is holomorphic and nonvanishing on $S_\pi$, with $D(0,c)=1$.
Let $\log D(\cdot,c)$ denote its holomorphic logarithm normalized to zero at the origin, and define
\begin{equation}\label{MC:eq:Across}
A_{\mathrm{cross}}(\theta,c)
:=A_\infty(\theta)+\log D(\theta,c)
-\frac1c\Lambda_\infty(-\theta).
\end{equation}
For every $K\Subset S_\pi$ and $0<c_-<c_+<\infty$, as $L\to\infty$,
\begin{equation}\label{MC:eq:critical-uniform}
\sup_{\substack{\theta\in K\\c\in[c_-,c_+]}}
\left|
\log\Ee e^{\theta N_{[1,cL]}(L)}
-L\Lambda_{[1,cL]}(\theta)
-\frac{\theta^2}{2\pi^2}\log L
-A_{\mathrm{cross}}(\theta,c)
\right|\longrightarrow0.
\end{equation}
In particular, if $b_L/L\to c\in(0,\infty)$, the same expansion holds with $\Lambda_{[1,b_L]}$ in the leading term and $A_{\mathrm{cross}}(\theta,c)$ in the order-one term, locally uniformly on $S_\pi$.
Moreover, for every $K\Subset S_\pi$,
\begin{equation}\label{MC:eq:large-c-bound}
\sup_{\theta\in K}
|A_{\mathrm{cross}}(\theta,c)-A_\infty(\theta)|
=O_K(c^{-1}),\qquad c\to\infty.
\end{equation}
If $b_L/L\to\infty$, then for every $K\Subset S_\pi$,
\begin{equation}\label{MC:eq:supercritical-match}
\sup_{\theta\in K}
\left|
\log\Ee e^{\theta N_{[1,b_L]}(L)}
-\log\Ee e^{\theta N_\infty(L)}
\right|
=O_K(L/b_L)=o(1).
\end{equation}
\end{theorem}

Before proving Theorem~\ref{MC:thm:critical}, we establish the compression and relative Fredholm determinant identities, together with the required continuity and convergence estimates.

\begin{lemma}\label{MC:lem:compression}
For $h>0$, $\mathsf P_h$ is a positive trace-class contraction
with $\Tr\mathsf P_h=(4\pi h)^{-1}$.
As $h\downarrow0$, $\mathsf P_h$ converges strongly to a positive
contraction $\mathsf P_0$, and
$0\le\mathsf P_c\le\mathsf P_h\le\mathsf P_0$ whenever $0<h<c$.
For $L>0$ and $b>1$, put $d=L^{-1}$.
Then, for $\theta\in S_\pi$,
\begin{equation}\label{MC:eq:scaling-dets}
\Ee e^{\theta N_\infty(L)}
=\det(\Id+q_\theta\mathsf P_d),
\qquad
\Ee e^{\theta N_{[1,b]}(L)}
=\det(\Id+q_\theta(\mathsf P_d-\mathsf P_{b/L})).
\end{equation}
Moreover, $c\mapsto\mathsf P_c$ is continuous in trace norm
on $(0,\infty)$, with
\begin{equation}\label{MC:eq:trace-cont-c}
\|\mathsf P_c-\mathsf P_{c'}\|_1
=\frac1{4\pi}\left|\frac1c-\frac1{c'}\right|,
\qquad c,c'>0.
\end{equation}
\end{lemma}

\begin{proof}
For compactly supported $f\in L^2(0,\infty)$, define
\[
(\mathsf U_-f)(x+\ii y)
=\frac1{\sqrt\pi}\int_0^\infty
\sqrt u\,e^{-\ii ux-uy}f(u)\dd u.
\]
Plancherel in $x$ and integration in $y$ give
$\|\mathsf U_-f\|_{L^2(\Hh)}=\|f\|_2$,
so $\mathsf U_-$ extends to an isometry.
Set $\Omega_h=[0,1]\times[h,\infty)$ for $h\ge0$.
Direct integration gives
$\mathsf P_h=\mathsf U_-^*\mathbf1_{\Omega_h}\mathsf U_-$
for $h>0$; use the same formula to define $\mathsf P_0$.
This representation proves positivity and contractivity.
The inclusions $\Omega_c\subset\Omega_h\subset\Omega_0$ for $0<h<c$
give the operator inequalities, while strong convergence follows
from $\mathbf1_{\Omega_h}\to\mathbf1_{\Omega_0}$ strongly.
For $h>0$, $\mathbf1_{\Omega_h}\mathsf U_-$ is Hilbert--Schmidt, with
\[
\Tr\mathsf P_h
=\|\mathbf1_{\Omega_h}\mathsf U_-\|_2^2
=\int_h^\infty\frac{\dd y}{4\pi y^2}
=\frac1{4\pi h}.
\]
The order relation and trace formula imply
\eqref{MC:eq:trace-cont-c}, hence trace-norm continuity.

For $a>0$, the Fourier--Laplace compression associated with
heights $[a,\infty)$ has kernel
\[
T_{a,L}^{\infty}(\xi,\eta)
=\frac{\sqrt{\xi\eta}\,e^{-a(\xi+\eta)}}{\pi(\xi+\eta)}
\int_0^L e^{\ii(\xi-\eta)x}\dd x.
\]
The unitary dilation $(U_Lf)(u)=L^{-1/2}f(u/L)$ gives
$U_LT_{a,L}^{\infty}U_L^{-1}=\mathsf P_{a/L}$.
Applying the Fredholm determinant formula to $T_{1,L}^{\infty}$
and $T_{1,L}^{\infty}-T_{b,L}^{\infty}$ proves
\eqref{MC:eq:scaling-dets}.
\end{proof}

\begin{lemma}\label{MC:lem:relative}
For $c>0$ and $0<d=L^{-1}<c$, we have, for $\theta\in S_\pi$,
\begin{equation}\label{MC:eq:finite-ratio}
\frac{\Ee e^{\theta N_{[1,cL]}(L)}}
{\Ee e^{\theta N_\infty(L)}}
=\det\!\left[
\Id-q_\theta(\Id+q_\theta\mathsf P_d)^{-1}\mathsf P_c
\right].
\end{equation}
As $d\downarrow0$, the right-hand side converges to
$D(\theta,c)\ne0$, locally uniformly on $S_\pi\times(0,\infty)$.
The holomorphic logarithms normalized to vanish at $\theta=0$
converge to $\log D(\theta,c)$ with the same uniformity.
\end{lemma}

\begin{proof}
Fix $K\Subset S_\pi$ and $0<c_-<c_+<\infty$.
Put $R_h(\theta)=(\Id+q_\theta\mathsf P_h)^{-1}$ for $h\ge0$,
and write $D_d(\theta,c)$ for the determinant in
\eqref{MC:eq:finite-ratio}.
Spectral calculus gives $\|R_h(\theta)\|\le C_K$ for
$h\ge0$ and $\theta\in K$.
For $0\le h<c$, the factorization
\[
\Id-q_\theta R_h(\theta)\mathsf P_c
=R_h(\theta)\bigl(\Id+q_\theta(\mathsf P_h-\mathsf P_c)\bigr)
\]
has invertible factors, since $\mathsf P_h$ and
$\mathsf P_h-\mathsf P_c$ are positive contractions.
For $h=d$, multiplicativity of Fredholm determinants and
\eqref{MC:eq:scaling-dets} give \eqref{MC:eq:finite-ratio}.
The cases $h=d$ and $h=0$ also prove that $D_d(\cdot,c)$ and
$D(\cdot,c)$ are nonvanishing.
Both are holomorphic on $S_\pi$ and equal one at zero,
so their normalized holomorphic logarithms exist.

Since $\mathsf P_d\to\mathsf P_0$ strongly, the resolvent identity
and the uniform inverse bounds give strong convergence
$R_d(\theta)\to R_0(\theta)$, uniformly for $\theta\in K$
on each fixed vector.
By \eqref{MC:eq:trace-cont-c},
$\{\mathsf P_c:c\in[c_-,c_+]\}$ is compact in trace norm.
Finite-rank approximation therefore gives
\[
\sup_{\substack{\theta\in K\\c\in[c_-,c_+]}}
\|(R_d(\theta)-R_0(\theta))\mathsf P_c\|_1
\longrightarrow0.
\]
Trace-norm continuity of Fredholm determinants yields
$D_d\to D$ locally uniformly.
On a compact neighborhood of the segments joining zero to $K$,
continuity and nonvanishing of $D$ imply that $D_d/D\to1$
uniformly for $c\in[c_-,c_+]$.
Its principal logarithm consequently tends uniformly to zero
and, by normalization at zero, equals the difference of the
normalized holomorphic logarithms.
This proves the final assertion.
\end{proof}

\begin{proof}[Proof of Theorem~\ref{MC:thm:critical}]
Fix $K\Subset S_\pi$ and $0<c_-<c_+<\infty$.
The holomorphy and nonvanishing of $D(\cdot,c)$, with $D(0,c)=1$, follow from Lemma~\ref{MC:lem:relative}.
The determinant identity \eqref{MC:eq:finite-ratio}, convergence of the normalized logarithms, and the strong half-line expansion give
\[
\log\Ee e^{\theta N_{[1,cL]}(L)}
=L\Lambda_\infty(\theta)
+\frac{\theta^2}{2\pi^2}\log L
+A_\infty(\theta)+\log D(\theta,c)+o(1)
\]
uniformly for $\theta\in K$ and $c\in[c_-,c_+]$.
On the same parameter sets, Lemma~\ref{lem:pressure-height} gives
\[
L\Lambda_{[1,cL]}(\theta)
=L\Lambda_\infty(\theta)
+c^{-1}\Lambda_\infty(-\theta)
+O_K((c^2L)^{-1}).
\]
Combining these expansions proves \eqref{MC:eq:critical-uniform}.
If $c_L=b_L/L\to c_0\in(0,\infty)$, uniformity of the expansion and continuity of $A_{\mathrm{cross}}(\theta,c)$ in $c$ give the asserted asymptotic with the actual pressure $\Lambda_{[1,b_L]}$.

Applying \eqref{eq:height-truncation-log} with $b=cL$ bounds the difference of the logarithmic moment generating functions by $C_K/c$, with $C_K$ independent of $c$ and $L$ whenever $cL>1$.
For each fixed $c>0$, let $L\to\infty$ and use Lemma~\ref{MC:lem:relative} to obtain $\sup_{\theta\in K}|\log D(\theta,c)|\le C_K/c$.
Equation~\eqref{MC:eq:Across} then gives \eqref{MC:eq:large-c-bound}.
Taking $b=b_L$ directly in \eqref{eq:height-truncation-log} proves \eqref{MC:eq:supercritical-match}.
\end{proof}

We next identify the mean and variance corrections encoded by the crossover determinant.

\begin{proposition}\label{MC:prop:checks}
For every $c>0$, $\partial_\theta\log D(0,c)=-(4\pi c)^{-1}$ and $A_{\mathrm{cross}}(0,c)=\partial_\theta A_{\mathrm{cross}}(0,c)=0$.
Moreover,
\[
\lim_{L\to\infty}
\bigl(\operatorname{Var}N_{[1,cL]}(L)
-\operatorname{Var}N_\infty(L)\bigr)
=\partial_\theta^2\log D(0,c)=:\Delta(c),
\]
where
\[
\Delta(c)
=-\frac1{4\pi c}
+\frac1{\pi^2}\int_c^{2c}
\frac{2(t-c)\arctan(1/t)}{t^3}\dd t
+\frac1{\pi^2}\int_{2c}^\infty
\frac{\arctan(1/t)}{t^2}\dd t.
\]
In particular,
\begin{align}
\Delta(c)
&=\frac1{8\pi c}+\frac{\log(c/2)}{\pi^2}
-\frac{c^2}{3\pi^2}+O(c^4),
&&c\downarrow0,
\label{MC:eq:Delta-small}\\
\Delta(c)
&=-\frac1{4\pi c}+\frac7{24\pi^2c^2}+O(c^{-4}),
&&c\to\infty.
\label{MC:eq:Delta-large}
\end{align}
Consequently,
\begin{equation}\label{MC:eq:Across-curvature}
A_{\mathrm{cross}}''(0,c)
=-\frac{\log2}{\pi^2}+\Delta(c)-\frac1{8\pi c},
\end{equation}
and $A_{\mathrm{cross}}''(0,c)=\pi^{-2}\log(c/4)+O(c^2)$ as $c\downarrow0$.
\end{proposition}

\begin{proof}
Differentiating \eqref{MC:eq:D-definition} at $\theta=0$ gives $\partial_\theta\log D(0,c)=-\Tr\mathsf P_c=-(4\pi c)^{-1}$ and
\[
\partial_\theta^2\log D(0,c)
=-\Tr\mathsf P_c
+2\Tr(\mathsf P_0\mathsf P_c)
-\Tr(\mathsf P_c^2).
\]
By Lemma~\ref{MC:lem:relative} and Cauchy's formula, the variance difference converges to $\partial_\theta^2\log D(0,c)$.
The strong half-line expansion, together with $\Ee N_\infty(L)=L/(4\pi)$ and $\Lambda_\infty'(0)=1/(4\pi)$, gives $A_\infty(0)=A_\infty'(0)=0$.
Equation~\eqref{MC:eq:Across} then yields the asserted normalization of $A_{\mathrm{cross}}$.

By trace cyclicity and the compression representation in Lemma~\ref{MC:lem:compression}, the quantity $2\Tr(\mathsf P_0\mathsf P_c)-\Tr(\mathsf P_c^2)$ equals the integral of $|\mathcal K_{\Hh}(z,w)|^2$ over $(\Omega_0\times\Omega_c)\cup(\Omega_c\times\Omega_0)$, where $\Omega_h=[0,1]\times[h,\infty)$.
For fixed $t=y+y'$, the admissible values of $y$ have total length $2(t-c)$ for $c<t<2c$, $t$ for $t\ge2c$, and zero for $t\le c$.
Together with \eqref{eq:HS-kernel}, the identity
\[
2\int_0^1\frac{1-r}{(r^2+t^2)^2}\dd r
=\frac{\arctan(1/t)}{t^3}
\]
gives the stated integral representation of $\Delta(c)$.

Integration by parts and expansion as $c\downarrow0$ give \eqref{MC:eq:Delta-small}.
Substituting $\arctan(1/t)=t^{-1}+O(t^{-3})$ into the integral representation gives \eqref{MC:eq:Delta-large}.
Finally, the half-line variance asymptotic gives $A_\infty''(0)=-\pi^{-2}\log2$, while $\Lambda_\infty''(0)=1/(8\pi)$.
Differentiating \eqref{MC:eq:Across} twice proves \eqref{MC:eq:Across-curvature}.
Substitution of \eqref{MC:eq:Delta-small} into \eqref{MC:eq:Across-curvature} gives the final small-$c$ expansion.
\end{proof}

The integral formula in Proposition~\ref{MC:prop:checks}
also determines the global behavior of the variance correction.

\begin{corollary}\label{MC:cor:variance-threshold}
There exists a unique $c_*\in(0,\infty)$ such that
$\Delta(c_*)=0$, and
\[
\Delta(c)>0\quad(0<c<c_*),\qquad
\Delta(c)<0\quad(c>c_*).
\]
Moreover, $c\mapsto A_{\mathrm{cross}}''(0,c)$
is strictly increasing on $(0,\infty)$.
\end{corollary}

\begin{proof}
Let $w(u)=2(u-1)u^{-3}$ for $1\le u\le2$
and $w(u)=u^{-2}$ for $u>2$.
Then $w\ge0$ and $\int_1^\infty w(u)\,\dd u=3/4$.
Changing variables $t=cu$ in the integral formula of
Proposition~\ref{MC:prop:checks} gives
\[
c\Delta(c)=-\frac1{4\pi}
+\frac1{\pi^2}\int_1^\infty
w(u)\arctan\frac1{cu}\,\dd u.
\]
The right-hand side is continuous and strictly decreasing
in $c>0$, with limits $1/(8\pi)$ and $-1/(4\pi)$
as $c\downarrow0$ and $c\to\infty$, respectively,
by dominated convergence.
This proves the first assertion.

Combining this representation with
\eqref{MC:eq:Across-curvature} and
$\arctan(1/t)=\pi/2-\arctan t$ for $t>0$ gives
\[
\pi^2 A_{\mathrm{cross}}''(0,c)
=-\log2-\frac1c\int_1^\infty
w(u)\arctan(cu)\,\dd u.
\]
For every $u>0$, the function
$c\mapsto c^{-1}\arctan(cu)$ is strictly decreasing,
since $\arctan t>t/(1+t^2)$ for $t>0$.
As $w>0$ almost everywhere on $(1,\infty)$,
the second assertion follows.
\end{proof}

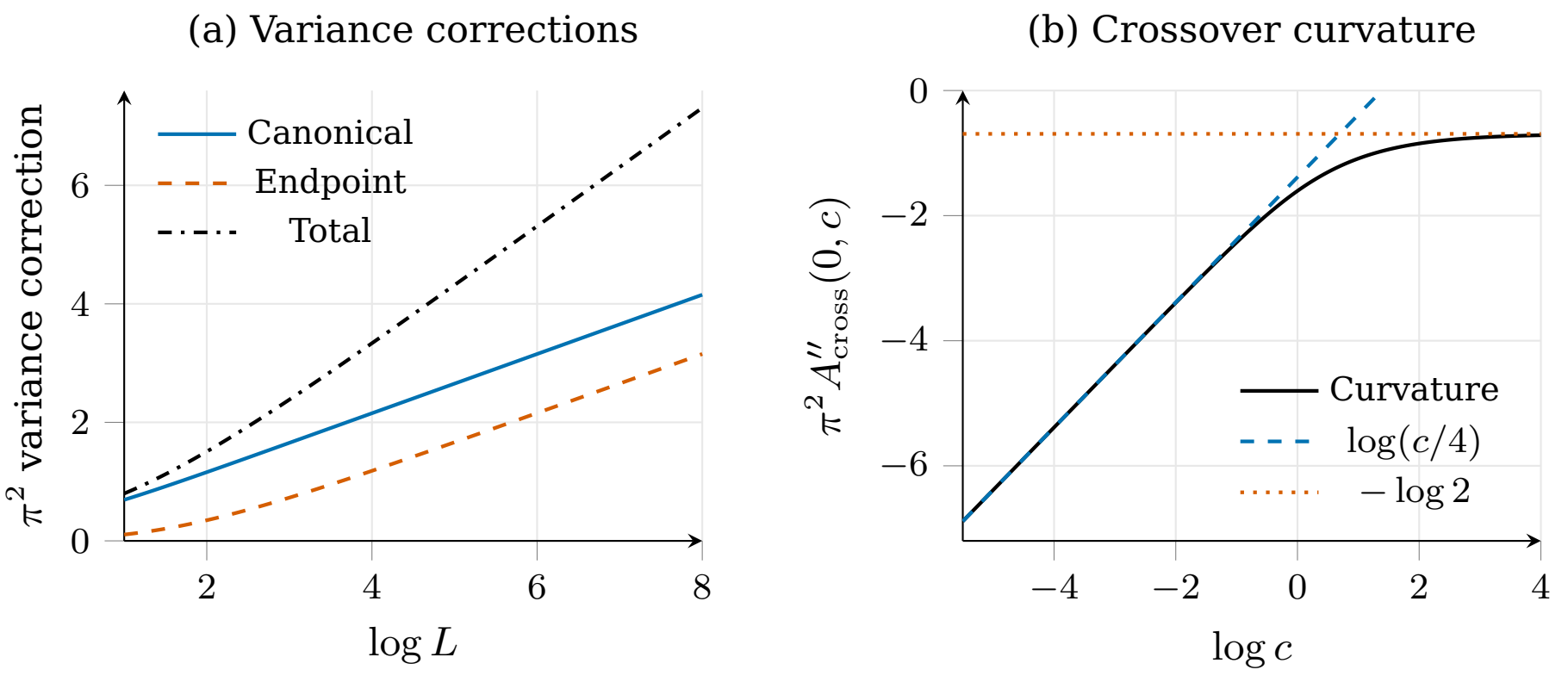
\begin{figure}[!htbp]
\centering
\begin{tikzpicture}
\pgfmxfpdeclarefunction{hscan}{1}{
  exp(#1)*atan(2*exp(-#1))/4+ln(1+exp(2*#1)/4)/4
}
\pgfmxfpdeclarefunction{hsend}{1}{
  exp(-#1)*atan(exp(#1)/2)+ln(1+exp(2*#1)/4)/4-1/2
}
\pgfmxfpdeclarefunction{hscurv}{1}{
  #1-ln(4)+ln(1+4*exp(2*#1))/2-ln(1+exp(2*#1))
  +1/2-exp(-#1)*atan(exp(#1))
  +exp(-#1)*atan(2*exp(#1))/4
  -exp(#1)*atan(1/(exp(-#1)+2*exp(#1)))
}
\begin{groupplot}[
 group style={group size=2 by 1,horizontal sep=2.0cm},
 width=0.43\textwidth,height=0.36\textwidth,
 axis lines=left,tick align=outside,
 tick label style={font=\footnotesize},label style={font=\small},
 title style={font=\small},
 legend style={font=\footnotesize,draw=none,fill=none,
   at={(0.03,0.98)},anchor=north west},
 grid=major,grid style={black!9},
 samples=201,enlargelimits=false]
\nextgroupplot[
 title={(a) Variance corrections},
 xlabel={$\log L$},
 ylabel={$\pi^2$ variance correction},
 domain=1:8,xmin=1,xmax=8,ymin=0,ymax=7.6,
 xtick={2,4,6,8},ytick={0,2,4,6}]
\addplot[HSblue,thick] {hscan(x)};
\addlegendentry{Canonical}
\addplot[HSorange,thick,dashed] {hsend(x)};
\addlegendentry{Endpoint}
\addplot[black,thick,dashdotted] {hscan(x)+hsend(x)};
\addlegendentry{Total}
\nextgroupplot[
 title={(b) Crossover curvature},
 xlabel={$\log c$},
 ylabel={$\pi^2 A_{\mathrm{cross}}''(0,c)$},
 domain=-5.5:4,xmin=-5.5,xmax=4,ymin=-7.2,ymax=0,
 xtick={-4,-2,0,2,4},ytick={-6,-4,-2,0},
 legend style={at={(0.97,0.03)},anchor=south east}]
\addplot[black,thick] {hscurv(x)};
\addlegendentry{Curvature}
\addplot[HSblue,dashed,thick,samples=2] {x-ln(4)};
\addlegendentry{$\log(c/4)$}
\addplot[HSorange,dotted,thick,samples=2] {-ln(2)};
\addlegendentry{$-\log 2$}
\end{groupplot}
\end{tikzpicture}
\caption{Variance corrections (left) and the crossover variance constant (right).}
\label{fig:variance-crossover}
\end{figure}

Numerically solving $\Delta(c)=0$ gives $c_*\approx0.15190$.
Figure~\ref{fig:variance-crossover} (left) shows the canonical and endpoint contributions to
$\pi^2[\operatorname{Var}N_\infty(L)-L/(8\pi)]$, together with their sum.
The right panel plots $\pi^2 A_{\mathrm{cross}}''(0,c)$,
with its small-$c$ asymptote $\log(c/4)$ and large-$c$ limit $-\log2$.

\FloatBarrier

The large-ratio estimate in Theorem~\ref{MC:thm:critical} admits the following first-order refinement.

\begin{proposition}\label{BC:prop:large}
For every $K\Subset S_\pi$, uniformly for $\theta\in K$,
\begin{equation}\label{BC:eq:large}
A_{\mathrm{cross}}(\theta,c)=A_\infty(\theta)
-\frac1c\left(\frac{e^\theta-1}{4\pi}
+\Lambda_\infty(-\theta)\right)+O_K(c^{-2}),
\qquad c\to\infty.
\end{equation}
\end{proposition}

\begin{proof}
Fix $K\Subset S_\pi$ and put $R_0(\theta)=(\Id+q_\theta\mathsf P_0)^{-1}$.
Since $\mathsf P_0$ is a positive contraction, spectral calculus gives $\sup_{\theta\in K}\|R_0(\theta)\|<\infty$.
By Lemma~\ref{MC:lem:compression}, $\|\mathsf P_c\|_1=(4\pi c)^{-1}$, so the normalized Fredholm logarithm admits the expansion
\[
\log D(\theta,c)
=-q_\theta\Tr(R_0(\theta)\mathsf P_c)+O_K(c^{-2})
\]
uniformly for $\theta\in K$ as $c\to\infty$.

Since $R_0(\theta)-\Id=-q_\theta R_0(\theta)\mathsf P_0$, spectral calculus and positivity of $\mathsf P_c$ give $|\Tr[(R_0(\theta)-\Id)\mathsf P_c]|\le C_K\Tr(\mathsf P_0\mathsf P_c)$.
By \eqref{MC:eq:Ph-kernel}, strong convergence $\mathsf P_h\to\mathsf P_0$, and $|X_1|\le1$,
\[
0\le\Tr(\mathsf P_0\mathsf P_c)
\le\frac1{\pi^2}\int_0^\infty\int_0^\infty
\frac{uv\,e^{-c(u+v)}}{(u+v)^2}\dd u\dd v
=\frac1{6\pi^2c^2}.
\]
Consequently, $\log D(\theta,c)=-q_\theta/(4\pi c)+O_K(c^{-2})$, and substitution into \eqref{MC:eq:Across} proves \eqref{BC:eq:large}.
\end{proof}

\subsection{Uniform boundary estimates}

We compare $\mathcal F_{b,L}$ with the fixed-height constant $\mathcal C_b=\lim_{L\to\infty}\mathcal F_{b,L}$ from \eqref{eq:intro-recalled-interval-strong}.
Keeping the height variable in the Fourier multiplier allows us to separate the two spatial boundaries uniformly in $b$.
For a separable Hilbert space $H$, put $(\mathcal F_Hu)(\xi)=(2\pi)^{-1/2}\int_{\Rr}e^{-\ii x\xi}u(x)\dd x$ and use the multiplier $g(-D)=\mathcal F_H^{-1}M_g\mathcal F_H$ on $L^2(\Rr;H)$, consistently with Section~\ref{sec:WH-halfline}.
All trace-ideal norms include $H$.

\begin{lemma}\label{UC:lem:gap}
Let $g\in W^{2,1}((0,\infty);\mathfrak S_1(H))$ with $g(0+)=0$, and extend $g$ by zero to negative frequencies.
For $h>0$, put $Q_-=\mathbf1_{(-\infty,0)}\otimes\Id_H$ and $Q_+^h=\mathbf1_{(h,\infty)}\otimes\Id_H$.
Then
\begin{equation}\label{UC:eq:gap-bound}
 \|Q_-g(-D)Q_+^h\|_1+\|Q_+^h g(-D)Q_-\|_1
 \le\frac{\|g'(0+)\|_1+\int_0^\infty\|g''(\xi)\|_1\dd\xi}{2\pi h}.
\end{equation}
Adding a constant scalar multiple of $\Id_H$ to the symbol leaves these operators unchanged.
\end{lemma}

\begin{proof}
The absolutely continuous representatives of $g,g'$ vanish at infinity, so two Bochner integrations by parts give
\[
 k(t):=\frac1{2\pi}\int_0^\infty e^{\ii t\xi}g(\xi)\dd\xi
 =\frac{g'(0+)+\int_0^\infty e^{\ii t\xi}g''(\xi)\dd\xi}
 {2\pi(\ii t)^2},\qquad t\ne0.
\]
Reflection and translation give off-diagonal kernels $k(\pm(u+v+h))$, $u,v>0$.
The scalar Hankel operator $H_h$ with kernel $(u+v+h)^{-2}$ is positive and trace class, since
\[
 (u+v+h)^{-2}=\int_0^\infty t e^{-ht}e^{-tu}e^{-tv}\dd t,
 \qquad \|H_h\|_1=\frac1{2h}.
\]
For $C\in\mathfrak S_1(H)$, $\|H_h\otimes C\|_1=\|C\|_1/(2h)$.
Multiplication of its kernel by $e^{\pm\ii\xi(u+v+h)}$ preserves this norm by left and right unitary multiplication.
The formula for $k$ is thus an absolutely convergent integral in $\mathfrak S_1$ and bounds each off-diagonal operator by the numerator in \eqref{UC:eq:gap-bound} divided by $4\pi h$.
Summation proves the estimate, and the constant symbol has zero off-diagonal compressions.
\end{proof}

\begin{lemma}\label{SB:lem:boundary-tail}
For every $K\Subset S_\pi$, there exists $C_K<\infty$ such that
\begin{equation}\label{SB:eq:uniform-constant-approach}
\sup_{\theta\in K}
|\mathcal F_{b,L}(\theta)-\mathcal C_b(\theta)|
\le C_K\frac bL,
\qquad L\ge b\ge2.
\end{equation}
\end{lemma}

\begin{proof}
Fix $K\Subset S_\pi$ and set $\varepsilon=b^{-1}$, $T=L/b$, $H=L^2(0,1)$, and $\mathcal H=L^2(\Rr;H)$.
For $\xi>0$, put $v_{\varepsilon,\xi}(y)=e^{-\xi y}\mathbf1_{[\varepsilon,1]}(y)$ and $\mathsf a_\varepsilon(\xi)=2\xi\,|v_{\varepsilon,\xi}\rangle\langle v_{\varepsilon,\xi}|$, and extend $\mathsf a_\varepsilon$ by zero to $\xi\le0$.
Its nonzero eigenvalue is $r_\varepsilon(\xi)=e^{-2\varepsilon\xi}-e^{-2\xi}\in[0,1]$, so $A=\mathsf a_\varepsilon(-D)$ is a positive contraction.
Let $P=\mathbf1_{(0,T)}\otimes\Id_H$ and regard $A_T=PAP$ on $P\mathcal H$.
Fourier inversion and planar dilation by $b^{-1}$ identify the nonzero spectrum of $A_T$ with that of the counting operator for $N_{[1,b]}(L)$.
The factorization $A_T=V^*V$, where $Vu=\mathsf a_\varepsilon^{1/2}\mathcal F_H(Pu)$, gives $\Tr A_T=T(\varepsilon^{-1}-1)/(4\pi)$.
Put $\ell_\theta(t)=\Logp(1+(e^\theta-1)t)$.
Since $\ell_\theta(t)/t$ has a removable singularity at zero, both $\ell_\theta(A_T)$ and $P\ell_\theta(A)P$ are trace class, and
\begin{equation}\label{UC:eq:height-trace}
\log\Ee e^{\theta N_{[1,b]}(L)}
=\Tr\ell_\theta(A_T),\qquad
\Tr(P\ell_\theta(A)P)
=\frac{T}{2\pi}\int_0^\infty
\ell_\theta(r_\varepsilon(\xi))\dd\xi
=L\Lambda_{[1,b]}(\theta).
\end{equation}
Choose a positively oriented contour $\Gamma$ enclosing $[0,1]$, whose closed interior lies in a common domain of holomorphy of $\ell_\theta$, $\theta\in K$.
For $\zeta\in\Gamma$, put $d=\operatorname{dist}(\Gamma,[0,1])>0$ and
\[
B=(\zeta\Id-A)^{-1},
\qquad
g_\zeta(\xi)
=(\zeta\Id_H-\mathsf a_\varepsilon(\xi))^{-1}
-\zeta^{-1}\Id_H.
\]
Since $\partial_\xi^j v_{\varepsilon,\xi}(y)=(-y)^jv_{\varepsilon,\xi}(y)$,
integration over $y\in[\varepsilon,1]$ gives, uniformly in $0<\varepsilon\le1/2$,
\[
\|\partial_\xi^j v_{\varepsilon,\xi}\|_{L^2(0,1)}
\le C_j(1+\xi)^{-j-1/2},\qquad j=0,1,2.
\]
Differentiating the rank-one symbol gives
$\|\mathsf a_\varepsilon'(\xi)\|_1\le C(1+\xi)^{-1}$
and $\|\mathsf a_\varepsilon''(\xi)\|_1\le C(1+\xi)^{-2}$.
For fixed $\varepsilon$, exponential decay at infinity gives
$g_\zeta\in W^{2,1}((0,\infty);\mathfrak S_1(H))$.
Writing $R_\zeta=(\zeta\Id_H-\mathsf a_\varepsilon)^{-1}$, we have
$g_\zeta'=R_\zeta\mathsf a_\varepsilon'R_\zeta$ and
$g_\zeta''=R_\zeta\mathsf a_\varepsilon''R_\zeta
+2R_\zeta\mathsf a_\varepsilon'R_\zeta\mathsf a_\varepsilon'R_\zeta$.
The uniform resolvent bound therefore yields
\begin{equation}\label{UC:eq:height-regularity}
g_\zeta(0+)=0,
\qquad
\|g_\zeta'(0+)\|_1
+\int_0^\infty\|g_\zeta''(\xi)\|_1\dd\xi
\le C_\Gamma.
\end{equation}
Let $Q_-=\mathbf1_{(-\infty,0)}\otimes\Id_H$ and $Q_+=\mathbf1_{(T,\infty)}\otimes\Id_H$.
Since $B-\zeta^{-1}\Id=g_\zeta(-D)$, Lemma~\ref{UC:lem:gap} gives
\begin{equation}\label{UC:eq:height-resolvent-gap}
\|Q_-BQ_+\|_1+\|Q_+BQ_-\|_1
\le C_\Gamma/T.
\end{equation}

Put $Q=Q_-+Q_+$, $E_\ell=\Id-Q_-$, and $E_r=\Id-Q_+$.
For an orthogonal projection $E$, write $B_E=EBE|_{E\mathcal H}$ and extend compressed operators and their inverses by zero.
Since every compression of $A$ is a positive contraction, Schur complementation gives $\|B_E^{-1}\|\le|\zeta|+1+d^{-1}$ and
\begin{equation}\label{UC:eq:boundary-defect}
\mathcal D_E
:=(\zeta\Id_{E\mathcal H}-EAE)^{-1}-EBE
=-EBE^\perp B_{E^\perp}^{-1}E^\perp BE.
\end{equation}
For fixed $\varepsilon$, the convolution kernel $k_\zeta$ of $g_\zeta(-D)$ is bounded in $\mathfrak S_2(H)$ near zero by integrability of $g_\zeta$.
Two integrations by parts using \eqref{UC:eq:height-regularity} give $\|k_\zeta(t)\|_2=O_\Gamma(|t|^{-2})$ at infinity.
Thus $\int_{\Rr}|t|\|k_\zeta(t)\|_2^2\dd t<\infty$, so the off-diagonal factors in \eqref{UC:eq:boundary-defect} are Hilbert--Schmidt for $E=P,E_\ell,E_r$.
The corresponding defects are therefore trace class and continuous in trace norm on $\Gamma$.
On $Q\mathcal H=Q_-\mathcal H\oplus Q_+\mathcal H$, \eqref{UC:eq:height-resolvent-gap}, the inverse identity, and the uniform inverse bounds imply
\[
\bigl\|B_Q^{-1}
-(B_{Q_-}^{-1}\oplus B_{Q_+}^{-1})\bigr\|_1
\le C_\Gamma/T.
\]
The Schur-complement identity \eqref{UC:eq:boundary-defect} gives
\[
\mathcal D_P-P(\mathcal D_{E_\ell}+\mathcal D_{E_r})P
=-PBQ\bigl[B_Q^{-1}-(B_{Q_-}^{-1}\oplus B_{Q_+}^{-1})\bigr]QBP.
\]
The preceding inverse bound therefore controls its trace norm by $C_\Gamma/T$.
Expanding $E_\ell=P+Q_+$ and $E_r=P+Q_-$ also gives
$\|\mathcal D_{E_j}-P\mathcal D_{E_j}P\|_1\le C_\Gamma/T$ for $j=\ell,r$, since every term in the expanded differences contains $Q_-BQ_+$ or $Q_+BQ_-$.
Consequently,
\begin{equation}\label{UC:eq:boundary-separation}
\|\mathcal D_P-\mathcal D_{E_\ell}-\mathcal D_{E_r}\|_1
\le C_\Gamma/T,
\end{equation}
uniformly in $0<\varepsilon\le1/2$ and $T\ge1$.

Holomorphic functional calculus and \eqref{UC:eq:height-trace} give
\[
\mathcal F_{b,L}(\theta)
=\frac1{2\pi\ii}\int_\Gamma
\ell_\theta(\zeta)\Tr\mathcal D_P(\zeta)\dd\zeta,
\]
with interchange of trace and integration justified by trace-norm continuity on $\Gamma$.
By \eqref{UC:eq:boundary-separation}, replacing $\mathcal D_P$ by $\mathcal D_{E_\ell}+\mathcal D_{E_r}$ changes the integral by at most $C_K/T$.
The resulting integral is independent of $T$: $E_\ell$ is fixed, and translation invariance of $A$ makes $\Tr\mathcal D_{E_r}$ independent of $T$.
For fixed $b$, letting $L\to\infty$ identifies this integral with $\mathcal C_b$ by \eqref{eq:intro-recalled-interval-strong}.
Since the error is uniform for $b\ge2$ and $T\ge1$, \eqref{SB:eq:uniform-constant-approach} follows.
\end{proof}

\begin{proposition}\label{UC:prop:bound}
For every $K\Subset S_\pi$, there exists $C_K<\infty$ such that
\begin{equation}\label{UC:eq:bound}
 \sup_{L\ge b\ge2}\sup_{\theta\in K}
 \left|\mathcal F_{b,L}(\theta)-\frac{\theta^2}{2\pi^2}\log(b/2)\right|
 \le C_K.
\end{equation}
\end{proposition}

\begin{proof}
Theorem~\ref{MC:thm:critical} at $c=1$ bounds $\mathcal F_{b,b}-\theta^2\log(b/2)/(2\pi^2)$ uniformly on $K$ for large $b$.
For $b$ in a bounded interval, the positive-contraction estimate $|\Tr\ell_\theta(A_T)|\le C_K\Tr A_T$ and the pressure integral give the same bound at $L=b$.
Apply \eqref{SB:eq:uniform-constant-approach} first at $L=b$ to bound $\mathcal C_b-\theta^2\log(b/2)/(2\pi^2)$ and then at general $L\ge b$ to obtain \eqref{UC:eq:bound}.
\end{proof}

The following consequence of the three-circles theorem gives quantitative bounds from the local estimates.

\begin{lemma}\label{UC:lem:three-circles}
Let $\{F_i:i\in\mathcal I\}$ be a locally uniformly bounded family of holomorphic functions on $S_\pi$, and let $\eta_i\in(0,1]$.
Suppose that $\sup_{\theta\in U}|F_i(\theta)|\le C\eta_i^p$ for a fixed neighborhood $U\subset S_\pi$ of zero, constants $C,p>0$, and all $i\in\mathcal I$.
For every $K\Subset S_\pi$, there exist $C_K<\infty$ and $\delta_K>0$ such that
\[
 \sup_{\theta\in K}|F_i(\theta)|\le C_K\eta_i^{\delta_K},
 \qquad i\in\mathcal I.
\]
\end{lemma}

\begin{proof}
Fix $K\Subset S_\pi$.
The map $\psi(\theta)=\tanh(\theta/4)$ takes $S_\pi$ biholomorphically onto the unit disk.
Choose $0<r_0<r_K<r_1<1$ such that $\psi^{-1}(\{|w|\le r_0\})\subset U$ and $\psi(K)\subset\{|w|\le r_K\}$.
Applying Hadamard's three-circles theorem \cite[Chapter~VI, Theorem~3.13]{Conway1978} to $F_i\circ\psi^{-1}$ and using local uniform boundedness, we obtain
\[
 \sup_{\theta\in K}|F_i(\theta)|
 \le C_K\eta_i^{p\alpha},
 \qquad
 \alpha=\frac{\log(r_1/r_K)}{\log(r_1/r_0)}\in(0,1),
\]
with $C_K$ independent of $i$.
Taking $\delta_K=p\alpha$ proves the claim.
\end{proof}

\subsection{Small-ratio boundary correction}

For $0<\rho<1$, set $D_\rho:=\{\Logp(1+q):|q|<\rho\}\subset S_\pi$.
For $c>0$ and $\theta\in S_\pi$, define
\[
R_c(\theta)=\log D(\theta,c)-c^{-1}\Lambda_\infty(-\theta)
-\frac{\theta^2}{2\pi^2}\log(c/2).
\]
The uniform bound of Proposition~\ref{UC:prop:bound} allows us to identify the small-ratio limit on $S_\pi$ from its trace coefficients.
Near $\theta=0$, write $q=e^\theta-1$ and $\theta=\Logp(1+q)$.
Recall that $\mathcal K_{\Hh}(z,w)=\pi^{-1}\{y+y'-\ii(x-x')\}^{-2}$ for $z=x+\ii y$ and $w=x'+\ii y'$.
For $n\ge2$, set $u_1=u_{n+1}=0$, $y_{n+1}=y_1$, $\boldsymbol u=(u_2,\ldots,u_n)$, and $\ell_n=\max_{1\le j\le n}u_j-\min_{1\le j\le n}u_j$.
Let $G_{A,n}(\boldsymbol u)$ and $G_{Q,n}(\boldsymbol u)$ be the integrals of
$\operatorname{Re}\prod_{j=1}^n
\mathcal K_{\Hh}(u_j+\ii y_j,u_{j+1}+\ii y_{j+1})$
over $\boldsymbol y\in(0,\infty)^n$ and $\boldsymbol y\in(0,1)^n$, respectively.
The densities are defined almost everywhere in $\mathbb{R}^{n-1}$.
With $\delta G_n=G_{Q,n}-G_{A,n}$, put
\[
C_n=\int_{\ell_n\le1}\ell_n\delta G_n\dd\boldsymbol u
+\int_{\ell_n>1}\ell_nG_{Q,n}\dd\boldsymbol u.
\]

\begin{theorem}[Small-ratio boundary correction]\label{BC:thm:local}
There exists a unique $H_*\in\Hol(S_\pi)$ such that, for every $K\Subset S_\pi$, there is $\delta_K>0$ with
\[
 \sup_{\theta\in K}|R_c(\theta)-H_*(\theta)|
 =O_K(c^{\delta_K}),
 \qquad c\downarrow0.
\]
Moreover, $H_*(0)=H_*'(0)=H_*''(0)=0$.
For some $0<\rho_*<1$,
\begin{equation}\label{BC:eq:main}
H_*(\theta)=\frac{\theta^2(\log 2-1)}{2\pi^2}
-\sum_{n=2}^\infty\frac{(-1)^{n-1}}n C_n(e^\theta-1)^n,
\qquad \theta\in D_{\rho_*},
\end{equation}
with absolute and locally uniform convergence.
For every $0<\rho'<\rho_*$, $R_c-H_*=O_{\rho'}(c^2)$ uniformly on $\overline{D_{\rho'}}$ as $c\downarrow0$.
\end{theorem}

Before proving Theorem~\ref{BC:thm:local}, we derive weighted integral bounds for the cyclic densities and an expansion of the relative trace coefficients.

\begin{lemma}\label{BC:lem:bound}
There exist constants $C_*,M_*>0$ such that, for $n\ge2$ and $T\ge1$,
\begin{align}
&\int_{\ell_n\le1}\ell_n|\delta G_n|\dd\boldsymbol u
\le C_*n^{12}M_*^n,\label{BC:eq:localbound}\\
&\int_{\ell_n>T}\ell_n|G_{Q,n}|\dd\boldsymbol u
\le C_*n^{12}M_*^nT^{-2},\label{BC:eq:tail}\\
&\int_{1\le\ell_n\le2}\ell_n|G_{A,n}|\dd\boldsymbol u
\le C_*n^{12}M_*^n.\label{BC:eq:angle}
\end{align}
In particular, $|C_n|\le C_*n^{12}M_*^n$.
\end{lemma}

\begin{proof}
Telescoping the height indicators in $\delta G_n$ gives at most $n$ terms, each with a designated height $y\ge1$.
Cyclic relabeling places this vertex at the horizontal anchor, preserves $\ell_n$, and has unit Jacobian.
The identity $\int_{\Rr}|\mathcal K_{\Hh}(x+\ii y,x'+\ii y')|\dd(x-x')=(y+y')^{-1}$ and the Carleman bound show that the integral operator with kernel $|\mathcal K_{\Hh}|$ has norm at most $\pi$ on $L^2(\Hh)$.
Since $\int_{\Hh}|\mathcal K_{\Hh}(z,w)|^2\ddA(w)=1/(4\pi y^2)$, the integrated absolute cycle with a designated height $y\ge1$ is bounded by $\pi^{n-2}\int_1^\infty(4\pi y^2)^{-1}\dd y$.
Restriction to $\ell_n\le1$ and summation over the designated vertices prove \eqref{BC:eq:localbound}.

On a shell $R\le\ell_n\le2R$, choose vertices with minimum and maximum horizontal coordinates.
Each of the two cyclic paths between them contains an edge of horizontal length at least $R/(n-1)$, with the two directed edges used when $n=2$.
There are at most $n^4$ choices of the vertices and selected edges.
For heights in $(0,1)$, each selected kernel is bounded by $Cn^2R^{-2}$.
Removing the selected edges leaves two paths, possibly isolated vertices when $n=2$.
With $u_1=0$, the remaining $n-2$ edge differences and one relative translation form horizontal coordinates with unit Jacobian; on the shell, the translation has range at most $4R$.
Integrating the edge differences over $\Rr$ using the horizontal-integral identity above gives two Carleman path integrals, bounded in product by $\pi^{n-2}$ on $L^2(0,1)$.
Together with the selected long-edge bounds, the at most $n^4$ choices, and the $O(R^2)$ translation and span factors, we obtain $Cn^8M^nR^{-2}$, and hence
\[
\int_{R\le\ell_n\le2R}\ell_n
\int_{(0,1)^n}
\left|\prod_{j=1}^n
\mathcal K_{\Hh}(u_j+\ii y_j,u_{j+1}+\ii y_{j+1})\right|
\dd\boldsymbol y\dd\boldsymbol u
\le Cn^{12}M^nR^{-2},
\qquad R\ge1.
\]
Dyadic summation proves \eqref{BC:eq:tail}.

For $1\le\ell_n\le2$, the selected kernels are bounded by $Cn^2V_{\mathrm C}(y,z)$, where $V_{\mathrm C}(y,z)=[1+(y+z)^2]^{-1}$.
The remaining horizontal integrations produce powers of the Carleman operator $H_{\mathrm C}$ with kernel $H_{\mathrm C}(y,z)=(y+z)^{-1}$.
The resulting height integrals satisfy
\[
\left|\Tr\bigl(V_{\mathrm C}H_{\mathrm C}^j
V_{\mathrm C}H_{\mathrm C}^{n-2-j}\bigr)\right|
\le\|V_{\mathrm C}\|_2^2\|H_{\mathrm C}\|^{n-2}
\le\frac{\pi^{n-2}}2,
\]
since $\|H_{\mathrm C}\|\le\pi$ and $\|V_{\mathrm C}\|_2^2=1/2$.
Together with the bounded translation range, this proves \eqref{BC:eq:angle}.
By homogeneity, the same bound for the weighted integral of the absolute kernel product holds on every full-height shell $R\le\ell_n\le2R$, independently of $R>0$.
Finally, \eqref{BC:eq:localbound} and \eqref{BC:eq:tail} with $T=1$ give $|C_n|\le C_*n^{12}M_*^n$.
\end{proof}

Put $T=c^{-1}$.
Let $\mathsf A_T$ and $\mathsf B_T$ be the Fourier--Laplace compressions, as in Lemma~\ref{MC:lem:compression}, associated with $[0,T]\times(0,\infty)$ and $[0,T]\times[1,\infty)$, respectively.
Set $\mathsf Q_T=\mathsf A_T-\mathsf B_T$ and $t_n(T)=\Tr(\mathsf Q_T^n-\mathsf A_T^n)$ for $n\ge1$.

\begin{lemma}\label{BC:lem:decomp}
For $n\ge2$ and $T\ge1$,
\begin{equation}\label{BC:eq:J}
t_n(T)
=-\frac{Th_n}{4\pi}-C_n
+\frac{h_{n-1}}{\pi^2}(\log T+1)+e_n(T),
\end{equation}
where $|e_n(T)|\le C_*n^{12}M_*^nT^{-2}$.
\end{lemma}

\begin{proof}
Since $\mathsf B_T$ is trace class, so is $\mathsf Q_T^n-\mathsf A_T^n$.
The estimates for absolute kernel products in the proof of Lemma~\ref{BC:lem:bound} justify the horizontal-overlap formula
$t_n(T)=\int(T-\ell_n)_+\delta G_n\dd\boldsymbol u$.
Fourier integration gives
\[
\int\delta G_n\dd\boldsymbol u
=\frac1{2\pi}\int_0^\infty
\bigl[(1-e^{-2\xi})^n-1\bigr]\dd\xi
=-\frac{h_n}{4\pi}.
\]
Set $a_n=(\log2)^{-1}\int_{1\le\ell_n\le2}\ell_nG_{A,n}\dd\boldsymbol u$.
Homogeneity gives
$\int_{\ell_n>1}\min(T,\ell_n)G_{A,n}\dd\boldsymbol u=a_n(\log T+1)$.
Using $(T-\ell_n)_+=T-\min(T,\ell_n)$ and splitting at $\ell_n=1$, the definition of $C_n$ yields
\[
t_n(T)
=-\frac{Th_n}{4\pi}-C_n+a_n(\log T+1)
+\int_{\ell_n>T}(\ell_n-T)G_{Q,n}\dd\boldsymbol u.
\]
The final integral defines $e_n(T)$, whose bound follows from \eqref{BC:eq:tail}.

Let $G_{B,n}$ be the cyclic density obtained by integrating over heights in $[1,\infty)^n$.
Fourier integration gives $\int G_{B,n}\dd\boldsymbol u=1/(4\pi n)$, so
\[
\frac{T}{4\pi n}-\Tr\mathsf B_T^n
=\int\min(T,\ell_n)G_{B,n}\dd\boldsymbol u.
\]
Dilating the shell $2^j\le\ell_n\le2^{j+1}$ to $1\le\ell_n\le2$ changes the lower height cutoff to $2^{-j}$.
The full-height estimate for absolute kernel products in the proof of Lemma~\ref{BC:lem:bound} therefore permits dominated convergence, giving
\[
\int_{2^j\le\ell_n\le2^{j+1}}
\ell_nG_{B,n}\dd\boldsymbol u
\longrightarrow a_n\log2.
\]
The estimate for cycles with a height at least one bounds the contribution from $\ell_n\le1$.
The absolute shell bound controls the incomplete shell and gives
$T\int_{\ell_n>T}|G_{B,n}|\dd\boldsymbol u=O_n(1)$.
Dyadic summation and Ces\`aro convergence consequently yield
\[
\Tr\mathsf B_T^n
=\frac{T}{4\pi n}-a_n\log T+o_n(\log T),
\qquad T\to\infty.
\]

On the other hand, substitute $\theta=\Logp(1+q)$ into Theorem~\ref{thm:strong-halfline} and compare coefficients with the trace-logarithm expansion for $\mathsf B_T$.
Using $\Logp^2(1+q)=2\sum_{n\ge2}(-1)^nh_{n-1}q^n/n$, locally uniform holomorphic convergence and Cauchy's formula give
\[
\Tr\mathsf B_T^n
=\frac{T}{4\pi n}
-\frac{h_{n-1}}{\pi^2}\log T+O_n(1).
\]
Comparison of the logarithmic coefficients gives $a_n=h_{n-1}/\pi^2$, completing the proof.
\end{proof}

\begin{proof}[Proof of Theorem~\ref{BC:thm:local}]
For fixed $0<c\le1$, Theorem~\ref{MC:thm:critical} gives
\[
\mathcal F_{cL,L}(\theta)
-\frac{\theta^2}{2\pi^2}\log(cL/2)
\longrightarrow A_\infty(\theta)+R_c(\theta)
\]
locally uniformly on $S_\pi$.
Since $L\ge cL\ge2$ for sufficiently large $L$, \eqref{UC:eq:bound} implies
$\sup_{0<c\le1}\sup_{\theta\in K}|R_c(\theta)|<\infty$
for every $K\Subset S_\pi$.

Put $T=c^{-1}$ and $q=e^\theta-1$.
Dilation identifies $D(\theta,1/T)$ with
$\det[\Id-q(\Id+q\mathsf A_T)^{-1}\mathsf B_T]$.
For $|q|<1$, its normalized logarithm satisfies
\[
\log D(\theta,1/T)
=\sum_{n\ge1}\frac{(-1)^{n-1}}n q^nt_n(T).
\]
The series converges absolutely because
$\|\mathsf Q_T^n-\mathsf A_T^n\|_1\le n\|\mathsf B_T\|_1$.
Choose $0<\rho_*<\min\{1,M_*^{-1}\}$.
By \eqref{BC:eq:J} and $t_1(T)=-T/(4\pi)$, the terms proportional to $T$ sum to $T\Lambda_\infty(-\theta)$.
Using also
$\sum_{n\ge2}(-1)^{n-1}h_{n-1}q^n/(n\pi^2)
=-\theta^2/(2\pi^2)$,
we obtain
\[
R_{1/T}(\theta)
=\frac{\theta^2(\log2-1)}{2\pi^2}
-\sum_{n\ge2}\frac{(-1)^{n-1}}n C_nq^n+\sum_{n\ge2}\frac{(-1)^{n-1}}n e_n(T)q^n,
\qquad \theta\in D_{\rho_*}.
\]
The bounds in Lemmas~\ref{BC:lem:bound} and~\ref{BC:lem:decomp} give absolute and locally uniform convergence, with the last series bounded by $O_{\rho'}(T^{-2})$ for $|q|\le\rho'<\rho_*$.
Thus $R_c$ converges locally uniformly on $D_{\rho_*}$ to the function given by the first two terms.
Local boundedness and Montel's theorem give a locally uniformly convergent subsequence along every sequence $c\downarrow0$.
All subsequential limits agree on $D_{\rho_*}$, so the identity theorem yields $R_c\to H_*\in\Hol(S_\pi)$ locally uniformly on $S_\pi$.
Consequently, \eqref{BC:eq:main} holds and $R_c-H_*=O_{\rho'}(c^2)$ on $\overline{D_{\rho'}}$ for every $0<\rho'<\rho_*$.

Equation~\eqref{BC:eq:main} gives $H_*(0)=H_*'(0)=0$.
By the definition of $R_c$, Proposition~\ref{MC:prop:checks}, and \eqref{MC:eq:Delta-small},
\[
R_c''(0)
=\Delta(c)-\frac1{8\pi c}
-\frac1{\pi^2}\log(c/2)
=O(c^2).
\]
Since $R_c\to H_*$ locally uniformly, Cauchy's formula gives $H_*''(0)=0$.
Finally, the first paragraph bounds $R_c-H_*$ locally uniformly on $S_\pi$, and the series representation gives a uniform $O(c^2)$ bound on a fixed neighborhood of zero.
Lemma~\ref{UC:lem:three-circles}, applied with $\eta_c=c$ and $p=2$ for $0<c\le1$, therefore gives the asserted power bound on every compact subset of $S_\pi$.
\end{proof}

\subsection{Large-height constants and subcritical asymptotics}
\label{SB:sec:static}

The uniform remainder estimate relates the fixed-height constants to the small-ratio limit.

\begin{theorem}\label{SB:thm:static-bridge}
With $H_*$ from Theorem~\ref{BC:thm:local},
\begin{equation}\label{SB:eq:static-bridge}
\mathcal C_b(\theta)
-\frac{\theta^2}{2\pi^2}\log(b/2)
\longrightarrow A_\infty(\theta)+H_*(\theta),
\qquad b\to\infty,
\end{equation}
locally uniformly on $S_\pi$.
\end{theorem}

\begin{proof}
Fix $K\Subset S_\pi$.
By Lemma~\ref{SB:lem:boundary-tail}, $\mathcal C_b-\mathcal F_{b,bT}\to0$ uniformly for $b\ge2$ and $\theta\in K$ as $T\to\infty$.
For fixed $T\ge1$, Theorem~\ref{MC:thm:critical} and \eqref{MC:eq:Across} give
\[
 \mathcal F_{b,bT}(\theta)-\frac{\theta^2}{2\pi^2}\log(b/2)
 \longrightarrow A_\infty(\theta)+R_{1/T}(\theta),
 \qquad b\to\infty,
\]
uniformly for $\theta\in K$.
Theorem~\ref{BC:thm:local} gives $R_{1/T}\to H_*$ uniformly on $K$ as $T\to\infty$.
The triangle inequality, taking first $b\to\infty$ and then $T\to\infty$, proves \eqref{SB:eq:static-bridge}.
\end{proof}

\begin{proof}[Proof of Theorem~\ref{thm:intro-moving-cutoff}]
The critical regime and the large-$c$ estimate follow from Theorem~\ref{MC:thm:critical}.
The supercritical regime follows by combining \eqref{MC:eq:supercritical-match} with Theorem~\ref{thm:strong-halfline}.
The small-$c$ matching limit is Theorem~\ref{BC:thm:local} together with \eqref{MC:eq:Across}.
Finally, if $b_L\to\infty$ and $b_L/L\to0$, Lemma~\ref{SB:lem:boundary-tail} gives $\mathcal F_{b_L,L}-\mathcal C_{b_L}\to0$ locally uniformly on $S_\pi$.
Combining this with \eqref{SB:eq:static-bridge} proves the subcritical expansion.
\end{proof}

\section{Sparse count-conditioned geometry}\label{sec:conditioning}

We prove the macroscopic profile and local Poisson limits under
sparse exact-count conditioning stated in
Theorem~\ref{thm:intro-conditioned-poisson}.
Weighted Bergman kernel comparison gives the tilted limits,
and reduced-Palm estimates transfer them to the conditional law.

\subsection{Weighted Bergman kernels and localization}\label{CG:sec:weighted}

Fix $a>0$ and recall $D_{a,L}=[0,L]\times[a,\infty)$ and $N_{a,L}=N_a^\infty(L)$ from Section~\ref{sec:intro}.
Let $K_{a,L}$ denote the integral operator on $L^2(D_{a,L},\ddA)$ with kernel $\mathcal K_{\Hh}$ restricted to $D_{a,L}\times D_{a,L}$.
For $s>0$, define $K_{a,L}^{(s)}=sK_{a,L}(\Id+(s-1)K_{a,L})^{-1}$.
The inverse exists since $K_{a,L}$ is a positive contraction.
This is the kernel on $D_{a,L}$ under the law obtained by weighting configurations by $s^{N_{a,L}}$ and normalizing.

Let $A^2(\Hh)$ denote the Bergman space on $\Hh$.
For $s>0$, set
\[
w_{L,s}=1-(1-s)\mathbf1_{D_{a,L}},
\qquad
w_s=1-(1-s)\mathbf1_{\{y\ge a\}}.
\]
Equip $A^2(\Hh)$ with the weighted norm $\|f\|_w^2=\int_\Hh |f(z)|^2w(z)\ddA(z)$ for $w=w_{L,s},w_s$, and denote the corresponding reproducing kernels by $B_{L,s}$ and $B_s$, respectively.
Set $S_{a,s}(z,w)=sB_s(z,w)$ for $\operatorname{Im}z,\operatorname{Im}w\ge a$.

\begin{lemma}\label{CG:lem:weighted-identities}
For $0<s<1$, the finite-window tilted operator has kernel
\begin{equation}\label{CG:eq:finite-Bergman}
K_{a,L}^{(s)}(z,w)=sB_{L,s}(z,w),
\qquad z,w\in D_{a,L}.
\end{equation}
For $z=x+\ii y$ and $w=x'+\ii y'$ in $\Hh$, the stationary weighted Bergman kernel is
\begin{equation}\label{CG:eq:Bs}
B_s(z,w)=\frac1\pi\int_0^\infty
\frac{\xi e^{-\xi(y+y')}e^{\ii\xi(x-x')}}
{1-(1-s)e^{-2a\xi}}\dd\xi.
\end{equation}
Moreover, $S_{a,s}-K_{a,L}^{(s)}$ is a positive semidefinite kernel on $D_{a,L}$, and
\begin{equation}\label{CG:eq:order}
0\le K_{a,L}^{(s)}
\le\mathbf1_{D_{a,L}}S_{a,s}\mathbf1_{D_{a,L}}
\end{equation}
as operators on $L^2(D_{a,L},\ddA)$.
\end{lemma}

\begin{proof}
Let $H=A^2(\Hh)$, let $U:H\to L^2(\Hh,\ddA)$ be the canonical inclusion, and write $P=\mathbf1_{D_{a,L}}$.
Since $\mathcal K_{\Hh}$ is the Bergman kernel of $\Hh$, we have $K_{a,L}=PUU^*P$.
For $\omega\in\{w_{L,s},w_s\}$, set $G_\omega=U^*M_\omega U$.
Then $s\Id_H\le G_\omega\le\Id_H$, and $UG_\omega^{-1}U^*$ has integral kernel $B_\omega$ with respect to $\ddA$.
Writing $\alpha=1-s$, the Neumann series gives
\[
PUG_{w_{L,s}}^{-1}U^*P
=\sum_{m\ge0}\alpha^mK_{a,L}^{m+1}
=K_{a,L}(\Id-\alpha K_{a,L})^{-1}.
\]
The series converges in operator norm since $0<\alpha<1$ and $\|U^*PU\|\le1$.
Multiplication by $s$ proves \eqref{CG:eq:finite-Bergman}.
For the horizontally invariant weight $w_s$, the Fourier--Laplace representation of $f\in H$ is
$f(x+\ii y)=(2\pi)^{-1/2}\int_0^\infty
g(\xi)e^{\ii\xi x-\xi y}\dd\xi$.
Plancherel and integration in the height variable give
\[
\|f\|_{w_s}^2
=\int_0^\infty |g(\xi)|^2
\frac{1-(1-s)e^{-2a\xi}}{2\xi}\dd\xi.
\]
The corresponding reproducing-kernel formula yields \eqref{CG:eq:Bs}.

Finally, $w_{L,s}\ge w_s$ implies $G_{w_{L,s}}^{-1}\le G_{w_s}^{-1}$.
Applying this inequality to finite linear combinations of evaluation vectors shows that $B_s-B_{L,s}$ is a positive semidefinite kernel.
By \eqref{CG:eq:finite-Bergman} and $S_{a,s}=sB_s$, the same holds for $S_{a,s}-K_{a,L}^{(s)}$ on $D_{a,L}$.
Compressing the inverse operator inequality to $L^2(D_{a,L})$ and multiplying by $s$ gives \eqref{CG:eq:order}.
\end{proof}

Direct integration of \eqref{CG:eq:Bs} gives the stationary horizontal intensity
\begin{equation}\label{CG:eq:x}
x_a(s)=\int_a^\infty S_{a,s}(x+\ii y,x+\ii y)\dd y
=\frac{s\log(1/s)}{4\pi a(1-s)},\qquad 0<s<1.
\end{equation}
By Corollary~\ref{cor:general-halfline-scaling} and \eqref{eq:Lambda-derivatives-halfline}, $x_a(s)=(\Lambda_a^\infty)'(\log s)$.
It is strictly increasing, with $x_a(s)\downarrow0$ as $s\downarrow0$ and $x_a(s)\uparrow1/(4\pi a)$ as $s\uparrow1$.
Set $x_a(1)=1/(4\pi a)$ by continuity.

\begin{lemma}\label{CG:lem:localization}
Let $0<s<1$ and $z=x+\ii y\in D_{a,L}$ with $0<x<L$ and $sy\le a$, and set $b_z=\min(x,L-x)$.
For every integer $m\ge1$ and $R\ge y$,
\begin{equation}\label{CG:eq:diagonal-error}
0\le1-\frac{K_{a,L}^{(s)}(z,z)}{S_{a,s}(z,z)}
\le C_{a,m}\left[
\frac yR+
\frac{y}{s^2b_z^2}
\left(\frac{R+y}{b_z}\right)^{2m}\right].
\end{equation}
The constant $C_{a,m}$ is independent of $L,s,z,R$.
\end{lemma}

\begin{proof}
Expanding the denominator in \eqref{CG:eq:Bs} as a geometric series and using $\sum_{k\ge0}s(1-s)^k=1$, we obtain
$|S_{a,s}(z,w)|\le
[\pi((x-x')^2+(y+y')^2)]^{-1}$
for $w=x'+\ii y'$.
Consequently,
\begin{equation}\label{CG:eq:leakage}
\int_{\{x'\notin[0,L],\,y'\ge a\}}
|S_{a,s}(w,z)|^2\ddA(w)
\le\frac1{\pi^2}\int_{|t|\ge b_z}
\int_0^\infty\frac{\dd v\dd t}{(t^2+v^2)^2}
=\frac1{4\pi b_z^2}.
\end{equation}
We also have $B_s(z,z)\ge c/(ay)$ for an absolute constant $c>0$ whenever $sy\le a$.
Indeed, on $1/(2y)\le\xi\le1/y$, $1-(1-s)e^{-2a\xi}\le s+2a\xi\le3a/y$, and the bound follows from \eqref{CG:eq:Bs}.

For $R\ge y$, define the bounded holomorphic multiplier
\[
\psi_z(\zeta)=
\left(\frac{R+y}{R-\ii(\zeta-x)}\right)^m.
\]
Then $\psi_z(z)=1$, $\|\psi_z\|_{L^\infty(\Hh)}\le(1+y/R)^m$, and $|\psi_z(x'+\ii y')|\le((R+y)/b_z)^m$ whenever $x'\notin[0,L]$.
Set $f_z=\psi_zB_s(\cdot,z)\in A^2(\Hh)$.
The reproducing property gives $\|B_s(\cdot,z)\|_{w_s}^2=B_s(z,z)$.
Since $w_{L,s}-w_s=(1-s)\mathbf1_{\{x'\notin[0,L],\,y'\ge a\}}$, the estimate \eqref{CG:eq:leakage} yields
\[
\|f_z\|_{w_{L,s}}^2
\le(1+y/R)^{2m}B_s(z,z)
+\frac1{4\pi s^2b_z^2}
\left(\frac{R+y}{b_z}\right)^{2m}.
\]
As $f_z(z)=B_s(z,z)$, the variational characterization of the reproducing kernel gives
\[
\frac{K_{a,L}^{(s)}(z,z)}{S_{a,s}(z,z)}
=\frac{B_{L,s}(z,z)}{B_s(z,z)}
\ge\left[
(1+y/R)^{2m}
+\frac1{4\pi s^2b_z^2B_s(z,z)}
\left(\frac{R+y}{b_z}\right)^{2m}
\right]^{-1}.
\]
The diagonal lower bound, $(1+y/R)^{2m}-1\le C_my/R$, and $1-(1+u)^{-1}\le u$ give the upper bound in \eqref{CG:eq:diagonal-error}.
Its nonnegativity follows from the positive-semidefinite kernel order in Lemma~\ref{CG:lem:weighted-identities}.
\end{proof}

The localization estimate yields uniform diagonal and local trace-norm comparisons with the stationary kernel away from the vertical sides of the observation window.

\begin{proposition}\label{CG:prop:uniform-comparison}
Let $s_L\in(0,1)$ satisfy $s_L\to0$ and $Lx_a(s_L)\to\infty$.
For every fixed $\delta\in(0,1/2)$ and $\beta\in(0,1)$,
\begin{equation}\label{CG:eq:relative-diagonal}
\sup_{\substack{\delta L\le x\le(1-\delta)L\\
a\le y\le as_L^{-\beta}}}
\left|
\frac{K_{a,L}^{(s_L)}(x+\ii y,x+\ii y)}
{S_{a,s_L}(x+\ii y,x+\ii y)}-1
\right|\longrightarrow0.
\end{equation}
Moreover, for every measurable subset $\mathcal U_L$ of this region,
\begin{equation}\label{CG:eq:trace-norm-comparison}
\bigl\|\mathbf1_{\mathcal U_L}
(S_{a,s_L}-K_{a,L}^{(s_L)})
\mathbf1_{\mathcal U_L}\bigr\|_1
\le o(1)\int_{\mathcal U_L}S_{a,s_L}(z,z)\ddA(z),
\end{equation}
where the $o(1)$ is uniform over all such $\mathcal U_L$.
\end{proposition}

\begin{proof}
Let $\varepsilon_L$ denote the supremum in \eqref{CG:eq:relative-diagonal}.
Set $\ell_L=\log(1/s_L)$ and $\kappa=(1-\beta)/2$, and choose an integer $m\ge1$ such that $m(1-\beta)>\beta$.
By \eqref{CG:eq:x}, the hypothesis implies $(Ls_L)^{-1}=o(\ell_L)$, so $(Ls_L)^{-1}\le\ell_L$ for all sufficiently large $L$.
Set $R=Ls_L^\kappa$.
Uniformly over the region in \eqref{CG:eq:relative-diagonal}, we have $s_Ly\le a$, $b_z\ge\delta L$, and $y/R\le a\ell_Ls_L^\kappa\to0$.
Thus $R\ge y$ and $R+y\le2R$ for sufficiently large $L$.
Applying Lemma~\ref{CG:lem:localization}, the second term in \eqref{CG:eq:diagonal-error} is bounded by $C_{a,m,\delta}\ell_L^2s_L^{m(1-\beta)-\beta}$.
Consequently,
\[
\varepsilon_L\le C_{a,m,\delta}
\left(\ell_Ls_L^\kappa+
\ell_L^2s_L^{m(1-\beta)-\beta}\right)=o(1),
\]
proving \eqref{CG:eq:relative-diagonal}.

For every measurable $\mathcal U_L$ as in the statement, the stationary kernel $S_{a,s_L}$ is positive definite and has integrable diagonal on $\mathcal U_L$, so its compression is trace class.
By \eqref{CG:eq:order}, the compressed difference is positive and trace class.
Its trace norm therefore equals the integral of its diagonal, and \eqref{CG:eq:relative-diagonal} gives \eqref{CG:eq:trace-norm-comparison}.
\end{proof}

\subsection{Stationary height distribution and exact saddle}
\label{CG:sec:saddle}

For $0<s<1$ and $\xi>0$, the horizontal Fourier fibre of $S_{a,s}$ on $L^2([a,\infty))$ has rank one, with normalized eigenfunction $v_\xi(y)=\sqrt{2\xi}e^{-\xi(y-a)}$ and eigenvalue
\[
q_s(\xi)=\frac{s e^{-2a\xi}}{1-(1-s)e^{-2a\xi}}\in(0,1).
\]
The fibres vanish for $\xi<0$, so $S_{a,s}$ is a positive contraction.
Its norm equals one since $q_s(\xi)\to1$ as $\xi\downarrow0$.
Its diagonal is
\begin{equation}\label{CG:eq:rho}
\rho_{a,s}(y)=S_{a,s}(x+\ii y,x+\ii y)
=\frac{s}{4\pi}\sum_{k\ge0}\frac{(1-s)^k}{(y+ak)^2},
\qquad y\ge a.
\end{equation}
By \eqref{CG:eq:x}, $\int_a^\infty\rho_{a,s}(y)\dd y=x_a(s)$.
Let $Y_s$ have density $\rho_{a,s}(y)/x_a(s)$ on $[a,\infty)$.

\begin{theorem}[Logarithmic height distribution]\label{CG:thm:height-escape}
As $s\downarrow0$,
\[
\sup_{0\le t\le1}
\left|\Pp\{\log(Y_s/a)/\log(1/s)>t\}-(1-t)\right|
=O(1/\log(1/s)).
\]
In particular, $\log(Y_s/a)/\log(1/s)\Longrightarrow\operatorname{Unif}(0,1)$.
\end{theorem}

\begin{proof}
Integration of \eqref{CG:eq:rho} and \eqref{CG:eq:x} gives, for $R\ge1$ and $\ell=\log(1/s)$,
\begin{equation}\label{CG:eq:height-survival}
\Pp\{Y_s>aR\}=\frac{1-s}{\ell}
\sum_{k\ge0}\frac{(1-s)^k}{R+k}.
\end{equation}
Put $M=\lfloor s^{-1}\rfloor$.
Since $0\le1-(1-s)^k\le sk$ and $\sum_{k>M}(1-s)^k/(R+k)\le[s(R+M+1)]^{-1}$, we obtain, uniformly for $1\le R\le s^{-1}$,
\[
\sum_{k\ge0}\frac{(1-s)^k}{R+k}
=\sum_{k=0}^M\frac1{R+k}+O(1)
=\log\frac1{sR}+O(1).
\]
The last equality follows by sum--integral comparison, since $s(R+M+1)$ is bounded above and below by positive constants.
Taking $R=s^{-t}$, $0\le t\le1$, gives
\[
\Pp\{\log(Y_s/a)/\ell>t\}
=(1-s)(1-t)+O(\ell^{-1})
=1-t+O(\ell^{-1})
\]
uniformly in $t$, as required.
\end{proof}

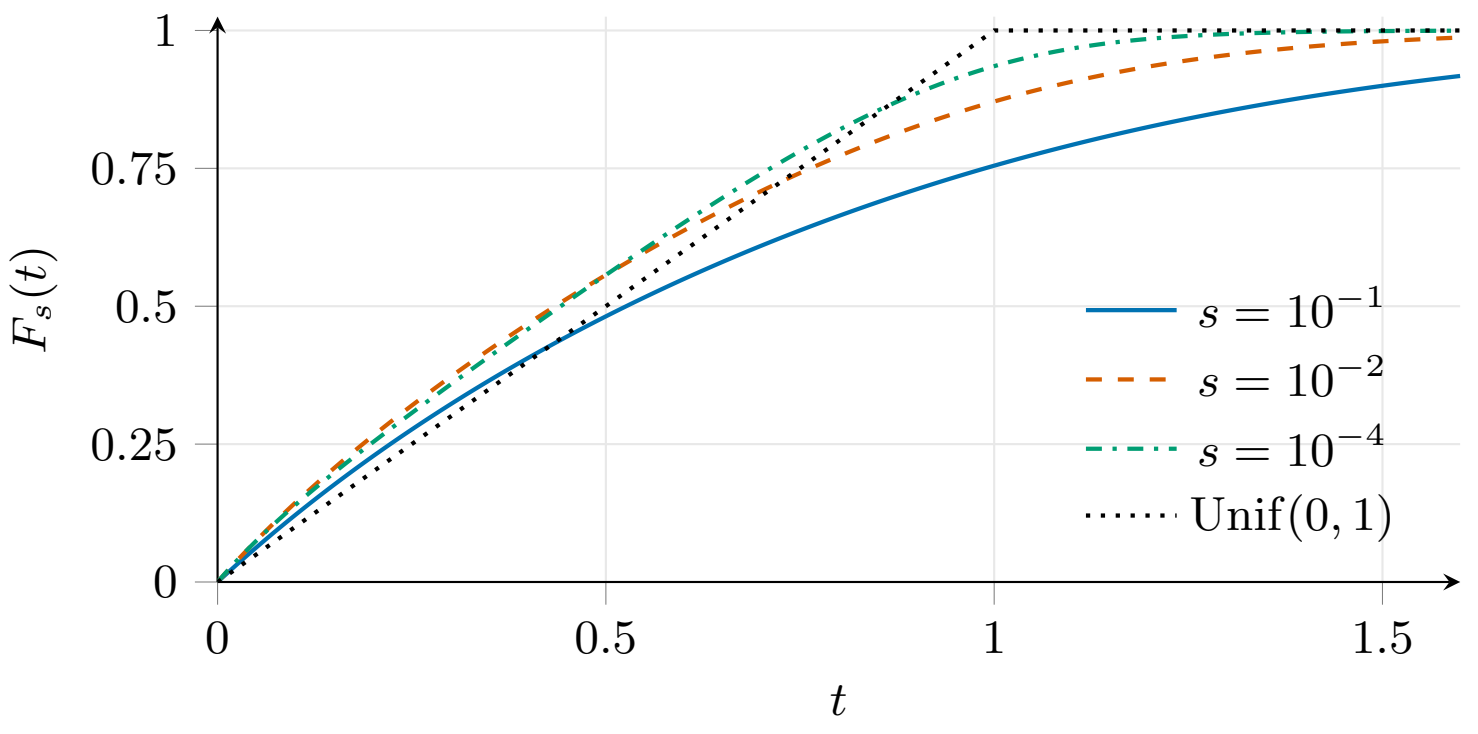
\begin{figure}[!htbp]
\centering
\begin{tikzpicture}
\ExplSyntaxOn
\fp_new:N \g_hs_s_fp
\fp_new:N \g_hs_q_fp
\fp_new:N \g_hs_log_fp
\fp_new:N \g_hs_alpha_fp
\fp_new:N \g_hs_pref_fp
\fp_new:N \g_hs_p_fp
\fp_gset:Nn \g_hs_s_fp {-1}
\fp_new:N \l_hs_r_fp
\fp_new:N \l_hs_v_fp
\fp_new:N \l_hs_z_fp
\fp_new:N \l_hs_e_fp
\fp_new:N \l_hs_w_fp
\fp_new:N \l_hs_f_fp
\int_new:N \l_hs_k_int
\tl_new:N \l_hs_cf_tl

\int_step_inline:nn {44}
  {
    \tl_set:Nn \l_hs_cf_tl {0}
    \int_step_inline:nnnn {#1} {-1} {1}
      {
        \tl_set:Nx \l_hs_cf_tl
          {\int_eval:n {##1*##1}/
            (\exp_not:N \l_hs_z_fp+\int_eval:n {2*##1+1}
              -\exp_not:V \l_hs_cf_tl)}
      }
    \cs_new:cpx {hs_height_cf_#1:}
      {1/(\exp_not:N \l_hs_z_fp+1-\exp_not:V \l_hs_cf_tl)}
  }

\cs_new_protected:Npn \hs_height_parameters:n #1
  {
    \fp_compare:nNnF {#1} = {\g_hs_s_fp}
      {
        \fp_gset:Nn \g_hs_s_fp {#1}
        \fp_gset:Nn \g_hs_q_fp {1-\g_hs_s_fp}
        \fp_gset:Nn \g_hs_log_fp {ln(\g_hs_s_fp)}
        \fp_compare:nNnTF {\g_hs_s_fp} > {0.1}
          {\fp_gset:Nn \g_hs_alpha_fp {-ln(\g_hs_q_fp)}}
          {
            \fp_gset:Nn \g_hs_alpha_fp {\g_hs_s_fp/(2-\g_hs_s_fp)}
            \fp_set:Nn \l_hs_w_fp {\g_hs_alpha_fp*\g_hs_alpha_fp}
            \fp_gset:Nn \g_hs_alpha_fp
              {2*\g_hs_alpha_fp*(1+\l_hs_w_fp*(1/3+\l_hs_w_fp*(1/5+
                \l_hs_w_fp*(1/7+\l_hs_w_fp*(1/9+\l_hs_w_fp*(1/11+
                \l_hs_w_fp/13))))))}
          }
        \fp_gset:Nn \g_hs_pref_fp {\g_hs_q_fp/\g_hs_log_fp}
        \fp_gset:Nn \g_hs_p_fp {\g_hs_q_fp*\g_hs_q_fp}
        \fp_gset:Nn \g_hs_p_fp {\g_hs_p_fp*\g_hs_p_fp}
        \fp_gset:Nn \g_hs_p_fp {\g_hs_p_fp*\g_hs_p_fp}
        \fp_gset:Nn \g_hs_p_fp {\g_hs_p_fp*\g_hs_p_fp}
      }
  }

\cs_new_protected:Npn \HSHeightCompute #1#2
  {
    \fp_compare:nNnTF {#1} > {0}
      {
        \hs_height_parameters:n {#2}
        \fp_set:Nn \l_hs_r_fp {exp(-(#1)*\g_hs_log_fp)}
        \fp_set:Nn \l_hs_v_fp {1/(\l_hs_r_fp+16)}
        \fp_set:Nn \l_hs_z_fp {\g_hs_alpha_fp/\l_hs_v_fp}
        \fp_compare:nNnTF {\l_hs_z_fp} > {2}
          {
            \int_set:Nn \l_hs_k_int
              {\fp_eval:n {ceil(80/\l_hs_z_fp)+4}}
            \fp_set:Nn \l_hs_e_fp
              {\use:c {hs_height_cf_\int_use:N \l_hs_k_int :}}
          }
          {
\fp_set:Nn \l_hs_e_fp {exp(\l_hs_z_fp)*(-0.5772156649015329-ln(\l_hs_z_fp)+\l_hs_z_fp*(1-1*\l_hs_z_fp/4*
  (1-2*\l_hs_z_fp/9*(1-3*\l_hs_z_fp/16*(1-4*\l_hs_z_fp/25*(1-5*\l_hs_z_fp/36*(1-6*\l_hs_z_fp/49*(1-7*
  \l_hs_z_fp/64*(1-8*\l_hs_z_fp/81*(1-9*\l_hs_z_fp/100*(1-10*\l_hs_z_fp/121*(1-11*\l_hs_z_fp/144*(1-12*
  \l_hs_z_fp/169*(1-13*\l_hs_z_fp/196*(1-14*\l_hs_z_fp/225*(1-15*\l_hs_z_fp/256*(1-16*\l_hs_z_fp/289*(1-
  17*\l_hs_z_fp/324*(1-18*\l_hs_z_fp/361*(1-19*\l_hs_z_fp/400*(1-20*\l_hs_z_fp/441*(1-21*\l_hs_z_fp/484*
  (1-22*\l_hs_z_fp/529*(1-23*\l_hs_z_fp/576*1))))))))))))))))))))))))}
          }
        \fp_set:Nn \l_hs_w_fp {\l_hs_v_fp*\l_hs_v_fp}
        \fp_set:Nn \l_hs_f_fp
          {1+\g_hs_pref_fp*(
            1/\l_hs_r_fp+\g_hs_q_fp*(1/(\l_hs_r_fp+1)+\g_hs_q_fp*(1/(\l_hs_r_fp+2)+\g_hs_q_fp*(1/
              (\l_hs_r_fp+3)+\g_hs_q_fp*(1/(\l_hs_r_fp+4)+\g_hs_q_fp*(1/(\l_hs_r_fp+5)+\g_hs_q_fp*(1/
              (\l_hs_r_fp+6)+\g_hs_q_fp*(1/(\l_hs_r_fp+7)+\g_hs_q_fp*(1/(\l_hs_r_fp+8)+\g_hs_q_fp*(1/
              (\l_hs_r_fp+9)+\g_hs_q_fp*(1/(\l_hs_r_fp+10)+\g_hs_q_fp*(1/(\l_hs_r_fp+11)+\g_hs_q_fp*(1/
              (\l_hs_r_fp+12)+\g_hs_q_fp*(1/(\l_hs_r_fp+13)+\g_hs_q_fp*(1/(\l_hs_r_fp+14)+\g_hs_q_fp*(1/
              (\l_hs_r_fp+15))))))))))))))))
            +\g_hs_p_fp*(\l_hs_e_fp+\l_hs_v_fp/2+\l_hs_w_fp*((1*\l_hs_z_fp+1)/12+\l_hs_w_fp*(-(((1*
              \l_hs_z_fp+3)*\l_hs_z_fp+6)*\l_hs_z_fp+6)/720+\l_hs_w_fp*((((((1*\l_hs_z_fp+5)*\l_hs_z_fp+
              20)*\l_hs_z_fp+60)*\l_hs_z_fp+120)*\l_hs_z_fp+120)/30240+\l_hs_w_fp*(-(((((((1*\l_hs_z_fp+
              7)*\l_hs_z_fp+42)*\l_hs_z_fp+210)*\l_hs_z_fp+840)*\l_hs_z_fp+2520)*\l_hs_z_fp+5040)*
              \l_hs_z_fp+5040)/1209600))))))}
      }
      {\fp_zero:N \l_hs_f_fp}
    \cs_set:Npx \HSresult {\fp_to_decimal:N \l_hs_f_fp}
  }
\ExplSyntaxOff

\pgfmathdeclarefunction{hsheight}{2}{%
  \begingroup
  \pgfkeys{/pgf/fpu=true,/pgf/fpu/output format=sci}%
  \pgfmathsetmacro\HSt{#1}%
  \pgfmathsetmacro\HSs{#2}%
  \HSHeightCompute{\HSt}{\HSs}%
  \edef\HSnext{\endgroup\noexpand\pgfmathparse{\HSresult}}%
  \HSnext
}

\begin{axis}[
  width=0.70\textwidth,
  height=0.38\textwidth,
  axis lines=left,
  tick align=outside,
  tick label style={font=\small},
  label style={font=\small},
  xlabel={$t$},
  ylabel={$F_s(t)$},
  xmin=0,xmax=1.6,
  ymin=0,ymax=1.025,
  xtick={0,0.5,1,1.5},
  ytick={0,0.25,0.5,0.75,1},
  legend style={
    font=\small,
    draw=none,
    fill=none,
    at={(0.97,0.04)},
    anchor=south east
  },
  grid=major,
  grid style={black!9},
  enlargelimits=false,
  domain=0:1.6,
  samples=201
]
\addplot[HSblue,thick] {hsheight(x,0.1)};
\addlegendentry{$s=10^{-1}$}

\addplot[HSorange,thick,dashed] {hsheight(x,0.01)};
\addlegendentry{$s=10^{-2}$}

\addplot[HSgreen,thick,dashdotted] {hsheight(x,0.0001)};
\addlegendentry{$s=10^{-4}$}

\addplot[black,thick,dotted,samples at={0,1,1.6}] {min(x,1)};
\addlegendentry{$\mathrm{Unif}(0,1)$}
\end{axis}
\end{tikzpicture}
\caption{Stationary tilted logarithmic-height distributions and their $\operatorname{Unif}(0,1)$ limit.}
\label{fig:log-height-distribution}
\end{figure}

For $0<s<1$, taking $R=s^{-t}$ in \eqref{CG:eq:height-survival} gives
\[
F_s(t)
=\Pp\!\left\{\frac{\log(Y_s/a)}{\log(1/s)}\le t\right\}
=1-\frac{1-s}{\log(1/s)}
\sum_{k=0}^{\infty}\frac{(1-s)^k}{s^{-t}+k},
\qquad t\ge0.
\]
Figure~\ref{fig:log-height-distribution} compares the distribution functions $F_s(t)$ with their $\operatorname{Unif}(0,1)$ limit.
The functions $F_s$ are independent of $a$ and satisfy $F_s(t)<1$ for $t>1$.

\FloatBarrier

\begin{lemma}\label{CG:lem:square-trace}
Let $0<s<1$ and let $B\subset\Rr\times[a,\infty)$ be measurable with horizontal projection of Lebesgue measure at most $b<\infty$.
Then
\[
\Tr[(\mathbf1_BS_{a,s}\mathbf1_B)^2]
\le\frac b{2\pi}\int_0^\infty q_s(\xi)^2\dd\xi
\le\frac{bs}{4\pi a}.
\]
\end{lemma}

\begin{proof}
Write $S=S_{a,s}$ and $P=\mathbf1_B$.
Plancherel and the normalization of the Fourier eigenfunctions $v_\xi$ give
\[
\|PS\|_2^2
=\int_B\int_{\Rr\times[a,\infty)}
|S(z,w)|^2\ddA(w)\ddA(z)
\le\frac b{2\pi}\int_0^\infty q_s(\xi)^2\dd\xi.
\]
Since $PSP$ is self-adjoint, $\Tr[(PSP)^2]=\|PSP\|_2^2\le\|PS\|_2^2$.
Finally, $q_s(\xi)=s/(s+e^{2a\xi}-1)\le s/(s+2a\xi)$, and integration of its square gives $\int_0^\infty q_s(\xi)^2\dd\xi\le s/(2a)$.
\end{proof}

The positive-contraction fibre representation above and the locally integrable diagonal of $S_{a,s}$ imply that $S_{a,s}$ is locally trace class.
The determinantal existence criterion (see, for example, \cite[Theorem~1.2]{ShiraiTakahashi2003}) therefore yields a determinantal point process $\mathcal Z_{a,s}$ on $\Rr\times[a,\infty)$ with kernel $S_{a,s}$.
The process is stationary under horizontal translations.

\Needspace{9\baselineskip}
\begin{proposition}\label{CG:prop:stationary-PPP}
For $0<s<1$, define
\[
\mathcal M_{a,s}
:=\left.\sum_{x+\ii y\in\mathcal Z_{a,s}}
\delta_{(x_a(s)x,\,\log(y/a)/\log(1/s))}\right|_{\Rr\times(0,1)}.
\]
Then, as $s\downarrow0$,
\[
\mathcal M_{a,s}\Longrightarrow\PPP(\dd u\,\dd v)
\quad\text{vaguely on }\Rr\times(0,1).
\]
\end{proposition}

\begin{proof}
By horizontal stationarity, the intensity measure of $\mathcal M_{a,s}$ is the restriction of
\[
 \dd u\,\Law\!\left(\frac{\log(Y_s/a)}{\log(1/s)}\right)(\dd v)
\]
to $\Rr\times(0,1)$.
Theorem~\ref{CG:thm:height-escape} therefore gives vague convergence of the intensity measures to $\dd u\,\dd v$.
For $f\in C_c^+(\Rr\times(0,1))$, let $E_s$ be the inverse image of $\operatorname{supp}f$ under the rescaling in the statement.
Its horizontal projection has length $O_f(x_a(s)^{-1})$.
Lemma~\ref{CG:lem:square-trace} and \eqref{CG:eq:x} give
\[
\Tr[(\mathbf1_{E_s}S_{a,s}\mathbf1_{E_s})^2]
=O_f(s/x_a(s))
=O_f(1/\log(1/s))\longrightarrow0.
\]

Set $\ell=\log(1/s)$ and $g_s(x+\ii y)=1-e^{-f(x_a(s)x,\,\log(y/a)/\ell)}$, where $f$ is extended by zero.
The positive trace-class operator $A_s=\sqrt{g_s}\,S_{a,s}\sqrt{g_s}$ satisfies $\Tr A_s\to\int_{\Rr\times(0,1)}(1-e^{-f})\dd u\dd v$ and
$\Tr A_s^2\le
\Tr[(\mathbf1_{E_s}S_{a,s}\mathbf1_{E_s})^2]\to0$.
Thus $\|A_s\|\to0$ and $\log\det(\Id-A_s)=-\Tr A_s+o(1)$.
The determinantal Laplace-functional identity now yields the Poisson Laplace functional, proving the assertion.
\end{proof}

For $s>0$, put $\mu_L(s)=\Tr K_{a,L}^{(s)}$ and $v_L(s)=\Tr\bigl(K_{a,L}^{(s)}-(K_{a,L}^{(s)})^2\bigr)$.

\begin{lemma}\label{CG:lem:spectral-basics}
The operator $K_{a,L}$ is positive trace class with $\Tr K_{a,L}=L/(4\pi a)$, has infinite rank, and satisfies $\|K_{a,L}\|<1$.
Moreover, $\Pp\{N_{a,L}=n\}>0$ for every integer $n\ge0$.
For each integer $n\ge1$, the equation $\mu_L(s)=n$ has a unique solution $s>0$.
\end{lemma}

\begin{proof}
Let $H=A^2(\Hh)$ and let $R:H\to L^2(D_{a,L},\ddA)$ be the restriction operator.
Since $\mathcal K_{\Hh}$ is the Bergman kernel of $\Hh$, we have $K_{a,L}=RR^*$ and
\[
\Tr K_{a,L}=\|R\|_2^2
=\int_{D_{a,L}}\mathcal K_{\Hh}(z,z)\ddA(z)
=\frac{L}{4\pi a}.
\]
The identity theorem implies that $R$ is injective, so $K_{a,L}$ has infinite rank.
If $\|K_{a,L}\|=1$, compactness of $R$ would give a nonzero $f\in H$ with $\|Rf\|=\|f\|$.
Then $f$ would vanish outside $D_{a,L}$, contradicting the identity theorem.
Thus $\|K_{a,L}\|<1$.

Let $(\lambda_j)$ be the positive eigenvalues of $K_{a,L}$.
The count $N_{a,L}$ is distributed as a sum of independent Bernoulli$(\lambda_j)$ variables.
Since $\sum_j\lambda_j<\infty$ and $\sup_j\lambda_j<1$, we have $\prod_j(1-\lambda_j)>0$.
The event that exactly the first $n$ Bernoulli variables equal one therefore has positive probability for every $n\ge0$.
Finally,
\[
\mu_L(s)=\sum_j\frac{s\lambda_j}
{1-\lambda_j+s\lambda_j}
\]
is continuous and strictly increasing on $(0,\infty)$.
The trace-class property and $\|K_{a,L}\|<1$ imply $\mu_L(s)\to0$ as $s\downarrow0$, while infinite rank gives $\mu_L(s)\to\infty$ as $s\to\infty$.
Hence $\mu_L(s)=n$ has a unique solution for every $n\ge1$.
\end{proof}

We next establish the mean and variance asymptotics for the tilted count.

\begin{proposition}\label{CG:prop:mean-variance}
Let $s_L\in(0,1)$ satisfy $s_L\to0$ and $Lx_a(s_L)\to\infty$.
Then
\[
\mu_L(s_L)\sim Lx_a(s_L),\qquad
v_L(s_L)\sim Lx_a(s_L).
\]
For every $0<s<1$,
\begin{equation}\label{CG:eq:mean-upper}
\mu_L(s)\le Lx_a(s),\qquad
\Tr[(K_{a,L}^{(s)})^2]\le\frac{Ls}{4\pi a}.
\end{equation}
\end{proposition}

\begin{proof}
The mean upper bound follows from \eqref{CG:eq:order} and \eqref{CG:eq:x}.
For the reverse asymptotic bound, write $s=s_L$ and fix $\delta\in(0,1/2)$ and $0<\alpha<\beta<1$.
Proposition~\ref{CG:prop:uniform-comparison} and Theorem~\ref{CG:thm:height-escape} give
\[
\mu_L(s)
\ge\int_{\delta L}^{(1-\delta)L}
\int_{as^{-\alpha}}^{as^{-\beta}}
K_{a,L}^{(s)}(x+\ii y,x+\ii y)\dd y\dd x
\ge(1-o(1))(1-2\delta)Lx_a(s)
(\beta-\alpha+o(1)).
\]
Letting $L\to\infty$ and then $\delta,\alpha\downarrow0$ and $\beta\uparrow1$ proves $\mu_L(s)\sim Lx_a(s)$.

By \eqref{CG:eq:order}, the eigenvalue monotonicity of positive compact operators, and Lemma~\ref{CG:lem:square-trace},
\[
\Tr[(K_{a,L}^{(s)})^2]
\le\Tr[(\mathbf1_{D_{a,L}}S_{a,s}\mathbf1_{D_{a,L}})^2]
\le\frac{Ls}{4\pi a}.
\]
Since $x_a(s)\sim s\log(1/s)/(4\pi a)$, this bound is $o(\mu_L(s))$.
Hence $v_L(s)\sim\mu_L(s)\sim Lx_a(s)$.
\end{proof}

\begin{proposition}\label{CG:prop:saddle-asymptotics}
Let $n_L$ be positive integers with $n_L\to\infty$ and $d_L:=n_L/L\to0$, and set $h_L=\log(1/d_L)$.
For all sufficiently large $L$, let $s_L\in(0,1)$ be the unique solution of $x_a(s_L)=d_L$, and let $\widehat s_L>0$ be the unique exact saddle satisfying $\mu_L(\widehat s_L)=n_L$.
Then
\[
\frac{\widehat s_L}{s_L}\longrightarrow1,\qquad
s_L\sim\frac{4\pi a n_L}{L\log(L/n_L)},\qquad
v_L(\widehat s_L)\sim n_L.
\]
Moreover, $\log(1/\widehat s_L)/h_L\to1$.
\end{proposition}

\begin{proof}
Since $x_a$ is strictly increasing and $d_L\to0$, the solution $s_L$ exists uniquely for all large $L$ and satisfies $s_L\to0$ and $Lx_a(s_L)=n_L\to\infty$.
By \eqref{CG:eq:mean-upper}, $\mu_L(s_L)\le n_L$, hence $\widehat s_L\ge s_L$.
For each fixed $\varepsilon>0$, Proposition~\ref{CG:prop:mean-variance} and \eqref{CG:eq:x} give
\[
\frac{\mu_L((1+\varepsilon)s_L)}{n_L}
\sim\frac{x_a((1+\varepsilon)s_L)}{x_a(s_L)}
\longrightarrow1+\varepsilon.
\]
Monotonicity of $\mu_L$ implies $s_L\le\widehat s_L\le(1+\varepsilon)s_L$ for all large $L$, so $\widehat s_L/s_L\to1$.
Consequently, $Lx_a(\widehat s_L)\sim n_L\to\infty$, and Proposition~\ref{CG:prop:mean-variance} gives $v_L(\widehat s_L)\sim n_L$.
Let $\ell_L=\log(1/s_L)$.
The equation $x_a(s_L)=d_L$ gives $h_L=\ell_L-\log\ell_L+\log(4\pi a)+\log(1-s_L)$, hence $\ell_L/h_L\to1$.
Therefore,
\[
s_L=\frac{4\pi a(1-s_L)d_L}{\ell_L}
\sim\frac{4\pi a d_L}{h_L}.
\]
Together with $\widehat s_L/s_L\to1$, this also yields $\log(1/\widehat s_L)/h_L\to1$.
\end{proof}

\subsection{Tilted limits and exact conditioning}
\label{CG:sec:conditioning}

Define the exact-saddle tilted law by
\begin{equation}\label{CG:eq:tilted-law}
\frac{\dd\Q_L}{\dd\Pp}
=\frac{\widehat s_L^{N_{a,L}}}{\Ee\widehat s_L^{N_{a,L}}}.
\end{equation}
Its restriction to $D_{a,L}$ is determinantal with kernel $K_{a,L}^{(\widehat s_L)}$.
We use the same formula for $\mathfrak M_L$ without conditioning, so its total mass is $N_{a,L}/n_L$.

\begin{proposition}\label{CG:prop:tilted-Poisson}
Under $\Q_L$, for every bounded continuous real-valued function $f$ on $[0,1]\times[0,\infty)$,
\begin{equation}\label{CG:eq:tilted-global-L2}
\Ee_{\Q_L}\left|
\langle\mathfrak M_L,f\rangle
-\int_0^1\!\int_0^1 f(u,v)\dd u\dd v
\right|^2\longrightarrow0.
\end{equation}
The joint local convergence in \eqref{eq:intro-conditioned-poisson} also holds under $\Q_L$, with convergence of all fixed mixed factorial moments on relatively compact continuity sets.
\end{proposition}

\begin{proof}
Set $s=\widehat s_L$, $d=n_L/L$, and $\ell=\log(1/s)$.
Proposition~\ref{CG:prop:saddle-asymptotics} gives $s\to0$, $x_a(s)/d\to1$, and $\ell/h_L\to1$.
Let $\nu_L=\Ee_{\Q_L}\mathfrak M_L$, and let $\overline\nu_L$ be the rescaled stationary intensity, the product of $(x_a(s)/d)\mathbf1_{[0,1]}(u)\dd u$ and the law of $\log(Y_s/a)/h_L$.
By \eqref{CG:eq:order} and $\mu_L(s)=n_L$,
\begin{equation}\label{CG:eq:tilted-intensity-TV}
0\le\nu_L\le\overline\nu_L,\qquad
\|\overline\nu_L-\nu_L\|_{\mathrm{TV}}
=\frac{Lx_a(s)-n_L}{n_L}\longrightarrow0.
\end{equation}
Theorem~\ref{CG:thm:height-escape} and \eqref{CG:eq:tilted-intensity-TV} imply $\nu_L\Longrightarrow\mathbf1_{[0,1]^2}(u,v)\dd u\dd v$.
Together with the determinantal variance bound
$\operatorname{Var}_{\Q_L}\langle\mathfrak M_L,f\rangle\le\|f\|_\infty^2/n_L$,
this proves \eqref{CG:eq:tilted-global-L2}.

Fix compact rectangles $C_j\Subset E$, one in each local observation space, and put $C=\bigsqcup_{j=1}^r C_j$.
Their inverse images under the local rescalings are disjoint for large $L$, since $d^{-1}=o(L)$.
Their union $\mathcal E_L$ has horizontal projection of length $O_C(d^{-1})$ and lies in $[\delta L,(1-\delta)L]\times[a,as^{-\beta}]$ for some fixed $\delta\in(0,1/2)$ and $\beta\in(0,1)$.
Its stationary mass is $O_C(x_a(s)/d)=O_C(1)$, so Proposition~\ref{CG:prop:uniform-comparison} gives
\begin{equation}\label{CG:eq:tilted-local-trace}
\bigl\|\mathbf1_{\mathcal E_L}
(K_{a,L}^{(s)}-S_{a,s})\mathbf1_{\mathcal E_L}\bigr\|_1
\longrightarrow0.
\end{equation}
Let $\widetilde K_L$ be the compression of $K_{a,L}^{(s)}$ to $\mathcal E_L$, transported unitarily to $L^2(C)$.
By Theorem~\ref{CG:thm:height-escape} and \eqref{CG:eq:tilted-local-trace}, its intensity converges to $\dd u\,\dd v$ in each local observation space.
Eigenvalue monotonicity, \eqref{CG:eq:order}, and Lemma~\ref{CG:lem:square-trace} give
\begin{equation}\label{CG:eq:tilted-square-trace}
\Tr\widetilde K_L=O_C(1),\qquad
\Tr\widetilde K_L^2
=O_C(s/d)=O_C(\ell^{-1})\longrightarrow0.
\end{equation}
The square-trace bound includes all cross terms between distinct windows.
For nonnegative $f_j\in C_c(E)$ supported in $C_j$, let $F=f_j$ on the $j$th component of $C$ and set $A_L=\sqrt{1-e^{-F}}\,\widetilde K_L\sqrt{1-e^{-F}}$.
The intensity convergence and \eqref{CG:eq:tilted-square-trace} yield
\[
\log\det(\Id-A_L)
=-\Tr A_L+o(1)
\longrightarrow
-\sum_{j=1}^r\int_E(1-e^{-f_j})\dd u\dd v.
\]
The determinantal Laplace-functional identity therefore proves the joint convergence to independent Poisson point processes.

Let $\mathcal N_L(C)$ denote the total count in $C$.
For every $t>0$, \eqref{CG:eq:tilted-square-trace} gives
\[
\Ee_{\Q_L}e^{t\mathcal N_L(C)}
\le\exp\{(e^t-1)\Tr\widetilde K_L\}
=O_{C,t}(1).
\]
Since the rectangles $C_j$ are arbitrary, these bounds give uniform integrability of all fixed mixed factorial counts on relatively compact sets.
Their convergence on continuity sets follows from the joint weak convergence.
\end{proof}

The fixed-order comparison below is also established in \cite[Proposition~4.3]{AiPalm2026}; we give a proof with explicit kernel and spectral hypotheses for its use here.

\begin{proposition}\label{CG:prop:canonical}
Let $X$ be a locally compact second-countable Hausdorff space with a nonatomic Radon measure $\mu$, and let $T$ be a positive trace-class contraction on $L^2(X,\mu)$ with continuous Hermitian correlation kernel.
Assume that an orthonormal eigenbasis $(\phi_j)$, with eigenvalues $(\lambda_j)$, satisfies
$T(x,y)=\sum_j\lambda_j\phi_j(x)\overline{\phi_j(y)}$
for all $x,y$ in a common conull set, including on the diagonal.
Let $\Pp_T$ be the determinantal law with kernel $T$, let $N$ be its total count, and set $m=\Tr T$ and $v=\Tr(T-T^2)$.
For integers $n$ with $\Pp_T\{N=n\}>0$, let $\rho_{T,n}^{(k)}$ denote the $k$th factorial-moment density under $\Pp_T(\,\cdot\mid N=n)$.
Fix $k\ge1$ and $C>0$, and write $T[\boldsymbol p]=[T(p_i,p_j)]_{i,j=1}^k$.
As $v\to\infty$,
\[
\operatorname*{ess\,sup}_{\boldsymbol p:\,\det T[\boldsymbol p]>0}
\left|
\frac{\rho_{T,n}^{(k)}(\boldsymbol p)}
     {\det T[\boldsymbol p]}-1
\right|
=o(1),
\]
uniformly over $T$ and integers $n$ satisfying $|n-m|\le C$ and $\Pp_T\{N=n\}>0$.
Moreover, $\rho_{T,n}^{(k)}$ vanishes almost everywhere on $\{\det T[\boldsymbol p]=0\}$.
\end{proposition}

\begin{proof}
For almost every $\boldsymbol p$ with $M=T[\boldsymbol p]$ positive definite, the reduced Campbell formula and \cite[Theorem~6.5 and Corollary~6.6]{ShiraiTakahashi2003} give
\[
\rho_{T,n}^{(k)}(\boldsymbol p)
=\det M\,
\frac{\Pp_{T'}\{N=n-k\}}{\Pp_T\{N=n\}},
\]
where $T'=T-R_{\boldsymbol p}$ is the reduced-Palm kernel and
$R_{\boldsymbol p}(x,y)=T(x,\boldsymbol p)M^{-1}T(\boldsymbol p,y)$.
Define $A:\ell^2\to L^2(X,\mu)$ by $Ae_j=\sqrt{\lambda_j}\phi_j$ and $W:\Cc^k\to\ell^2$ by $We_i=(\sqrt{\lambda_j}\,\overline{\phi_j(p_i)})_j$.
The diagonal spectral expansion makes $W$ well defined, and $AA^*=T$ and $W^*W=M$.
Thus $R_{\boldsymbol p}=A\Pi A^*$, where $\Pi=WM^{-1}W^*$ is an orthogonal projection of rank $k$.
Consequently, $0\le R_{\boldsymbol p}\le T\le\Id$ and $\Tr R_{\boldsymbol p}\le k$.
Writing $m',v'$ for the mean and variance under $\Pp_{T'}$, we obtain $0\le m-m'=\Tr R_{\boldsymbol p}\le k$.
Since $T$ and $T'$ are positive contractions,
\[
|v'-v|
=\left|\Tr\!\left(R_{\boldsymbol p}(T+T'-\Id)\right)\right|
\le\Tr R_{\boldsymbol p}\le k.
\]
In particular, $v'\to\infty$ and $|(n-k)-m'|\le C+k$.

Both counts admit Bernoulli spectral representations.
Removing the deterministic summands with parameter $1$ preserves the variances and centering errors.
The local central limit theorem \cite[Theorem~1]{ForresterLebowitz2014}, whose proof depends only on the Bernoulli spectra, therefore gives
\[
\Pp_T\{N=n\}
=\frac{1+o(1)}{\sqrt{2\pi v}},
\qquad
\Pp_{T'}\{N=n-k\}
=\frac{1+o(1)}{\sqrt{2\pi v'}}.
\]
The errors are uniform in $T,n,\boldsymbol p$: otherwise, a violating sequence would have divergent variances and bounded centering errors, contradicting the cited theorem.
Since $v'/v\to1$ uniformly, the probability ratio tends to one uniformly, proving the comparison.
Finally,
$\Pp_T\{N=n\}\rho_{T,n}^{(k)}(\boldsymbol p)\le\det T[\boldsymbol p]$
almost everywhere, which proves the assertion on the zero set.
\end{proof}

\begin{proof}[Proof of Theorem~\ref{thm:intro-conditioned-poisson}]
The density in \eqref{CG:eq:tilted-law} is constant on $\{N_{a,L}=n_L\}$, so conditioning under $\Pp$ or $\Q_L$ gives the same law.
To transfer the tilted limits, we first verify the pointwise spectral hypothesis of Proposition~\ref{CG:prop:canonical}.
Set $s=\widehat s_L$ and let $J_s$ be the map $f\mapsto\sqrt{s}\,f|_{D_{a,L}}$ from the $w_{L,s}$-weighted Bergman space to $L^2(D_{a,L},\ddA)$.
Since $B_{L,s}(z,z)\le s^{-1}B_1(z,z)$,
\[
\|J_s\|_2^2=s\int_{D_{a,L}}B_{L,s}(z,z)\ddA(z)
\le\frac{L}{4\pi a}.
\]
The map $J_s$ is injective by the identity theorem.
As $K_{a,L}^{(s)}=J_sJ_s^*$, the singular-value decomposition of $J_s$ and Parseval for point evaluations give the required spectral expansion, including the diagonal.
The tilted count has mean $n_L$ and variance $v_L(s)\sim n_L\to\infty$, so Proposition~\ref{CG:prop:canonical} applies uniformly in the insertion points at each fixed order.

For the global limit, fix bounded continuous real-valued $f$ and set $X_{L,f}=n_L\langle\mathfrak M_L,f\rangle$.
The uniform comparisons at orders $k=1,2$ give
\[
\Ee[X_{L,f}\mid N_{a,L}=n_L]
=\Ee_{\Q_L}X_{L,f}+o(n_L),\qquad
\Ee[X_{L,f}^2\mid N_{a,L}=n_L]
=\Ee_{\Q_L}X_{L,f}^2+o(n_L^2).
\]
Here $k=2$ treats distinct pairs in the second moment, while $k=1$ treats its diagonal contribution; the corresponding tilted factorial moments are $O(n_L^2)$ and $O(n_L)$.
Proposition~\ref{CG:prop:tilted-Poisson} therefore yields
\[
\Ee\left[\left|
\langle\mathfrak M_L,f\rangle
-\int_0^1\!\int_0^1 f(u,v)\dd u\dd v
\right|^2\,\middle|\,N_{a,L}=n_L\right]\longrightarrow0.
\]
Applying this conclusion to a countable convergence-determining family and to continuous height cutoffs gives the asserted weak convergence in probability of the conditional empirical measures.

For the joint local limit, integrate the uniform comparison in Proposition~\ref{CG:prop:canonical} over products of pairwise disjoint relatively compact continuity sets in the disjoint union of the local observation spaces.
At every fixed total order, the conditional mixed factorial moments differ by $o(1)$ from their tilted counterparts, which converge to those of independent Poisson counts by Proposition~\ref{CG:prop:tilted-Poisson}.
The Poisson count vectors are determined by these moments, so the conditional count vectors converge in distribution.
The order-one comparison also bounds expected counts uniformly on compact sets, giving tightness of the local point processes in the vague topology.
These count limits identify the joint limit as independent Poisson point processes.
Overlapping continuity sets are handled by partitioning them into finitely many disjoint continuity cells, which also establishes the stated convergence of mixed factorial moments.
This proves \eqref{eq:intro-conditioned-poisson}.
\end{proof}

At the exact saddle, the tilted count has mean $n_L$ and variance $v_L(\widehat s_L)\sim n_L\to\infty$.
Since $K_{a,L}^{(\widehat s_L)}$ has spectrum in $[0,1)$, the determinantal local central limit theorem \cite[Theorem~1]{ForresterLebowitz2014} and the change of measure \eqref{CG:eq:tilted-law} give
\begin{equation}\label{CG:eq:saddle-probability}
\Q_L\{N_{a,L}=n_L\}\sim\frac1{\sqrt{2\pi n_L}},\qquad
\Pp\{N_{a,L}=n_L\}
=\frac{\Ee\widehat s_L^{N_{a,L}}\,\widehat s_L^{-n_L}}
{\sqrt{2\pi n_L}}(1+o(1)).
\end{equation}
Section~\ref{sec:endpoint} uses this estimate to derive relative count asymptotics throughout the diverging sublinear count regime.

\section{Voids and dilute relative count asymptotics}\label{sec:endpoint}

Fix $a>0$.
We determine the void exponent and control its remainder through a finite-window mean deficit.
This comparison also yields relative logarithmic asymptotics for diverging sublinear counts.
Spectral estimates then give a sharper description of small-count probabilities.
Let $K_{a,\Rr}$ be the compression of the Bergman projection with kernel $\mathcal K_{\Hh}$ to $L^2(\Rr\times[a,\infty),\ddA)$; its compression to $D_{a,L}=[0,L]\times[a,\infty)$ is $K_{a,L}$.

\subsection{Void probability and mean deficit}

The following trace bound will be used for both the void probability and the small-count estimates.

\begin{lemma}\label{VC:lem:moment-bound}
For every integer $m\ge1$,
\[
0\le\Tr K_{a,L}^m\le\frac{L}{4\pi am}.
\]
\end{lemma}

\begin{proof}
Let $(\lambda_j,f_j)$ be the eigenpairs of $K_{a,L}$, with $f_j$ extended by zero to $\Rr\times[a,\infty)$.
Since $K_{a,L}$ is a compression of $0\le K_{a,\Rr}\le\Id$, Jensen's inequality gives $\lambda_j^m\le\langle f_j,K_{a,\Rr}^mf_j\rangle$.
The horizontal Fourier fibres of $K_{a,\Rr}$ have rank one for positive frequencies, with nonzero eigenvalue $e^{-2a\xi}$.
Summing over $j$ and using positivity yields
\[
\Tr K_{a,L}^m
\le\Tr(\mathbf1_{D_{a,L}}K_{a,\Rr}^m\mathbf1_{D_{a,L}})
=\frac{L}{2\pi}\int_0^\infty e^{-2am\xi}\dd\xi
=\frac{L}{4\pi am}.
\]
\end{proof}

For $z\in\Cc$, set $G_{a,L}(z):=\Ee z^{N_{a,L}}$; in particular, $G_{a,L}(0)=\Pp\{N_{a,L}=0\}$.
With $x_a$ and $\mu_L$ as in Section~\ref{CG:sec:saddle}, define $\Delta_{a,L}(s):=Lx_a(s)-\mu_L(s)$ for $0<s\le1$, and set
\begin{equation}\label{VC:eq:B-def}
\mathcal B_{a,L}(s)
:=L\left[\Lambda_a^\infty(\log s)+\frac{\pi}{24a}\right]
=\frac{L}{4\pi a}\left[\Li_2(s)+\log s\log(1-s)\right],
\qquad 0<s<1.
\end{equation}
The function $\mathcal B_{a,L}$ extends continuously to $[0,1]$, with $\mathcal B_{a,L}(0)=0$ and $\mathcal B_{a,L}(1)=\pi L/(24a)$.
Define the void remainder by
\begin{equation}\label{VC:eq:R-def}
\mathcal R_a(L):=\log G_{a,L}(0)+\frac{\pi L}{24a}.
\end{equation}

\begin{lemma}\label{VC:lem:mean-deficit}
For $0<s\le1$,
\begin{equation}\label{VC:eq:mean-deficit-bridge}
\mathcal B_{a,L}(s)-\log\frac{G_{a,L}(s)}{G_{a,L}(0)}
=\int_0^s\frac{\Delta_{a,L}(u)}u\dd u\ge0.
\end{equation}
In particular,
\begin{equation}\label{VC:eq:R-deficit}
\mathcal R_a(L)=\int_0^1\frac{\Delta_{a,L}(u)}u\dd u.
\end{equation}
\end{lemma}

\begin{proof}
By Lemma~\ref{CG:lem:spectral-basics}, $G_{a,L}(0)>0$.
For $0<u<1$, the Fredholm product and \eqref{VC:eq:B-def} give $u(\log G_{a,L})'(u)=\mu_L(u)$ and $u\mathcal B_{a,L}'(u)=Lx_a(u)$.
Integrating their difference and using continuity at zero proves \eqref{VC:eq:mean-deficit-bridge} for $0<s<1$.
The integrand is nonnegative by \eqref{CG:eq:mean-upper}.
Continuity extends the identity to $s=1$, where $G_{a,L}(1)=1$ and $\mathcal B_{a,L}(1)=\pi L/(24a)$, giving \eqref{VC:eq:R-deficit}.
\end{proof}

\begin{theorem}[Void exponent and remainder]\label{VC:thm:void}
For fixed $a>0$, as $L\to\infty$,
\begin{equation}\label{VC:eq:void-rate}
\log\Pp\{N_{a,L}=0\}=-\frac{\pi}{24a}L+o(L).
\end{equation}
Moreover, the remainder in \eqref{VC:eq:R-def} satisfies
\[
\mathcal R_a(L)\ge0,\qquad \mathcal R_a(L)=o(L),\qquad
\frac{\mathcal R_a(L)}{\log(L/a)}\longrightarrow\infty.
\]
\end{theorem}

\begin{proof}
By Lemma~\ref{CG:lem:spectral-basics}, $\|K_{a,L}\|<1$, so Lemma~\ref{VC:lem:moment-bound} gives
\[
-\log G_{a,L}(0)
=\sum_{m\ge1}\frac{\Tr K_{a,L}^m}{m}
\le\frac{L}{4\pi a}\sum_{m\ge1}\frac1{m^2}
=\frac{\pi L}{24a}.
\]
For each fixed $t>0$, $G_{a,L}(0)\le G_{a,L}(e^{-t})$, so Theorem~\ref{thm:thermo-halfline} and exact scaling give $\limsup_{L\to\infty}L^{-1}\log G_{a,L}(0)\le\Lambda_a^\infty(-t)$.
Since $\Lambda_a^\infty(-t)\to-\pi/(24a)$ as $t\to\infty$, the two bounds prove \eqref{VC:eq:void-rate}, $\mathcal R_a(L)\ge0$, and $\mathcal R_a(L)=o(L)$.

For any fixed $T>0$, \eqref{VC:eq:mean-deficit-bridge} and the strong half-line expansion at $\theta=-T$ yield
\[
\mathcal R_a(L)
\ge\log G_{a,L}(e^{-T})-L\Lambda_a^\infty(-T)
=\frac{T^2}{2\pi^2}\log(L/a)+A_\infty(-T)+o_T(1).
\]
Consequently,
$\liminf_{L\to\infty}\mathcal R_a(L)/\log(L/a)
\ge T^2/(2\pi^2)$.
Since $T>0$ is arbitrary, the final assertion follows.
\end{proof}

\subsection{Probabilities of diverging dilute counts}

The mean-deficit identity also yields the following comparison as $s\downarrow0$.

\begin{proposition}\label{VC:prop:thermo-bridge}
Let $s_L\in(0,1)$ satisfy $s_L\to0$ and $Lx_a(s_L)\to\infty$.
Then
\[
\log\frac{G_{a,L}(s_L)}{G_{a,L}(0)}
\sim\mathcal B_{a,L}(s_L)\sim Lx_a(s_L).
\]
\end{proposition}

\begin{proof}
Fix $0<\varepsilon<1$.
Since $x_a(\varepsilon s_L)\sim\varepsilon x_a(s_L)$ and $x_a$ is increasing, Proposition~\ref{CG:prop:mean-variance} gives $\mu_L(u)=Lx_a(u)(1+o(1))$ uniformly for $\varepsilon s_L\le u\le s_L$.
Indeed, every sequence $u_L$ in this interval satisfies $u_L\to0$ and $Lx_a(u_L)\ge Lx_a(\varepsilon s_L)\to\infty$.

To compare the generating functions, split \eqref{VC:eq:mean-deficit-bridge} at $\varepsilon s_L$.
Using $0\le\Delta_{a,L}(u)\le Lx_a(u)$ yields
\[
0\le\mathcal B_{a,L}(s_L)
-\log\frac{G_{a,L}(s_L)}{G_{a,L}(0)}
\le\mathcal B_{a,L}(\varepsilon s_L)
+o(1)\bigl[\mathcal B_{a,L}(s_L)
-\mathcal B_{a,L}(\varepsilon s_L)\bigr].
\]
By \eqref{VC:eq:B-def}, $\mathcal B_{a,L}(s)\sim Lx_a(s)$ as $s\downarrow0$, so
$\mathcal B_{a,L}(\varepsilon s_L)/
\mathcal B_{a,L}(s_L)\to\varepsilon$.
Dividing by $\mathcal B_{a,L}(s_L)$ and letting $L\to\infty$, followed by $\varepsilon\downarrow0$, proves the assertion.
\end{proof}

\begin{proof}[Proof of Theorem~\ref{thm:intro-dilute-relative}]
Let $s_L$ solve $Lx_a(s_L)=n_L$, and let $\widehat s_L$ be the exact saddle.
Proposition~\ref{CG:prop:saddle-asymptotics} gives $\widehat s_L/s_L\to1$ and $Lx_a(\widehat s_L)\sim n_L$.
Proposition~\ref{VC:prop:thermo-bridge} and \eqref{CG:eq:saddle-probability} yield
\[
\log\frac{\Pp\{N_{a,L}=n_L\}}{\Pp\{N_{a,L}=0\}}
=\mathcal B_{a,L}(\widehat s_L)
-n_L\log\widehat s_L+o(n_L),
\]
since $\log n_L=o(n_L)$.

To replace the exact saddle by $s_L$, note that $\mu_L(s_L)\le Lx_a(s_L)=n_L=\mu_L(\widehat s_L)$.
The monotonicity of $\mu_L$ implies $s_L\le\widehat s_L$.
Since $u\mathcal B_{a,L}'(u)=Lx_a(u)$ and $x_a$ is increasing,
\[
0\le
\mathcal B_{a,L}(\widehat s_L)-n_L\log\widehat s_L
-\mathcal B_{a,L}(s_L)+n_L\log s_L
\le
\bigl[Lx_a(\widehat s_L)-n_L\bigr]
\log(\widehat s_L/s_L)=o(n_L).
\]
The saddle equation $(\Lambda_a^\infty)'(\log s_L)=n_L/L$ and Legendre duality give
\[
\mathcal B_{a,L}(s_L)-n_L\log s_L
=L\bigl[J_a^\infty(0)-J_a^\infty(n_L/L)\bigr].
\]
This proves \eqref{eq:intro-dilute-relative}.
The explicit expansion \eqref{eq:intro-dilute-relative-explicit} follows from Corollary~\ref{cor:low-density-rate}.
\end{proof}

\subsection{Spectral estimates for small counts}

The strict contraction property of $K_{a,L}$ allows us to define
\begin{equation}\label{VC:eq:W-def}
\mathsf V_{a,L}:=K_{a,L}(\Id-K_{a,L})^{-1},
\qquad W_a(L):=\Tr\mathsf V_{a,L}.
\end{equation}
The Fredholm product gives the exact factorization
\begin{equation}\label{VC:eq:pgf-factor}
G_{a,L}(z)=G_{a,L}(0)\det(\Id+z\mathsf V_{a,L}).
\end{equation}

\begin{lemma}\label{VC:lem:W-bounds}
For fixed $a>0$, as $L\to\infty$,
\[
W_a(L)\sim\frac{L}{4\pi a}\log(L/a).
\]
\end{lemma}

\begin{proof}
For an upper bound, put $u=L/(4\pi a)$, $m_0=\lceil4u\rceil$, and $\lambda_*=\|K_{a,L}\|$.
Lemma~\ref{VC:lem:moment-bound} gives $\lambda_*^{m_0}\le\Tr K_{a,L}^{m_0}\le u/m_0\le1/4$.
Thus $(1-\lambda_*)^{-1}\le(1-4^{-1/m_0})^{-1}=O(L/a)$.
The trace-class Neumann series, split at $m_0$, then gives
\[
W_a(L)=\sum_{m\ge1}\Tr K_{a,L}^m
\le u\sum_{m=1}^{m_0-1}\frac1m
+\frac{\Tr K_{a,L}^{m_0}}{1-\lambda_*}
\le\frac{L}{4\pi a}\log(L/a)+O(L/a).
\]

For the matching lower bound, take $s=a/L$.
Then $s\to0$ and $Lx_a(s)=\log(L/a)/[4\pi(1-a/L)]\to\infty$.
Proposition~\ref{CG:prop:mean-variance} therefore gives $\mu_L(s)\sim Lx_a(s)$.
Since $K_{a,L}^{(s)}=s\mathsf V_{a,L}(\Id+s\mathsf V_{a,L})^{-1}\le s\mathsf V_{a,L}$,
\[
W_a(L)\ge\frac{\mu_L(a/L)}{a/L}
\sim\frac{L}{4\pi a}\log(L/a).
\]
The two bounds prove the asserted asymptotic.
\end{proof}

The Fourier--Laplace compression gives stronger spectral estimates at the zero-count endpoint.

\begin{lemma}\label{VC:lem:refined-spectral}
For fixed $a>0$, put $T=L/a$.
As $T\to\infty$,
\[
 T\bigl(1-\|K_{a,L}\|\bigr)\longrightarrow\infty,
 \qquad
 \|\mathsf V_{a,L}\|=o(T),
 \qquad
 \Tr\mathsf V_{a,L}^2=o(T^2).
\]
\end{lemma}

\begin{proof}
Put $h=T^{-1}$ and recall $\mathsf P_h$ and $\mathsf U_-$ from Lemma~\ref{MC:lem:compression}.
The nonzero spectra of $K_{a,L}$ and $\mathsf P_h$ agree, and the compression representation gives $0\le\mathsf P_h\le\mathsf P_0\le\Id$.
Allowing all horizontal coordinates in that representation also gives $\mathsf P_h\le M_h^2$, where $(M_hf)(\xi)=e^{-h\xi}f(\xi)$.
For every $R>0$, $\mathbf1_{(0,R)}\mathsf P_0$ is Hilbert--Schmidt: its kernel is bounded for $0<\xi<R$, $0<\eta\le2R+2$, and is $O_R(\sqrt\xi\,\eta^{-3/2})$ otherwise, by \eqref{MC:eq:Ph-kernel} at $h=0$.
Moreover, $\mathsf P_0$ has no eigenvalue $1$.
Indeed, equality in the compression would force the anti-holomorphic image $\mathsf U_-f$ of a corresponding eigenvector to vanish outside $[0,1]\times(0,\infty)$, so the identity theorem and the isometry give $f=0$.

If the first assertion failed, there would be $h_n\downarrow0$ with $\|\mathsf P_{h_n}\|\ge1-Ch_n$.
Choose normalized top eigenvectors $f_n$ of these positive trace-class operators.
The two operator bounds give
\[
 \int_0^\infty(1-e^{-2h_n\xi})|f_n(\xi)|^2\dd\xi\le Ch_n,
 \qquad
 \|f_n-\mathsf P_0f_n\|_2^2\le Ch_n.
\]
Thus $\limsup_{n\to\infty}\|\mathbf1_{(R,\infty)}f_n\|_2^2\le C/(2R)$ for every $R>0$.
Local compactness of $\mathbf1_{(0,R)}\mathsf P_0$ and $f_n-\mathsf P_0f_n\to0$ give convergent subsequences on bounded frequency intervals.
A diagonal extraction and the tail bound yield a subsequence converging strongly to $f$ with $\|f\|_2=1$ and $\mathsf P_0f=f$, a contradiction.
Hence $T(1-\|K_{a,L}\|)\to\infty$, and the asserted norm bound follows from $\|\mathsf V_{a,L}\|=\lambda_* /(1-\lambda_*)$, where $\lambda_*=\|K_{a,L}\|$.

Write $K=K_{a,L}$, $V=\mathsf V_{a,L}$ and $u=T/(4\pi)$.
The positive series gives $\Tr V^2=\sum_{m\ge2}(m-1)\Tr K^m$.
For $r=\lfloor m/2\rfloor$, Lemma~\ref{VC:lem:moment-bound} yields
\[
 (m-1)\Tr K^m
 \le (m-1)\lambda_*^{m-r}\frac{u}{r}
 \le2u\lambda_*^{m-r}.
\]
Summation gives $\Tr V^2\le2u\lambda_*(1+\lambda_*)/(1-\lambda_*)=o(T^2)$, proving the last assertion.
\end{proof}

The preceding spectral bounds yield a sharper estimate for diverging counts up to logarithmic order.

\begin{corollary}\label{VC:cor:log-count-scale}
Fix $a>0$, and let $n_L$ be positive integers with $n_L\to\infty$ and $n_L=O(\log(L/a))$.
Let $\widehat s_L$ be the exact saddle satisfying $\mu_L(\widehat s_L)=n_L$.
Then
\[
\log\frac{\Pp\{N_{a,L}=n_L\}}{\Pp\{N_{a,L}=0\}}
=-n_L\log\widehat s_L+n_L
-\frac12\log(2\pi n_L)+o(1).
\]
\end{corollary}

\begin{proof}
Put $T=L/a$, $s=\widehat s_L$ and $V=\mathsf V_{a,L}$.
Proposition~\ref{CG:prop:saddle-asymptotics} gives $s=O(T^{-1})$.
By Lemma~\ref{VC:lem:refined-spectral}, $s\|V\|=o(1)$ and $s^2\Tr V^2=o(1)$.
Thus the eigenvalues $r_j$ of $K_{a,L}^{(s)}=sV(\Id+sV)^{-1}$ satisfy $\sup_jr_j\to0$, $\sum_jr_j=n_L$ and $\sum_jr_j^2\le s^2\Tr V^2=o(1)$.
The factorization \eqref{VC:eq:pgf-factor} gives
\[
0\le\log\frac{G_{a,L}(s)}{G_{a,L}(0)}-n_L
=\sum_j[-\log(1-r_j)-r_j]
\le C\sum_jr_j^2=o(1).
\]
The result follows from \eqref{CG:eq:saddle-probability}.
\end{proof}

The same estimates give relative probability equivalents uniformly up to logarithmic counts.

\begin{theorem}[Small-count asymptotics]
\label{VC:thm:near-zero}
Fix $a>0$, and put $T=L/a$.
For every $R>0$, as $L\to\infty$,
\[
 \log\frac{G_{a,L}(\tau/W_a(L))}{G_{a,L}(0)}
 =\tau+o_R((\log T)^{-2})
\]
uniformly for $|\tau|\le R$, with the logarithm normalized at $\tau=0$.
For every fixed $M>0$,
\[
\Pp\{N_{a,L}=k\}
=\Pp\{N_{a,L}=0\}\frac{W_a(L)^k}{k!}(1+o(1))
\]
uniformly for integers $0\le k\le M\log T$ as $L\to\infty$.
More generally, the equivalent holds for every sequence of nonnegative integers $k=k_L$ satisfying $k_L^2\Tr\mathsf V_{a,L}^2/W_a(L)^2\to0$.
\end{theorem}

\begin{proof}
Write $V=\mathsf V_{a,L}$ and $W=W_a(L)$.
Lemmas~\ref{VC:lem:W-bounds} and~\ref{VC:lem:refined-spectral} give
\[
 \frac{\|V\|}{W}=o((\log T)^{-1}),
 \qquad
 \frac{\Tr V^2}{W^2}=o((\log T)^{-2}).
\]
For $|\tau|\le R$, we therefore have $|\tau|\|V\|/W<1$ for all sufficiently large $L$.
The normalized logarithm is then well defined, and \eqref{VC:eq:pgf-factor} and the trace-logarithm expansion yield
\[
 \log\frac{G_{a,L}(\tau/W)}{G_{a,L}(0)}
 =\tau+O_R\!\left(\frac{\Tr V^2}{W^2}\right)
 =\tau+o_R((\log T)^{-2}).
\]

For the count probabilities, let $(\mu_j)$ be the eigenvalues of $V$ and put $p_j=\mu_j/W$.
By \eqref{VC:eq:pgf-factor}, $\Pp\{N_{a,L}=k\}/G_{a,L}(0)=e_k((\mu_j))$, the $k$th elementary symmetric function.
Hence $k!e_k((\mu_j))/W^k$ is the probability that $k$ independent draws from $(p_j)$ have distinct indices.
The union bound gives
\[
0\le
1-\frac{k!\,\Pp\{N_{a,L}=k\}}
{G_{a,L}(0)W^k}
\le\binom{k}{2}\sum_jp_j^2
=\binom{k}{2}\frac{\Tr V^2}{W^2}.
\]
The right-hand side tends to zero uniformly for $0\le k\le M\log T$, and for every sequence satisfying the stated criterion.
\end{proof}

For fixed $k$, Lemma~\ref{VC:lem:W-bounds} also gives the explicit equivalent
\[
\Pp\{N_{a,L}=k\}
=\Pp\{N_{a,L}=0\}\frac1{k!}
\left(\frac{L}{4\pi a}\log(L/a)\right)^k(1+o(1)).
\]
For growing $k=k_L$, replacing $W_a(L)$ by its leading equivalent requires
\[
k_L\log\frac{4\pi aW_a(L)}{L\log(L/a)}\longrightarrow0.
\]
This condition does not follow from the equivalent for $W_a(L)$ alone.

\end{document}